\documentclass[11pt]{amsart}

\usepackage{hyperref,wasysym,cite}
\hypersetup{colorlinks=true, citecolor=BurntOrange, linkcolor=ForestGreen}

\usepackage{amssymb}
\usepackage{graphicx}
\usepackage{amsmath}
\usepackage[usenames,dvipsnames]{xcolor}
\usepackage[text={6.5in,9in},centering]{geometry}

\newcommand{\eq}[2]{\begin{equation}\begin{split}#1\end{split}\label{#2}\end{equation}}
\newcommand{\eqnn}[1]{\begin{equation}\begin{split}#1\end{split}\nonumber\end{equation}}

\newcommand{\red}[1]{{\color{red}{#1}}}

\newcommand{\intr}{\int_{\mathbb{R}}}

\newcommand{\LT}{L^2({\mathbb{R})}}

\newcommand{\realpart}[1]{\operatorname{\rm Re}\!\left(#1\right)}
\newcommand{\impart}[1]{\operatorname{\rm Im}\!\left(#1\right)}
\newcommand{\argum}[1]{\operatorname{\rm arg}\left(#1\right)}
\newcommand{\argumc}[1]{\operatorname{\rm arg}\!\left(#1\right)}

\newtheorem{Lemma}{Lemma}[section]

\newtheorem{Definition}[Lemma]{Definition}

\newtheorem{Proposition}[Lemma]{Proposition}
\newtheorem{Remark}[Lemma]{Remark}
\newtheorem{Theorem}[Lemma]{Theorem}

\newenvironment{Proof}[1][.]%
  {\begin{trivlist}\item[]\textbf{Proof#1 }}%
  {\hspace*{\fill}$\rule{.3\baselineskip}{.35\baselineskip}$\end{trivlist}}

\def\Re{\mathop\mathrm{Re}\nolimits}      

\newcommand{\C}{\mathbb{C}}               
\newcommand{\R}{\mathbb{R}}               

\newcommand{\eps}{\epsilon} 
\newcommand{\cK}{\mathcal{K}}

\begin{document}

\numberwithin{equation}{section}

\title{ Stability of fronts  in the diffusive Rosenzweig-MacArthur  model.}


 
\date{\today}

\begin{abstract}{We consider a diffusive Rosenzweig-MacArthur predator-prey model in the situation when the prey diffuses at a rate much smaller than that of the predator. In a certain parameter regime, the existence of fronts in the system is known:  the underlying dynamical system in a singular limit is reduced to a scalar Fisher-KPP equation and the fronts supported by the full system are small perturbations of the Fisher-KPP fronts. The existence proof is based on the application of the Geometric Singular Perturbation Theory with respect to two small parameters. This paper is focused on the stability of the fronts. We {{show that, for some parameter regime, the fronts are spectrally and asymptotically stable using}}  energy estimates, exponential dichotomies,  the Evans function calculation, and a technique that involves constructing unstable augmented bundles.  The energy estimates provide bounds on the unstable spectrum which depend on the small parameters of the system; the bounds are inversely proportional to these parameters. We further improve these estimates by showing that the eigenvalue problem is a small perturbation of some limiting (as the modulus of the eigenvalue parameter goes to infinity) system and that the limiting system has exponential dichotomies. Persistence of the exponential dichotomies then leads to bounds uniform in the small parameters. The main novelty of this approach is related to the fact that the limit of the eigenvalue problem is not autonomous.  We then use the concept of the unstable augmented bundles and by treating these as multiscale topological structures with respect to the same two small parameters consequently as in the existence proof,  we show that the stability of the fronts is also governed by the scalar Fisher-KPP equation. 
{{Furthermore, we perform numerical computations of the Evans function to {{explicitly identify regions in the parameter space where the fronts are spectrally stable.}}}}}

\end{abstract}
\maketitle

\begin{center}    Anna Ghazaryan \textsuperscript{a},   St\'ephane Lafortune \textsuperscript{b}, Yuri Latushkin \textsuperscript{c}, Vahagn Manukian \textsuperscript{a,d}
 \end{center}
 
 \textsuperscript{a} \address{Department of Mathematics, Miami University,  Oxford, OH 45056, USA,}  
 \email{ghazarar@miamioh.edu.} \par
  \textsuperscript{b} \address{Department of Mathematics, College of Charleston, Charleston, SC 29401, USA,
  }   \email{LafortuneS@cofc.edu.}\par
   \textsuperscript{c} \address{Department of Mathematics, The University of Missouri, Columbia, MO 65211, USA,}  
 \email{latushkiny@missouri.edu}\par
 \textsuperscript{d} \address{Department of Mathematical and Physical Sciences, Miami University, Hamilton, OH 45011, USA,}  \email{manukive@miamioh.edu.}\par 


 
\section{Rosenzweig-MacArthur model  \label{s:intro}}


In 1963 Rosenzweig and MacArthur made an observation that influenced mathematical modeling of many predator-prey systems. They suggested that the rate at which predators consume prey stays bounded regardless of the density of the prey population. The situations where  the predator's population cannot grow without bounds, but instead   has a ``ceiling''  may be described using   the Rosenzweig-MacArthur  model   \cite{RM} 
\begin{eqnarray}\label{e:ODE0}
U_{t}&=&\mathcal AU\left(1-\frac{U}{\mathcal K}\right)-\frac{\mathcal B UW}{1+\mathcal E U},\notag\\
W_{t}&=&- CW+ \frac{\mathcal D UW}{1+\mathcal E  U}.
\end{eqnarray} 
 Here $t$ is the time variable,  $U$ and $W$ are the population densities of the prey and predator, respectively. 
 Parameter $\mathcal A>0$ is the linear growth factor for the prey species, $\mathcal K>0$ is the carrying capacity of the prey species. The death rate for the predator without prey is $\mathcal C>0$;   $\mathcal B>0$ and $\mathcal D >0$ are the interaction rates for the two species, ($B$ is is the rate of the nonlinear loss for the prey and $D$ is the rate for the nonlinear gain for the predator). The MacArthur-Rosenzweig observation is taken into the account through  the expression $\frac{U}{1+ \mathcal E U}$, where $\mathcal E>0$ is a constant.
 
In predator-prey systems,   the Rosenzweig-MacArthur  model may be used to capture a  variety of specific situations. It is typically used to reflect  predator's satiation. Another  interpretation for it is that the predators need time to handle the prey thus  causing saturation in consumption. Another mechanism captured by the Rosenzweig-MacArthur model is  related to the fact that higher predation rates  are often associated with a decrease in the prey mortality. Yet another example   is  related to   the ability of the prey to  take an environmental refuge which reflects negatively on the predation rates. 

 An interesting extension of the Rosenzweig-MacArthur  model is  formulated in  \cite{AVOOMR} where a model is studied that replaces the classical Rosenzweig-MacArthur term with a ``refuge function''  $\frac{ (U-U_r) }{\mathcal K+(U-U_r)}$  \cite{AVGOGY}. Indeed when  a quantity $U_r$ of prey  takes a refuge, then predator have access only to   $U-U_r$ of the  prey. In this case $\mathcal K$ is the number of prey  that represents the half of the maximum capacity of the refuge.  In \cite{AVOOMR}  the  stability of  the physically meaningful equilibrium points is studied and  the existence of limit cycles is shown. 
Obviously,  when $U_r$ is a proportion of the prey population then the refuge function is equivalent to the MacArthur term $\frac{D U}{1+ \mathcal E U}$.
 
One could consider  a system where the environment provides  not only  protection (advantage)  to the prey population, but also provides advantage or disadvantage  to the predator. This argument is brought up in a different model suggested in  \cite{AADO} where the low density prey reduces the linear growth of predator  by  a nonlinear factor  $\frac{ W}{1+ \mathcal E_2 U}$,  where $\mathcal E_2$ is different than the constant in the nonlinear reduction  $\frac{W}{1+\mathcal E_1 U}$ term  of the linear growth rate  of prey population. This happens,  for example, in pest-spider predation  \cite{AADO}, where
 the boundedness  of the solutions  is proved and  the co-existence equilibrium is shown to be globally stable    \cite{AADO}.  We also refer readers to \cite{AW1,SGO} for more examples of population systems that may indicate $\mathcal E_1\neq \mathcal E_2$. 
Incorporating  this idea in \eqref{e:ODE0},  one obtains
\eq{
U_{t}&=\mathcal AU\left(1-\frac{U}{\mathcal K}\right)-\frac{\mathcal B UW}{1+\mathcal E_1U},\\
W_{t}&=-\mathcal CW+ \frac{\mathcal D UW}{1+ \mathcal E_2U},
}{e:010}
where    $\mathcal E_1>0$ and $\mathcal E_2>0$ are not necessarily equal constants  that describe  how much  protection the environment provides to the  prey  and predator, respectively.   Another argument supporting the assumption that $\mathcal E_2> \mathcal E_1$  is that the predator may have another food source, in addition to the prey considered in the model, or kill the prey for reasons other than consumption.  

 A  diffusive version of  the Rosenzweig-MacArthur ODE \eqref{e:ODE0} is given by
\eq{ 
U_{\tau}&=\epsilon_u U_{xx}+\mathcal AU\left(1-\frac{U}{\mathcal K}\right)-\frac{\mathcal B UW}{1+ \mathcal EU},\\
W_{\tau}&=\epsilon_wW_{xx}-\mathcal CW+ \frac{\mathcal D UW}{ 1+\mathcal E U}, 
}{dunbar} 
where $x$ is a one-dimensional spatial variable, has also been studied. For example,   in \cite{SS} periodic wavetrains  in such system were discussed. Fronts, on the other hand, which are traveling waves  that that propagate  with constant velocity without changing their shape and that asymptotically connect distinct equilibria,   were  studied in  the   paper \cite{Dunbar}. More precisely,  the existence of the periodic traveling wavetrains and  fronts connecting an equilibrium to a periodic orbit  was proven in \cite{SS} in the case  when $\epsilon_u=0$.  In  case of  strictly positive  $\epsilon_u$ and $\epsilon_w$,  the existence of traveling waves was demonstrated numerically in \cite{OL}.  The existence of traveling wave solutions and small amplitude traveling wavetrains  was also investigated in  \cite{HLR}.

In  biological applications,  traveling fronts are meaningful  solutions. One of the most important examples  of such solutions are fronts in the  Fisher-KPP equation  \cite{Fisher}, \cite{KPP}, where KPP stands for Kolmogorov-Petrovski-Piskunov.  We  also refer the reader to the  works \cite{Fife},  \cite{M},  and \cite{VVV}. The fronts in the Rosenzweig-MacArthur system \eqref{dunbar} have been observed numerically for different parameter regimes  and their stability was studied   in \cite{DS}. The system  
\eq{
U_{\tau}&=\epsilon_u U_{xx}+\mathcal AU\left(1-\frac{U}{\mathcal K}\right)-\frac{\mathcal B UW}{1+ \mathcal E_1U},\\
W_{\tau}&=\epsilon_wW_{xx}-\mathcal CW+ \frac{\mathcal D UW}{ 1+\mathcal E_2 U}, 
}{CGMeqn}  {{related to \eqref{e:010}}}, has been  studied in \cite{CGM}, in the  parameter regimes  that were not covered in \cite{DS}.  It was proven in \cite{CGM} that in some parameter regimes  fronts exist  which are in some sense related to  a Fisher-KPP equation ``embedded'' in  \eqref{CGMeqn} and in other parameter regimes there are fronts which are not associated with any Fisher-KPP equations.  

In this paper we study the stability of fronts, existence of which  has been shown in \cite{CGM}.
In Section \ref{RM}, we review the results about the existence of the front solutions and introduce the two small parameters $\epsilon$ and $\delta$, with respect to which singular perturbation analyses are being performed. In the subsequent section, we obtain the linearization about the fronts and study its properties. In Section \ref{SFEP}, {{we state the main results about the spectral and asymptotic stability and}} we consider the eigenvalue problem associated with the linearization. We  use the two small parameters $\epsilon$ and $\delta$ coming from the existence results, one after another,  and we obtain  ``slow'' and ``fast'' in $\epsilon$ and  then in $\delta$ versions of the eigenvalue problem. In Section \ref{EssSect}, we calculate the right-most boundary  of  the essential spectrum of the linearization and show that, while it intersects with the right side of the complex plane, one can introduce exponential weights that ``move'' the spectrum to the left, i.e., stabilize the essential spectrum.  In the following section, we study the point spectrum by computing a bound on any unstable eigenvalue. Note that both the essential spectrum and the bounds are obtained in the linearization of the  systems with a  positive  $\epsilon$ and in the case of $\epsilon=0$, which partly requires a different approach, {and also in the case  when $\delta$ is set to be zero  in the system first obtained by taking  limit  as $\epsilon=0$, where the system is  reduced to the Fisher-KPP equation through singular perturbation theory.}
 It is proved  in Section \ref{uniform} that there exist  bounds on the point spectrum that are uniform with respect to $\epsilon$ and $\delta$. The bounds obtained in Section \ref{ESSECT} are used later in the numerical computations, while the uniform bounds obtained in Section \ref{uniform} are used in the proofs of our analytical results concerning the relation between the eigenvalues of the slow and fast versions of the linearization. The theorems about the relation between the eigenvalues of the problem when the small parameters are zero and when they are small are obtained in Section \ref{DS}, using the concepts Evans function and the augmented unstable bundle as defined in \cite{AGJ}. {{We end Section \ref{DS} with concluding arguments to prove the main results of the paper.}} We perform numerical computations of the Evans function in Section \ref{S:Numerics} to {{explicitly identify regions in the parameter space where the fronts are spectrally stable in the appropriate weighted space for both cases $\epsilon$ zero and nonzero.}}

\section{The existence of fronts\label{RM}}

We consider a    diffusive Rosenzweig-MacArthur system  studied  in  \cite{CGM}: 
\eq{
u_{t}&=\epsilon_uu_{xx} +\frac{1}{\delta}\left(u\left(\beta-u\right)-\frac{uw}{1+u}\right),\\
w_{t}&=\epsilon_ww_{xx} + \frac{ w\left(u-\alpha\right)}{\eta+u}.
}{e:11_0}
Here, with an abuse of notations, we use the same variables $t$ and $x$ for the rescaled versions of $t$ and $x$  in \eqref{CGMeqn}. The parameter $\eta$ is related to the protection rates and $\beta$ is  related to the carrying capacity of the prey population. 

We also, for simplicity,   assume $\beta=1$,  so
 the system reads
\eq{
u_{t}&=\epsilon_uu_{xx} +\frac{1}{\delta}\left(u\left(1-u\right)-\frac{uw}{1+u}\right),\notag\\
w_{t}&=\epsilon_ww_{xx} + \frac{ w\left(u-\alpha\right)}{\eta+u},
}{e:11}
where  $\alpha$,  $\eta>0$ and $\epsilon_u$, $\epsilon_w$, $\delta\geq 0$. 

In this paper,  we  consider  the traveling fronts  in  \eqref{e:11} with
\begin{equation} \label{parameters}\text{$0<\alpha< 1$,  $\epsilon_u/\epsilon_w := \epsilon \ll 1$ and $0<\epsilon \ll \delta \ll 1$.} \end{equation}

  Since $0<\alpha<1$, we have three physically relevant equilibria $(u,w)$ of the system \eqref{e:11} given by 
 \begin{equation}A=(\alpha,1-\alpha^2), \,\,\,\, B=(1, 0), \,\,\,\,O=(0,0). \label{equilibria}\end{equation}

We replace $x$ with $z=x/\sqrt{\epsilon_w}$, and get
\eq{
u_{t}&=\epsilon u_{zz} +\frac{1}{\delta}\left(u\left(1-u\right)-\frac{uw}{1+u}\right),\\
w_{t}&= w_{zz} + \frac{ w\left(u-\alpha\right)}{\eta+u}.
}{e:12}
 The existence of the fronts is proven in \cite{CGM} for a variety of parameter regimes. 
 In particular,  in the regime \eqref{parameters}, \cite{CGM} explores the relation of the traveling fronts  to the fronts in  the Fisher-KPP equation  \cite{Fisher, KPP}. In this paper we show that the stability of these fronts are also inherited from the relevant dynamics of the Fisher-KPP equation. The general formulation of the Fisher-KPP equation is
\begin{equation}u_t =u_{zz} +f(u).\label{KPP}\end{equation}
We assume that  there are two values of $u$, say $u=0$ and $u=a$, such that
\begin{equation}  f(0)=f(a)=0, \quad f^{\prime} (0)>0,\quad  f^{\prime} (a) <0, \quad f^{\prime\prime}(u)<0. \label{KPPnon}\end{equation}
 
 We here briefly recap the existence proof.  The following theorem was proved in \cite{CGM} for system \eqref{e:12} in the case where $\epsilon\ll \delta$.
 \begin{Theorem} \label{T:1}
For every fixed  $0<\alpha<1$, $\eta >0$,  and every  $c > 2\sqrt{\frac{1-\alpha}{\eta+1}  }$,  there exists  $\delta_0=\delta_0(\alpha,\eta,c) >0$    such that for every $0<\delta <\delta_0$  there exists $\epsilon_0(\alpha, \eta,c,\delta)>0$  such that for each $0<\epsilon <\epsilon_0$  there is a  translationally invariant family of  fronts  of the system \eqref{e:12} which   move with speed $c$,  converge to the  equilibrium  $A=(\alpha, 1-\alpha^2)$  at $-\infty$  and  to   the equilibrium   $B =(1,0)$ at $+\infty$, and, moreover,  which have positive components $u$ and $w$.  
\end{Theorem}

 Fronts  of  \eqref{e:12}  that move with velocity $c>0$ are solutions of the following ODE system  
\eq{
0&={\epsilon} u_{\zeta\zeta} +c u_{\zeta}+ \frac{1}{\delta}\left(u\left(1-u\right)-\frac{uw}{1+u}\right),\\
0&=  w_{\zeta\zeta}+ c w_{\zeta}+\frac{ w\left(u-\alpha\right)}{\eta +u},}{e:tw} 
 where $\zeta = z-ct$. The   heteroclinic orbits of the associated system of the first order equations   represent fronts in the original PDE  \eqref{e:12}.
 In \cite{CGM} the system \eqref{e:tw}  {
was  treated as a dynamical system 
 \eq{
\delta \frac{du_1}{d\zeta}&=u_2,\\
\epsilon\frac{du_2}{d\zeta}&=-c u_2+\frac{u_1w_1}{1+u_1}- u_1(1-u_1),\\
\frac{dw_1}{d\zeta}&=w_2,\\
\frac{dw_2}{d\zeta}&=-c w_2-{}  \frac{w_1\left(u_1-\alpha\right)}{\eta+u_1}.
}{e:31_0}
In the limit $\epsilon \to 0$   the system  \eqref{e:31_0} is  reduced to a 3-dimensional system. The latter  as   $\delta \to 0$  is reduced to a system equivalent to a scalar ODE of the second order.}
\begin{equation}\frac{d^2w_1}{d\zeta^2 }+c \frac{dw_1}{d\zeta} + \frac{w_1\left(\sqrt{1-w_1}-\alpha\right)}{\eta+\sqrt{1-w_1}} =0. \label{reduction}\end{equation}
The equation \eqref{reduction} is a traveling wave equation for is a Fisher-KPP type equation 
\begin{equation}\frac{\partial w_1}{\partial t}=\frac{\partial^2w_1}{\partial  z^2 } + \frac{w_1\left(\sqrt{1-w_1}-\alpha\right)}{\eta+\sqrt{1-w_1}}. \label{reductionKPP}\end{equation}
As such it supports heteroclininc orbits connecting the equilibrium at $w=0$ with the equilibrium at $w= 1-\alpha^2$.  The dynamics of the limiting equation \eqref{reduction} is restricted in the phase space of \eqref{e:31_0} with $\epsilon=0$ to the nullcline $u_1\left(1-u_1\right)-\frac{u_1w_1}{1+u_1}=0$.

The analysis in \cite{CGM} shows that the heteroclinic orbits  that are known to exist in  the Fisher-KPP type equation \eqref{reduction} persist as a heteroclinic orbit in the full system \eqref{e:tw} with sufficiently small,  positive  $\delta$ and $\epsilon$,  { $\epsilon\ll\delta$}.

The following properties of the front will be useful in what follows. 
\begin{Lemma}  \label{L:d} Consider \eqref{e:12} when $\epsilon =0$. For every fixed  $0<\alpha<1$, $\eta >0$,  and every  $c > 2\sqrt{\frac{1-\alpha}{\eta+1}  }$,  there exists  $\delta_0=\delta_0(\alpha,\eta,c) >0$    such that for every $0<\delta <\delta_0$   the following estimates on the components of the front hold:
\begin{eqnarray}
&i)& 0< w_f(\zeta)\leq M_{w}, \text{ where } M_{w}= 1-\alpha^2 +o(\delta)  >0;\\
 &ii)&m_{u} < u_f(\zeta)\leq  1, \text{ where }  m_{u} = \alpha +o(\delta) >0.
\end{eqnarray}
\end{Lemma}
The proof of this lemma follows from the geometric construction of the front in  \cite{CGM} as a $o(\delta)$ perturbation of the front that exists in the singular limit $\delta \to 0$ and  the simple observations that (1) $w=0$ is an invariant set and (2) the derivative of $u$ is positive  for any value of $u>1$, so the equilibrium $(u,w)=(1, 0)$ can be reached only  when the $u$ component is increasing toward $1$. The front is constructed as a heteroclinic orbit in the dynamical system  associated with the traveling wave equation using the Geometric Singular Perturbation Theory \cite{Fenichel79, Jones94, kuehn}.

The  consequent construction of the front as a perturbation of that front  with a small diffusion  $\epsilon>0$   in \cite{CGM} gives the similar estimates.  Here we still use that  $w=0$ is an invariant set  and that the $o(\epsilon)$-perturbation of the front near the equilibrium $(u,w)=(1, 0)$ follows its  stable manifold which as a small perturbation of the stable manifold in case $\epsilon =0$, thus   stays in the region $w>0$, $u<1$ of the $(u,w)$-plane. 

\begin{Lemma}  \label{L:de}  For every fixed  $0<\alpha<1$, $\eta >0$,  and every  $c > 2\sqrt{\frac{1-\alpha}{\eta+1}  }$,  there exists $\delta_0=\delta_0(\alpha,\eta,c) >0$  such that for every $0<\delta <\delta_0$  there exists $\epsilon_0(\alpha, \eta,c,\delta)>0$  such that for each $0<\epsilon <\epsilon_0$ the following estimates on the components of the front hold:
\eq{
i)& \;0< w_f(\zeta)\leq M_{w}, \text{ where } M_{w}= 1-\alpha^2 +o(\delta)+o(\epsilon)  >0;\\
 ii)& \;m_{u} < u_f(\zeta)\leq 1,  \text{ where }  m_{u} = \alpha +o(\delta)+o(\epsilon) >0. 
}{frB}
{ In addition, the rate with which both components of  front converges to the origin  at $+\infty$  to the leading order may be described as
 \begin{equation}
 \label{myre}
 \frac{1}{2}\left(-c +\sqrt{c^2- 4\frac{1-\alpha}{1+\eta}}\right) +o(\delta)+o(\epsilon).
 \end{equation}}
\end{Lemma}

\section{Linearization about the front}

To investigate the stability of fronts  we first rewrite the  \eqref{e:12} in the coordinate $\zeta$, 
\eq{
u_{t}&=\epsilon u_{\zeta\zeta} + cu_{\zeta}+\frac{1}{\delta}f(u,w),\\
w_{t}&= w_{\zeta\zeta} +cw_{\zeta}+ g(u,w),
}{PDE_zeta}
where for brevity 
\begin{equation}f(u,w) = u\left(1-u\right)-\frac{uw}{1+u}, \qquad g(u,w)=\frac{ w\left(u-\alpha\right)}{\eta+u}.\end{equation}
We denote the front solution  $(u_f, w_f)$
and consider the eigenvalue problem for the operator  of the linearization of \eqref{PDE_zeta} about $(u_f, w_f)$
\eq{
\delta \lambda U&=\epsilon \delta U_{\zeta\zeta} +\delta cU_{\zeta}+f_u(u_f,w_f)U + f_w(u_f,w_f)W,\\
 \lambda  W&= W_{\zeta\zeta} +cW_{\zeta}+ g_u(u_f,w_f)U + g_w(u_f,w_f)W,
}{evdelta2}
where
\eq{ &f_u(u,w) = 1-  2u - \frac{w}{(1+u)^2},    \qquad f_w(u,w) = -\frac{u}{1+u},\\
  &g_u(u,w)=\frac{ (\eta +\alpha)w  }{(\eta+u)^2},  \,\, \,\qquad\qquad \quad  g_w(u,w)=\frac{ u-\alpha}{\eta+u}.
}{fderiv}
We note that the derivatives  belong to $\LT$.
The limiting values of these derivatives at the equilibria $A=(\alpha,1-\alpha^2)$  at $-\infty$ and $ B=(1, 0)$  at $+\infty$ are:
\eq{ &f_u(\alpha,1-\alpha^2) = - \frac{ 2\alpha^2}{1+\alpha},    \quad f_w(\alpha,1-\alpha^2) = -\frac{\alpha}{1+\alpha}, \quad f_u(1, 0) = -1,    \quad f_w(1, 0) = -\frac{1}{2},\\
 &  g_u(\alpha,1-\alpha^2) =\frac{1-\alpha^2}{\eta+\alpha}, \quad\quad g_w(\alpha,1-\alpha^2) =0, \quad \qquad g_u(1, 0)=0, \quad\quad  g_w(1, 0)=\frac{ 1-\alpha}{1+\eta}.
}{fderivlim}

In the analysis of the linearized operator,  in order to obtain bounds on its spectrum,  we shall  use the following estimates.

\begin{Lemma}\label{negativeBound1} Consider \eqref{e:12} when $\epsilon =0$. For every fixed  $0<\alpha<1$, $\eta >0$,  and every  $c > 2\sqrt{\frac{1-\alpha}{\eta+1}  }$,  there exists  $\delta_0=\delta_0(\alpha,\eta,c) >0$    such that for every $0<\delta <\delta_0$   there exists a constant $B>0$ such that
\begin{equation}f_u(u_f, w_f) \leq-\alpha^2+o(\delta)\leq  -B, \text{ for } \zeta\in \mathbb{R},\label{fuB}\end{equation}uniformly in $\delta$.\end{Lemma}
\begin {Proof}
It follows from the geometric construction of the front in the case $\epsilon=0$ and sufficiently small that $w_f=1-u_f^2+o(\delta)$,  therefore using notations of Lemma~\ref{L:d}
\begin{eqnarray}f_u(u_f,w_f) &=& 1-  2u_f - \frac{w_f}{(1+u_f)^2} =1-2u_f- \frac{1-u_f^2}{(1+u_f)^2} - \frac{o(\delta)}{(1+u_f)^2}\notag \\&=& 1-2u_f- \frac{1-u_f}{1+u_f} - \frac{o(\delta)}{(1+u_f)^2} = \frac{1-2u_f +u_f-2u_f^2 - 1+u_f}{1+u_f} - \frac{o(\delta)}{(1+u_f)^2} \notag \\
&=&  \frac{-2u_f^2}{1+u_f} - \frac{o(\delta)}{(1+u_f)^2}\leq  \frac{-2m_u^2}{2} +o(\delta)=  -\alpha^2+o(\delta), 
\end{eqnarray}
from where  the statement of the lemma follows.
\end{Proof}
The same argument, but based on Lemma~\ref{L:de} instead of Lemma~\ref{L:d} implies the following result.
\begin{Lemma}\label{negativeBound2} For every fixed  $0<\alpha<1$, $\eta >0$,  and every  $c > 2\sqrt{\frac{1-\alpha}{\eta+1}  }$,  there exists  $\delta_0=\delta_0(\alpha,\eta,c) >0$    such that for every $0<\delta <\delta_0$  there exists $\epsilon_0(\alpha, \eta,c,\delta)>0$  such that for each $0<\epsilon <\epsilon_0$ $B>0$ such that
\begin{equation}f_u(u_f, w_f) \leq -\alpha^2+o(\delta)+o(\epsilon) \leq -B, \text{ for } \zeta\in \mathbb{R},\end{equation}
uniformly in $\delta$ and $\epsilon$.\end{Lemma}

\section{{{Formulation of the Stability Result}}}
\label{SFEP}

{{We start this section by presenting the result about stability we prove in this paper. The proof of that theorem will be concluded in Section \ref{DS}.
 \begin{Theorem} \label{T:2}
For every fixed  $0<\alpha<1$, $\eta >0$,  and every  $c > 2\sqrt{\frac{1-\alpha}{\eta+1}  }$,  there exists  $\delta_0=\delta_0(\alpha,\eta,c) >0$    such that for every $0<\delta <\delta_0$  there exists $\epsilon_0(\alpha, \eta,c,\delta)>0$  such that for each $0<\epsilon <\epsilon_0$  the    fronts  of the system \eqref{e:12} of Theorem \ref{T:1} are spectrally stable in $\LT$ equipped with a weight
given by a strictly positive smooth function $w$ satisfying
$
w(\zeta)e^{-\sigma \zeta}\to 1,\;\;{\rm{as}}\;\;\zeta\to\infty$ and $w(\zeta)\to 1,\;\;{\rm{as}}\;\;\zeta\to-\infty,
$
for any $\sigma>0$ satisfying inequality \eqref{KPPE}, with the limits being attained exponentially fast. Moreover, the fronts are asymptotically stable in the $C^1$ weighted norm. 
\end{Theorem}
A simple example 
of a weight described in Theorem \ref{T:2} is the one used by Sattinger in \cite{Sattinger76}, of the general form $w=1+e^{\sigma \zeta}$, with the choice $\sigma=c/2$ producing the largest gap between the essential spectrum and the imaginary axis.  

Furthermore, applying the results from Sattinger \cite{Sattinger76}[Theorem 4.3], we have that the result on the spectral stability of each front implies asymptotic stability
as defined in \cite{Sattinger76}[Definition 2.1]. 

We point out that asymptotic stability is achieved because the weight used in Theorem \ref{T:2} removes the eigenvalue caused by the translational symmetry (at the origin). Indeed the eigenvector corresponding to the zero eigenvalue does not belong to the weighted space (see Lemma \ref{L:de}). As a consequence, the spectrum is bounded away from the imaginary axis on the left side of the complex plane.
}}

The idea for our reasoning is roughly as follows. If we consider \eqref{evdelta2} coupled to the traveling wave system  \eqref{e:31_0} in carefully chosen scaling of the moving coordinate. We then take the limits as  $\epsilon$ and then  $\delta$ approach zero and we expect that the limiting system is the eigenvalue problem for the  equation  \eqref{e:31_0}. We intend to rigorously show that  the eigenvalues of the system with sufficiently  small, positive $\epsilon$ and $\delta>0$, where the ``smallness'' will be defined later, are small  perturbations of the  eigenvalues of the limiting problem. 
 However, when we take the limit in $\epsilon$ and later $\delta$,  the outcome of the interplay between $\epsilon$, $\delta$, and $\lambda$ is not clear or obvious. To illuminate this issue to the reader, we show when it arises the first time in this process of the reduction. 

We start by rewriting  \eqref{evdelta2} as a first order dynamical system. To make it autonomous we  couple it  to the system representing the traveling wave equations \eqref{reductionKPP}.
\eq{
  \frac{dU_1}{d\zeta}&=U_2,\\
\epsilon  \frac{dU_2}{d\zeta}& =(\lambda - \frac{1}{\delta } f_u(u_1,w_1) )U_1  -cU_{2}-  \frac{1}{\delta }f_w(u_1,w_1)W_1,\\
 \frac{dW_1}{d\zeta}&=W_2,\\
 \frac{dW_2}{d\zeta}&= -g_u(u_1,w_1)U_1   +(\lambda - g_w(u_1,w_1))W_1-cW_{2},\\
  \frac{du_1}{d\zeta}&=u_2,\\
\epsilon\frac{du_2}{d\zeta}&=- c u_2+\frac{1}{\delta}f(u_1,w_1),\\
\frac{dw_1}{d\zeta}&=w_2,\\
\frac{dw_2}{d\zeta}&=-c w_2-g(u_1,w_1).
}{fastd}
In the slow scale $\xi =\zeta/\epsilon$, we obtain
\eq{
  \frac{dU_1}{d\xi}&=\epsilon U_2,\\
 \frac{dU_2}{d\xi}& = (\lambda - \frac{1}{\delta } f_u(u_1,w_1) )U_1  -cU_{2}-  \frac{1}{\delta }f_w(u_1,w_1)W_1,\\
 \frac{dW_1}{d\xi}&=\epsilon W_2,\\
 \frac{dW_2}{d\xi}&=  \epsilon \left( -g_u(u_1,w_1)U_1 +(\lambda - g_w(u_1,w_1))W_1 -cW_{2} \right),\\
 \frac{du_1}{d\xi}&=u_1,\\
\frac{du_2}{d\xi}&=- c u_2+\frac{1}{\delta}f(u_1,w_1),\\
\frac{dw_1}{d\xi}&= \epsilon w_2,\\
\frac{dw_2}{d\xi}&=\epsilon(-c w_2-g(u_1,w_1)).
}{slowd}
{The last four equations in \eqref{fastd} and \eqref{slowd} are decoupled from the eigenvalue problem and have been studied in detail in \cite{CGM}. The components $(u_1,w_1)$ represent the front solution $(u_f,w_f)$ and are $\epsilon$ and $\delta$ dependent, so we here denote them  as $(u_f[\epsilon, \delta],w_f[\epsilon, \delta])$. We will not write  the description of the independent variable in $(u_f,w_f)$ assuming that it is always consistent with the scaling of the variables. Moreover, we know from \cite{CGM} that the limit of  $(u_f[\epsilon, \delta],w_f[\epsilon, \delta])$  as $\epsilon\to 0$  as well as the limit when  $\delta\to 0$ exist}.  Let us focus on the first four equations  of the system \eqref{fastd} that correspond to the eigenvalue problem. In the limit when $\epsilon=0$,  the equation \eqref{fastd} becomes 
\eq{
  \frac{dU_1}{d\zeta}&=\frac{1}{c}\left(  (\lambda - \frac{1}{\delta } f_u(u_f[0, \delta],w_f[0, \delta]) )U_1 -  \frac{1}{\delta }f_w(u_f[0, \delta],w_f[0, \delta])W_1     \right),\notag\\
 \frac{dW_1}{d\zeta}&=W_2,\notag\\
 \frac{dW_2}{d\zeta}&=  -g_u(u_f[0, \delta],w_f[0, \delta])U_1 +(\lambda - g_w(u_f[0, \delta],w_f[0,t \delta]))W_1  -cW_{2} ,
}{fastd0}
on $$ U_{2} =\ \frac{1}{c}\left( (\lambda - \frac{1}{\delta } f_u(u_f[0, \delta],w_f[0, \delta]) )U_1 -  \frac{1}{\delta }f_w(u_f[0, \delta],w_f[0, \delta])W_1\right)$$
 and  the corresponding equations in \eqref{slowd} become
\eq{
  \frac{dU_1}{d\xi}&=0,\notag\\
 \frac{dU_2}{d\xi}& = (\lambda - \frac{1}{\delta } f_u(u_f[0, \delta],w_f[0, \delta]) )U_1-cU_{2} -  \frac{1}{\delta }f_w(u_f[0, \delta],w_f[0, \delta])W_1,\notag\\
 \frac{dW_1}{d\xi}&=0,\notag\\
 \frac{dW_2}{d\xi}&=  \lim_{\epsilon\to 0}(\epsilon \lambda) W_1,}
{slowd0}
if the finite $\lim_{\epsilon\to 0}(\epsilon \lambda)$ exists. In order to be able to carry our analysis forward, we need to know the value of this limit. Similar issues arise with $\delta$ and $\lambda$.
We will address this issue in Sections \ref{ESSECT} and \ref{uniform}. 

In the next section we find the location of the right most boundary of the essential spectrum  which is related to the dynamics near the rest states of the front. Then we spend the rest of the paper working on  the discrete spectrum.


 \section{Essential spectrum}
 \label{EssSect}
We consider operator $L$ defined by the right hand side of \eqref{evdelta2}
\begin{equation}\label{L}
L\begin{pmatrix}U\\W\end{pmatrix} : = \begin{pmatrix} 
 \epsilon  U_{\zeta\zeta} +cU_{\zeta}+\frac{1}{\delta} f_u(u_f,w_f)U &  \frac{1}{\delta} f_w(u_f,w_f)W \\
g_u(u_f,w_f)U&  W_{\zeta\zeta} +cW_{\zeta} + g_w(u_f,w_f)W
\end{pmatrix}
\end{equation}
The resolvent  of a  closed, densely defined operator $L$ is a set of $\lambda$ in the complex plane such that  $(L-\lambda I)^{-1}$  exists and is bounded.  The complement of the resolvent set is the spectrum. The latter may  contain isolated eigenvalues of finite algebraic multiplicity as well as continuous spectrum, or, the so-called essential spectrum.  The boundaries of the essential spectrum  are among the curves of the spectrum of constant coefficient operators obtained by taking the limits $L_{\pm\infty}$ of   \eqref{L} as $\zeta \to \pm\infty$. 
Recall that the heteroclinic orbit asymptotically converges to the equilibrium $B = (1, 0)$ at $+\infty$ and to he equilibrium $A = ( \alpha, 1-\alpha^2)$ at $-\infty$, and  so using \eqref{fderiv} we obtain
\begin{equation}\label{L+}
L_{+}\begin{pmatrix}U\\W\end{pmatrix} : = \begin{pmatrix} 
 \epsilon  U_{\zeta\zeta} +cU_{\zeta}-\frac{1}{\delta} U & - \frac{1}{2 \delta} W \\
0 &  W_{\zeta\zeta} +cW_{\zeta} + \frac{ 1-\alpha}{\eta+1} W
\end{pmatrix} 
\end{equation}
and 
\begin{equation}\label{L-}
L_{-}\begin{pmatrix}U\\W\end{pmatrix} : = \begin{pmatrix} 
 \epsilon  U_{\zeta\zeta} +cU_{\zeta}-\frac{1}{\delta} \frac{ 2\alpha^2 }{1+\alpha}  U &  - \frac{1}{\delta}\frac{\alpha}{1+\alpha}W \\
\frac{ (1-\alpha^2) }{(\eta+\alpha)}U&  W_{\zeta\zeta} +cW_{\zeta} 
\end{pmatrix} 
\end{equation}



The spectra of these constant-coefficient  differential operators  are obtained through Fourier analysis. At $+\infty$, the spectrum of $L_{+}$ is defined by 
\eq{
 \lambda= -\epsilon  k^2 +cik -\frac{1}{\delta} ,\\
\lambda= -k^2 +cik + \frac{ 1-\alpha}{\eta+1}.
}{curve+}

These are two parabolas in the complex $\lambda$-plane 
\eq{
\mathrm{Re}\, \lambda= -\frac{\epsilon}{c^2} (\mathrm{Im}\, \lambda)^2 -\frac{1}{\delta} ,\\
\mathrm{Re}\, \lambda= -\frac{1}{c^2} (\mathrm{Im}\, \lambda)^2  + \frac{ 1-\alpha}{\eta+1}.
}{curve++}
Note that in the case $\epsilon=0$, the first equation defines a vertical line on the left side of the complex plane.
The spectrum  is unstable due to the second curve intersecting with the right side of the complex plane, since $0<\alpha<1$ as specified by Theorem \ref{T:1}. 

At $-\infty$,   the spectrum of $L_{-}$ is defined by 
\begin{equation}\label{curve-}
\rm{det} \begin{pmatrix} 
- \epsilon k^2 +cki -\frac{1}{\delta} \frac{ 2\alpha^2 }{1+\alpha}  -\lambda  &  - \frac{1}{\delta}\frac{\alpha}{1+\alpha} \\
\frac{1-\alpha^2}{\eta+\alpha}&  -k^2 +cik -\lambda 
\end{pmatrix} =0.
\end{equation}
We insert $\lambda=\lambda_r+ i\lambda_i$ into \eqref{curve-} and separate  the real and imaginary parts of the determinant. We find that the imaginary part if the determinant in \eqref{curve-} is zero when either  
\begin{equation} \lambda_r=-\frac{\left({\epsilon}+1\right) k^2}{2}-\frac{\,{\alpha}
^{2}}{ \left( \alpha+1 \right) \delta}
\label{first_2}\end{equation}
or 
\begin{equation}
\lambda_i=c\,k.
\label{first_1}\end{equation}
We thus have two cases to consider separately:  \eqref{first_2} or \eqref{first_1}.
We note that in the case of \eqref{first_2}  no unstable essential spectrum is produced since $\lambda_r$ is strictly negative. We proceed to find the corresponding curves of the spectrum. 

\begin{figure}[h]
\vspace*{0mm}
\hspace{-0cm}
\scalebox{.4}{\includegraphics{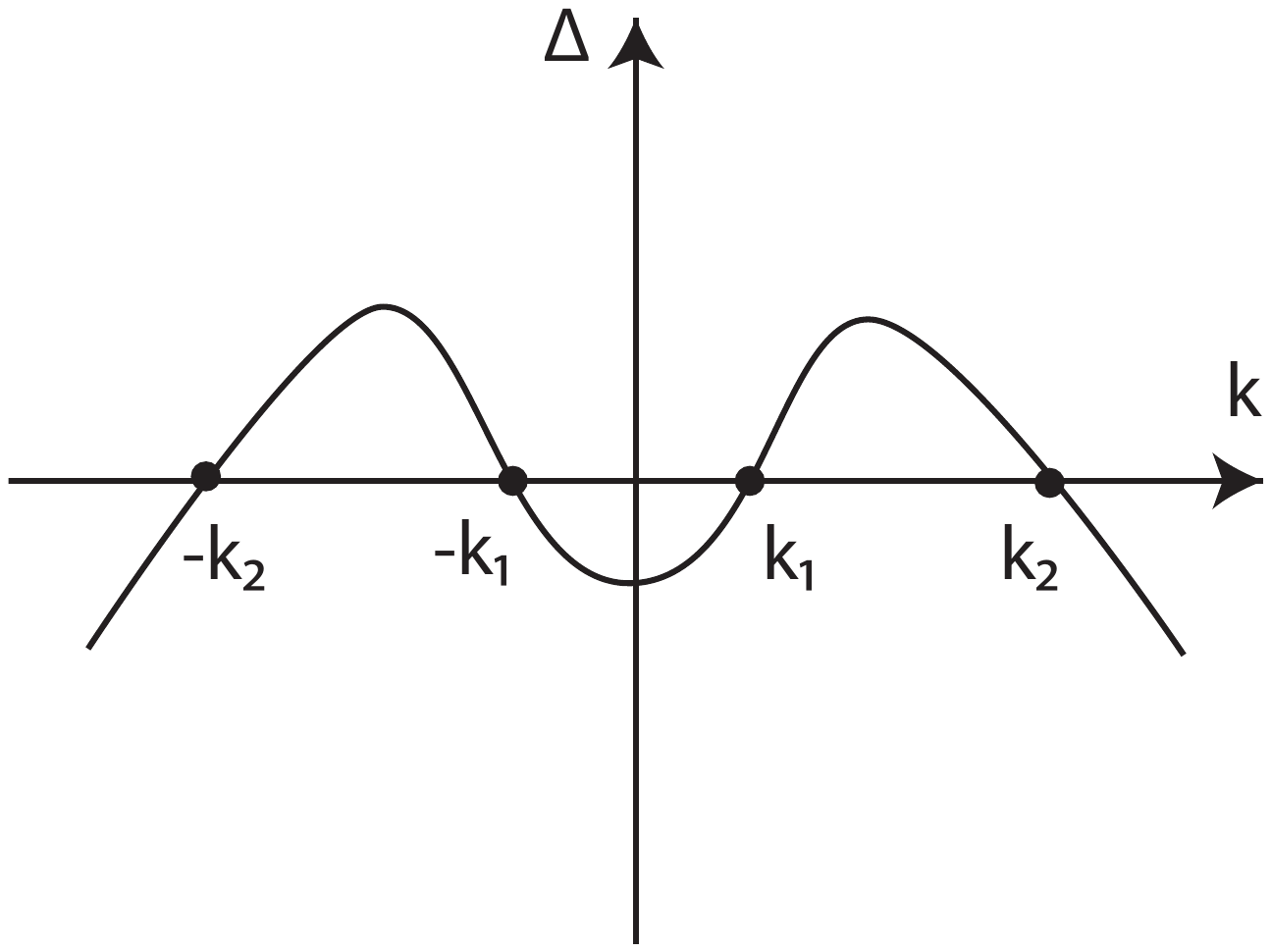}}
\vspace*{00mm}
\caption{\label{Del}  Graph of the discriminant $\Delta$ given in \eqref{disc}, where $k=\pm k_1, \pm k_2$,  $0<k_1<k_2$, are 
its four real roots.}
\end{figure}

When \eqref{first_2} holds, we find that the real part of \eqref{curve-}  gives the
following condition on $\lambda_i$:
\eq{
-\lambda_i^2+2ck\lambda_i
-\frac{\left( {\epsilon}-1 \right) ^{2}{k}^{4}}{4}-{\frac { \left( {c}^{2}
\delta(1+\alpha)+{\alpha}^{2}({\epsilon}-1)
 \right) {k}^{2}}{\delta\, \left( 1+\alpha \right) }}-{\frac {\alpha\,
 \left( {\alpha}^{3}(\alpha+\eta)-\delta(1-\alpha^2)(\alpha^2+1)\right)}
{{\delta}^{2} \left( 1+\alpha
 \right) ^{2} \left( \eta+\alpha \right)}}=0.
}{firstc}
The discriminant of the equation above as a polynomial in $\lambda_i$ can be computed to be
\eq{
\Delta=- \left( {\epsilon}-1 \right) ^{2}{k}^{4}+4\,{\frac { \left( 1-{\epsilon} \right) {\alpha}^{2}{k}^{2}}{\delta\, \left( 1+\alpha \right) }}-4\,
{\frac {\alpha\, \left( {\alpha}^{3}(\alpha+\eta)-\delta(1-\alpha^2)(\alpha+1) \right) }{{\delta}^{2}
 \left( 1+\alpha \right) ^{2} \left( \eta+\alpha \right) }},
}{disc}
whose discriminant as a quadratic in $k^2$ is given by
\eq{
16\,{\frac { \left( {\epsilon}-1 \right) ^{2} \left( 1-\alpha \right) \alpha}{\delta\, \left( \eta+\alpha \right) }}.
}{disc2}
Since the last discriminant is positive, \eqref{disc} has two real roots for $k^2$. Given that
$\delta$ is small enough to satisfy
\eq{
 \delta<\frac{{
\alpha}^{3} \left( \eta+\alpha \right)}
{\left( 1+\alpha \right)  \left(1-{\alpha}^{2}\right)},
 }{kzeroc}
 then the two roots of \eqref{disc} are positive values of $k^2$. Let us denote the corresponding roots of \eqref{disc} as $k=\pm k_1$ and $k=\pm k_2$ with $0<k_1<k_2$. The discriminant \eqref{disc} is thus positive only on the two open intervals $(k_1,k_2)$ and $(-k_2,-k_1)$ (see Figure \ref{Del}). On those two intervals,
 \eqref{firstc} has two real roots for every value of $k$. Furthermore, given that the coefficient of $\lambda_i^2$ of the polynomial in \eqref{firstc} is $-1$ and the constant coefficient is negative, the real solutions of  \eqref{firstc} have the same sign as $k$.

 \begin{Lemma}\label{3} Under Condition \eqref{kzeroc}, and the conditions described in Theorem \ref{T:1}, the discriminant \eqref{disc} has four roots $k=\pm k_1$ and $k=\pm k_2$, $0<k_1<k_2$. For each value of $k$ on $(k_1,k_2)$ or $(-k_2,-k_1)$, \eqref{firstc} has two solutions $\lambda_i=\lambda_i(k)$, $i=1,2$, corresponding to the imaginary part of the essential spectrum curve for $\zeta\rightarrow -\infty$. Furthermore, $\lambda_i(k)>0$ on $(k_1,k_2)$ and $\lambda_i(k)<0$ on $(-k_2,-k_1)$.  The real part is given by the equation of \eqref{first_2} and thus is negative.
 \end{Lemma}

If we now choose the solution \eqref{first_1}, we find that the real part of \eqref{curve-} gives rise to the equation
\eq{
\lambda_r^2+\left((\epsilon+1)k^2+\frac{2\alpha^2}{\delta(1+\alpha)}\right)\lambda_r+\epsilon k^4+\frac{2\alpha^2k^2}{\delta(1+\alpha)}+\frac{(1-\alpha)\alpha}{(\eta+\alpha)\delta}=0.}{roots}
Since its coefficients are all positive, the real roots of this polynomial are all negative. Furthermore, this polynomial has as its discriminant $-\Delta$, with $\Delta$ given in \eqref{disc}. Thus the discriminant of \eqref{roots} is positive for the values of $k$ in the complement of the set described in Lemma \ref{3}.  
\begin{Lemma}
\label{4}
Consider the roots $k_i$ described in Lemma \ref{3}. Then 
for every values of $k$ on the intervals
$(-\infty,-k_2)$, $(-k_1,k_1)$, and $(k_2,\infty)$, the polynomial \eqref{roots} has two negative real roots $\lambda_r=\lambda_r(k)$
corresponding to the real part of the essential spectrum curve for $\zeta\rightarrow -\infty$. The imaginary part of the curve is given by 
$\lambda_i=ck$. 
\end{Lemma}

\begin{figure}[h]
\vspace*{0mm}
\hspace{-0cm}
\scalebox{.36}{{\includegraphics{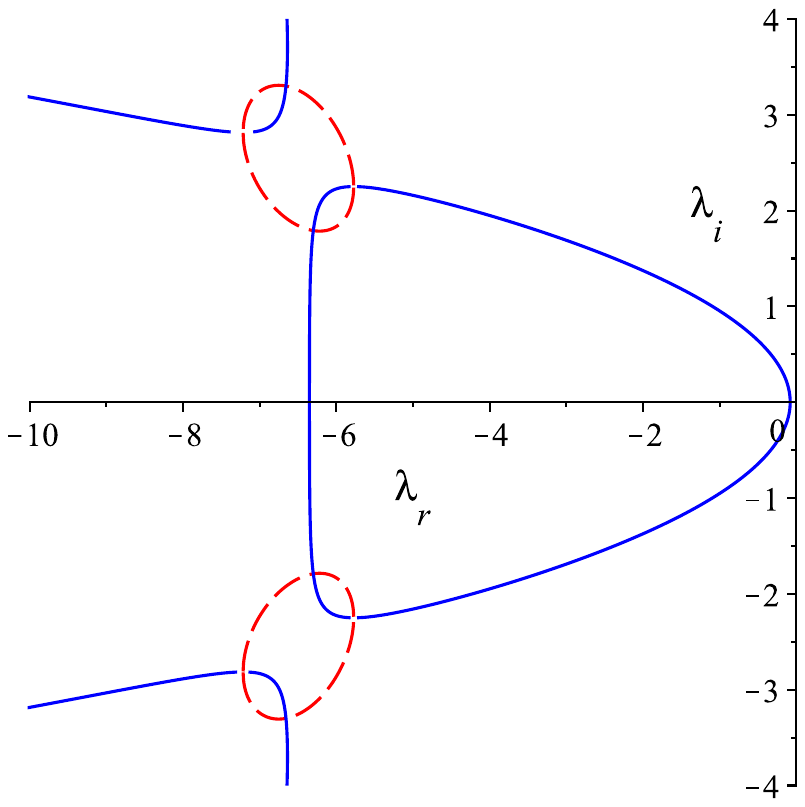}}}
\scalebox{.36}{{\includegraphics{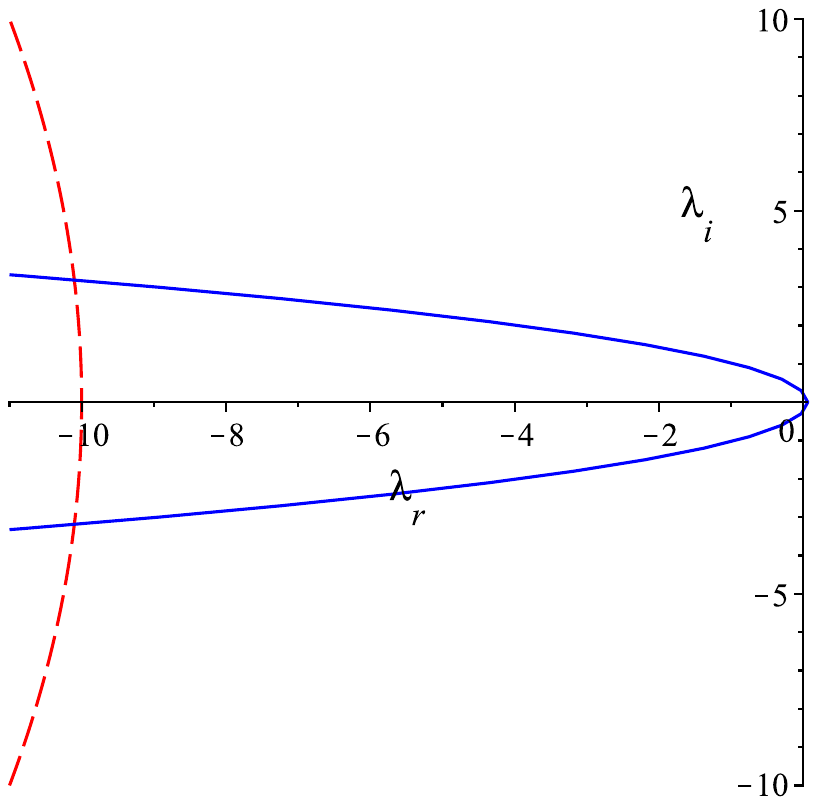}}}\\
\vspace*{-45mm}
\scalebox{.36}{{\includegraphics{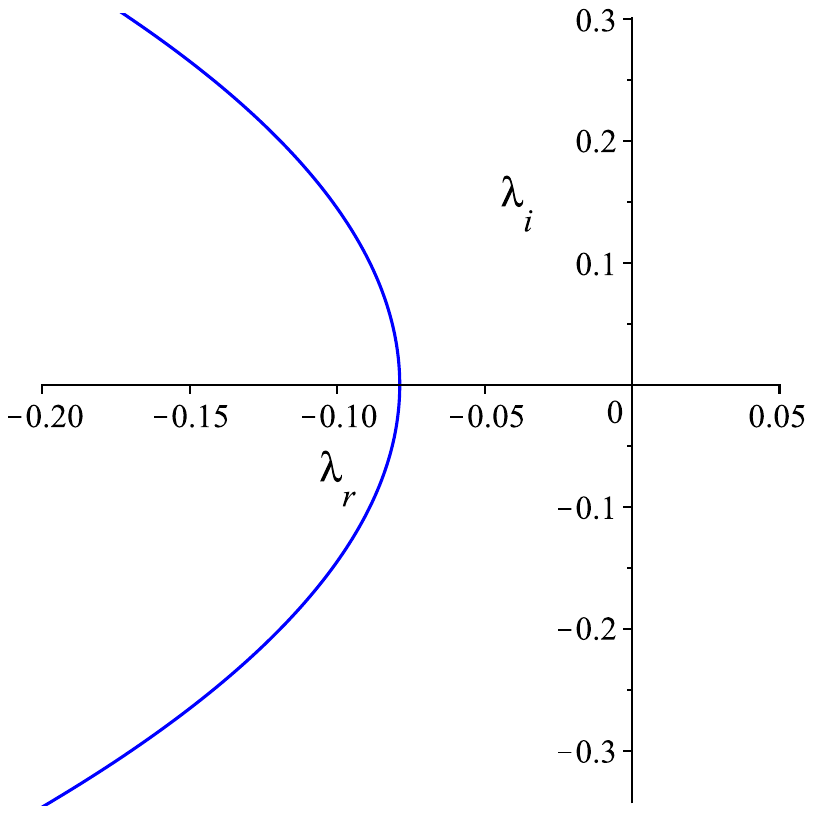}}}
\scalebox{.36}{{\includegraphics{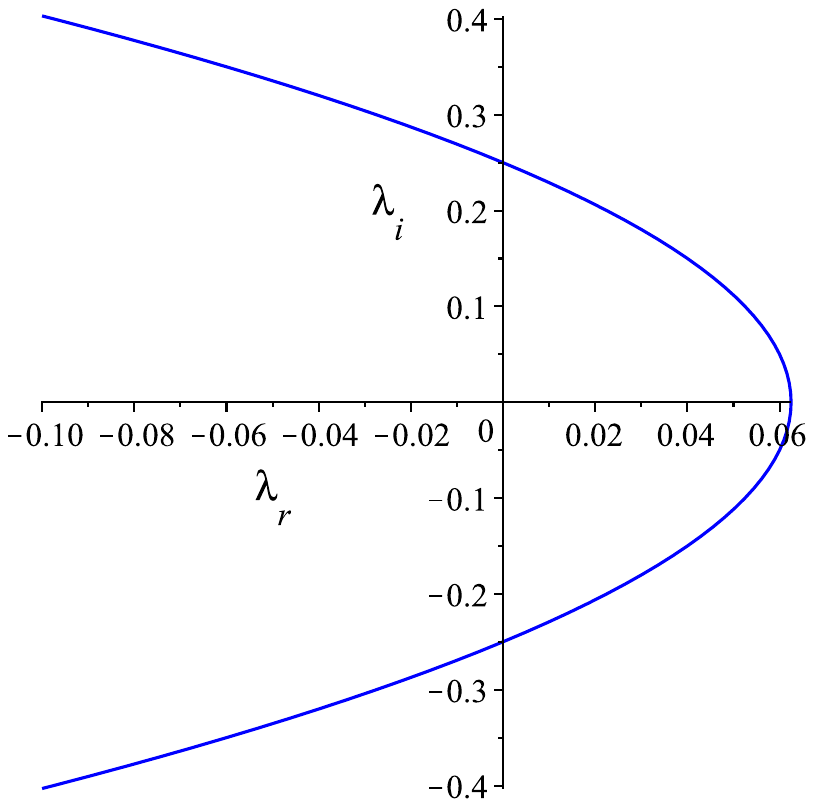}}}
\vspace*{-40mm}
\caption{\label{Ess1} The boundaries of the essential spectrum in the case  
$\alpha=0.75$, $\delta=0.1$, $\epsilon=0.01$, $\eta=3$ and $c=1$. The top left figure illustrates the 
spectrum boundaries for the case $\zeta\rightarrow -\infty$ 
 as described by Lemmas \ref{3} and \ref{4}, corresponding to $k_1=2.2506$ and $k_2=2.8146$. The two dashed 
closed curves correspond to Lemma \ref{3}, while the rest to Lemma \ref{4}. The bottom left shows a zoom around the origin. The gap between the curve and the imaginary axis for the figures on the left is found to be 0.0787 by taking the absolute value of the expression given in \eqref{spgap}. The top right figure illustrates the 
spectrum boundaries for the case $\zeta\rightarrow \infty$ as described by  equations \eqref{curve++}.
The dashed curve corresponds to the first equation in \eqref{curve++}, while the solid curve corresponds to the second equation.
The bottom right shows a zoom around the origin.
}
\end{figure}

Figure \ref{Ess1} shows the essential spectrum curves from \eqref{curve++} (for the contribution from $+\infty$) and from lemmas \ref{3} and \ref{4} (for $-\infty$) in the case  
$\alpha=0.75$, $\delta=0.1$, $\epsilon=0.01$, $\eta=3$ and $c=1$.

{We can be more precise and show that the closest point of the curves defined by lemmas \ref{3} and \ref{4} occurs when we set $k=0$, thus at one of the two points at which the curve crosses the real axis. To prove this, we take $k\in(0,k_1)$ and write the largest real solution to \eqref{roots} as
\eq{
\lambda_r=\frac{-\left((\epsilon+1)k^2+\frac{2\alpha^2}{\delta(1+\alpha)}\right)+\sqrt{-\Delta}}{2},
}{realg}
where $\Delta$ is given in \eqref{disc}, $-\Delta$ being the discriminant to the polynomial \eqref{roots}. On the interval $(0,k_1)$, the quantity $\Delta$ is increasing (see Figure \ref{Del}), thus the expression given above in \eqref{realg} decreases on that interval. From that fact and the fact that equation  \eqref{roots} is symmetric in $k$, one concludes that the maximum value of $\lambda_r$ on the interval $(-k_1,k_1)$ is at $k=0$. At $k=0$, the curve crosses the real axis and the corresponding value of the real part  is found by solving \eqref{roots} for $k=0$, which results in the expression}
\begin{equation}\label{spgap}
\lambda_r(0)={\frac {-{\alpha}^{2}(\eta+{\alpha})+\sqrt {\alpha\, \left( \eta+
\alpha \right) 
 \left( {\alpha}^3\left( \alpha+\eta \right)-\delta(1-\alpha^2)(1+\alpha)\right) }}{\delta
\, \left( 1+\alpha \right)  \left( \eta+\alpha \right) }}.
\end{equation}

For the interval $(k_2,\infty)$, we use implicit differentiation of \eqref{roots} with respect to $k$ to obtain
\eq{
{\lambda_r'=-\frac{k\left(2(\epsilon+1)\lambda_r+4k^2\epsilon +\frac{4\alpha^2}{\delta(1+\alpha)}\right)
}{2\lambda_r+(\epsilon+1)k^2+\frac{2\alpha^2}{\delta(1+\alpha)}}.}
}{imp}
Using equation \eqref{realg} for the largest root of \eqref{roots}, we find that the denominator in \eqref{imp} is positive  for $k\in(k_2,\infty)$. Using the same formula \eqref{realg} and the fact that $\Delta=0$ at $k=k_2$, we find that the numerator in \eqref{imp} is negative for $k=k_2$. Thus, by continuity, the numerator in \eqref{imp} is negative $k>k_2$ close enough to $k_2$, that is
\eq{
\lambda_r\leq-2\left(\frac{\epsilon k^2 +\frac{\alpha^2}{\delta(1+\alpha)}}{\epsilon+1}\right).
}{In1}
The derivative thus initially satisfies $\lambda_r'>0$ on the interval $(k_2,\infty)$ and $\lambda_r$ reaches the value given 
by the RHS of \eqref{In1} evaluated at some  $k=k'\in (k_2,\infty)$. Since the sign of $\lambda_r'$ changes if Inequality \eqref{In1} is not satisfied, $\lambda_r$ reaches its maximum at $k=k'$. It thus follows that
\eq{
\lambda_r\leq-{\frac{2\alpha^2}{\delta(\epsilon+1)(1+\alpha)}}
}{In2}
for all $k\in(k_2,\infty)$. Since the value on the RHS of \eqref{In2} is smaller than $\lambda_r(0)$ given in \eqref{spgap} if $\epsilon<1$, the maximum value of $\lambda_r$ in Lemma \ref{4} occurs at $k=0$.

The value given in \eqref{spgap} is also larger than any value of the real part of any point on the curve defined 
in Lemma \ref{3}. Indeed, the real part of the curve defined in Lemma \ref{3} is given in equation \eqref{first_2}. It is easy to verify that to ensure
that the value 
given above in  \eqref{spgap} is greater than the value given in \eqref{first_2} for any real value of $k$, it suffices that
$$
{{\frac {\sqrt {\alpha\, \left( \eta+
\alpha \right) 
 \left( {\alpha}^3\left( \alpha+\eta \right)-\delta(1-\alpha^2)(1+\alpha)\right) }}{\delta
\, \left( 1+\alpha \right)  \left( \eta+\alpha \right) }}>0.}
$$
Under Condition \eqref{kzeroc} and the conditions described in Theorem \ref{T:1}, the inequality above is satisfied. 

We thus have the following lemma
concerning the distance of the essential spectrum curve at $+\infty$ and the imaginary axis
\begin{Lemma}\label{5} Under Condition \eqref{kzeroc}, and the conditions described in Theorem \ref{T:1}, the essential spectrum curve defined by Lemmas 
 \ref{3} and \ref{4} is on the open left side of the complex plane and it  is the closest to the imaginary axis at one of the two points where it crosses the real axis. The point of intersection of the curve and the real axis closest to the imaginary axis is at the value given in \eqref{spgap}.
 \end{Lemma}


In order to stabilize the essential spectrum at $+\infty$ as determined by \eqref{curve+}, we use the  {{weight introduced in Theorem \ref{T:2}  given by a strictly positive smooth function $w$ satisfying
\eq{
w(\zeta)e^{-\sigma \zeta}\to 1,\;\;{\rm{as}}\;\;\zeta\to\infty,\\
w(\zeta)\to 1,\;\;{\rm{as}}\;\;\zeta\to-\infty,
}{disW}
for some $\sigma>0$, where the limits are attained exponentially fast. A simple example 
of such a weight is the one used by Sattinger in \cite{Sattinger76}, of the general form $w=1+e^{\sigma \zeta}$.

To study the essential spectrum in the weighted space, we apply the change of variable $(\widetilde{U},\widetilde{V})=w(\zeta)(U,V)$ to the eigenvalue problem for the operator $L$  given in \eqref{L}. We obtain the following problem on $\LT$ (now without weight) 
\eq{\left(L-\lambda\right) \begin{pmatrix}\widetilde{U}\\\widetilde{W}\end{pmatrix}+
\begin{pmatrix} 
 \epsilon \left( \left(\frac{2w'^2-ww''}{w^2}\right)\widetilde{{{{{{{U}}}}}}} -\frac{2w'}{w}\widetilde{U}_\zeta\right)-\frac{cw'}{w}\widetilde{U} & 0 \\
0&    \left(\frac{2w'^2-ww''-cww'}{w^2}\right)\widetilde{W} -\frac{2w'}{w}\widetilde{W}_{\zeta}
\end{pmatrix}=0.
}{LP}
As we did to compute the curves \eqref{curve+}, we compute the boundary of the essential spectrum of the eigenvalue problem \eqref{LP}
at $\zeta\to \infty$ to obtain the following relations 
\eq{
 \lambda= -\epsilon  k^2 +i(c-2\epsilon \sigma)k +\sigma(\epsilon \sigma-c)-\frac{1}{\delta} ,\\
\lambda= -k^2 +i(c-2\sigma)k + \sigma(\sigma-c)+\frac{ 1-\alpha}{\eta+1} .
}{curve+s}}}
These two curves are restricted to the left side of the complex plane if 
\eq{
\sigma(\epsilon \sigma-c)-\frac{1}{\delta}<0 ,\;\;{\rm{and}}\;\;
\sigma(\sigma-c)+\frac{ 1-\alpha}{\eta+1}<0. 
}{curve+in}
Using the conditions on the parameters as described in Theorem \ref{T:1}, the solutions to
\eq{
\sigma(\sigma-c)+\frac{ 1-\alpha}{\eta+1}=0
}{q1}
 are both real negative while the solutions to
 \eq{
\sigma(\epsilon \sigma-c)-\frac{1}{\delta}=0
}{q2}
 are real and of opposite sign. Thus in order to choose a $\sigma$ satisfying both inequalities, one needs to choose $\sigma$ between the largest zero of \eqref{q1} and both the smallest zeroes of \eqref{q1} and  \eqref{q2}. This gives rise to the condition
\eq{
\frac{c}{2}-\frac{1}{2}\sqrt{c^2-4\frac{1-\alpha}{\eta+1}}<\sigma<\max\left(\frac{c}{2}+\frac{1}{2}\sqrt{c^2-4\frac{1-\alpha}{\eta+1}},\;\frac{c}{2\epsilon}+\frac{1}{2\epsilon}\sqrt{c^2+\frac{4\epsilon}{\delta}}\right).
}{sigma+c}
It is straightforward to show that the inequality \eqref{sigma+c} above has a nonempty solution set if $\epsilon<1$. Actually, if $\epsilon$ is small enough to satisfy
$$
0<\epsilon<\frac{c\delta \sqrt{c^2-4\frac{1-\alpha}{\eta+1}}+c^2\delta+2}{\delta\left(c\sqrt{c^2-4\frac{1-\alpha}{\eta+1}}+c^2-2\frac{1-\alpha}{\eta+1}\right)},
$$
then inequality \eqref{sigma+c} becomes
\eq{
\frac{c}{2}-\frac{1}{2}\sqrt{c^2-4\frac{1-\alpha}{\eta+1}}<\sigma<\frac{c}{2}+\frac{1}{2}\sqrt{c^2-4\frac{1-\alpha}{\eta+1}}.
}{KPPE}

In the example $\alpha=0.75$, $\delta=0.1$, $\epsilon=0.01$, $\eta=3$ and $c=1$ considered above, one finds the 
interval given by \eqref{sigma+c} to be
$$
0.067<\sigma<0.93.
$$

\begin{figure}[h]
\vspace*{0mm}
\hspace{-0cm}
\scalebox{.36}{{\includegraphics{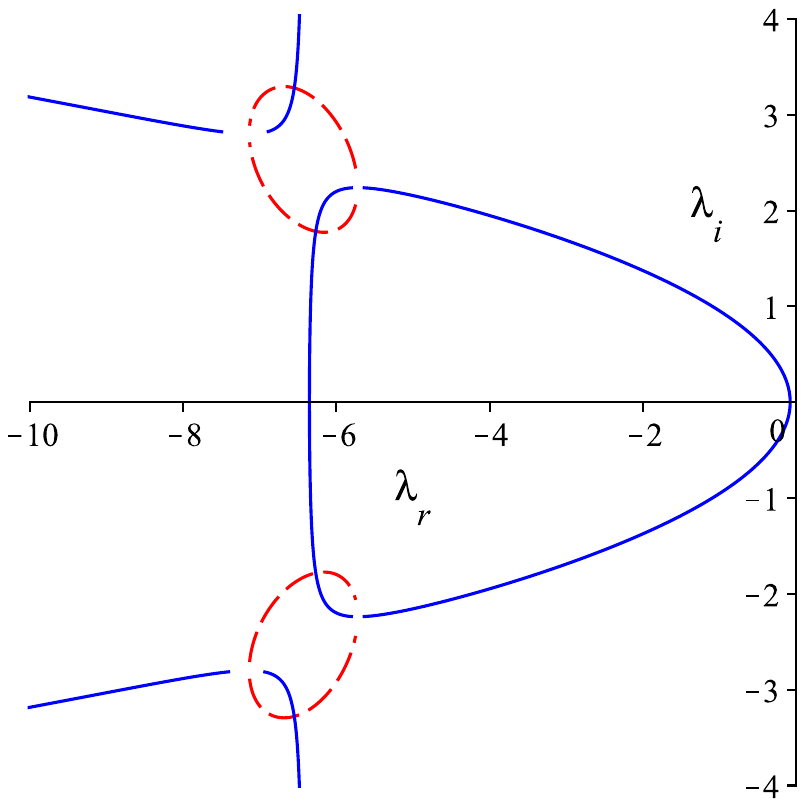}}}
\scalebox{.36}{{\includegraphics{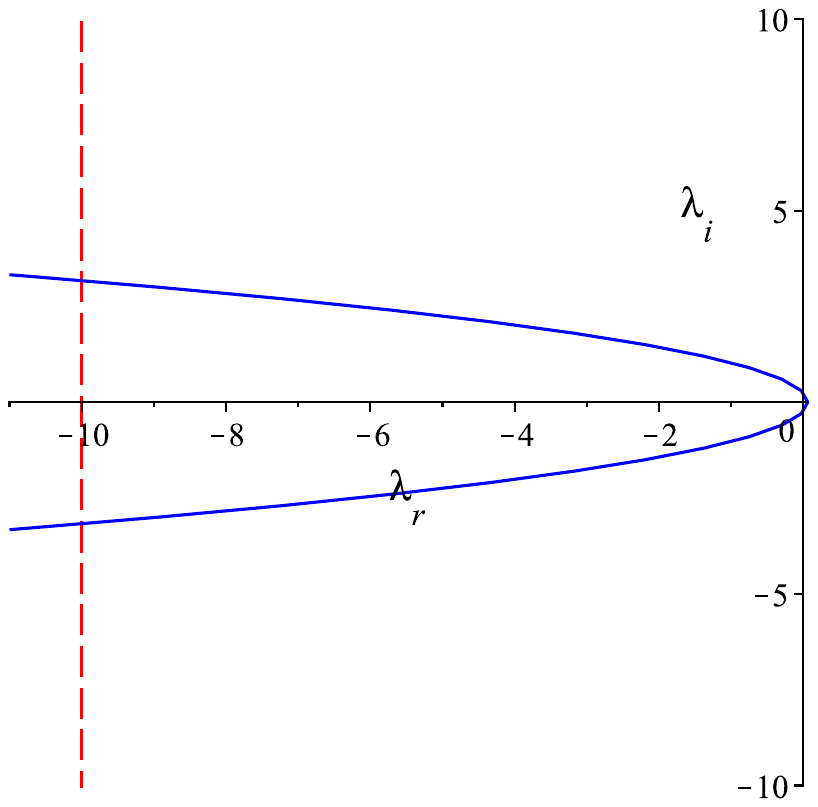}}}\\
\vspace*{-45mm}
\scalebox{.36}{{\includegraphics{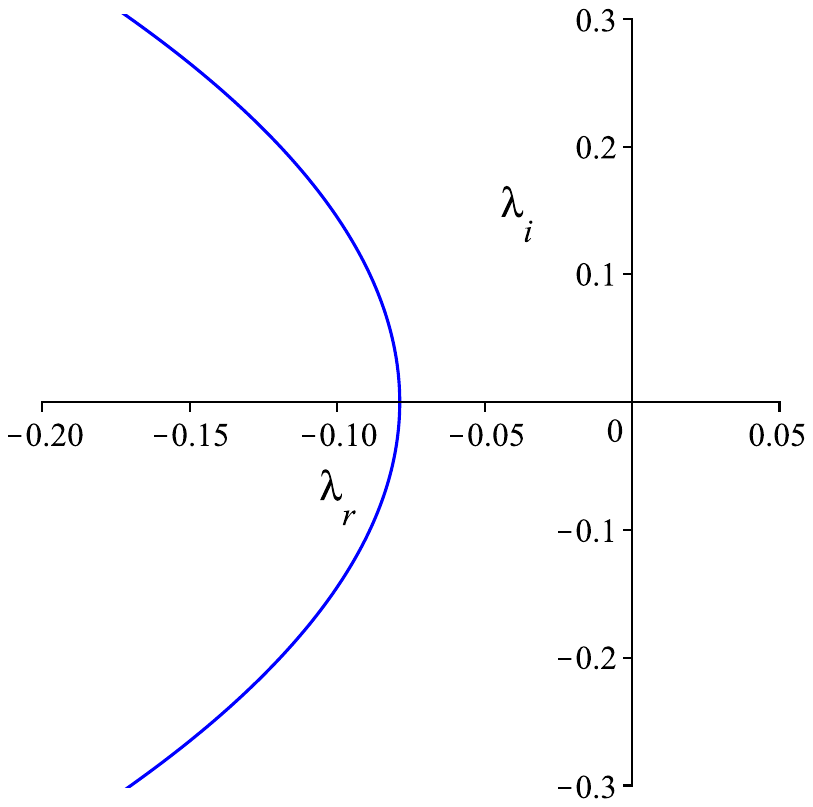}}}
\scalebox{.36}{{\includegraphics{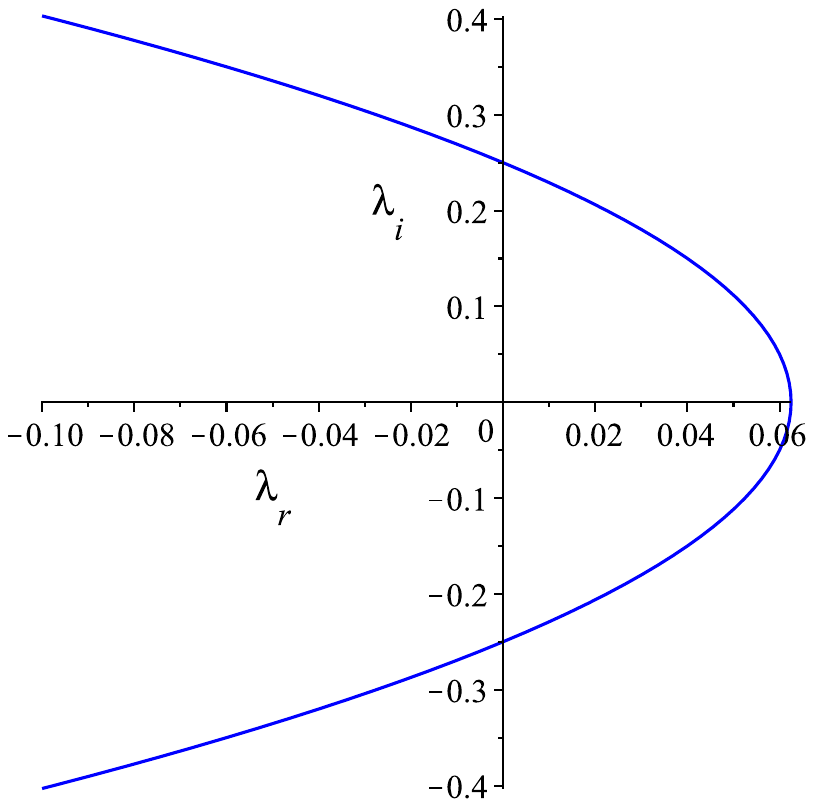}}}
\vspace*{-40mm}
\caption{\label{Ess2} The boundaries of the essential spectrum for the same parameter values as in Figure \ref{Ess1} ($\alpha=0.75$, $\delta=0.1$, $\eta=3$ and $c=1$.), except that $\epsilon=0$.
 The top left figure illustrates the 
spectrum boundaries for the case $\zeta\rightarrow -\infty$ 
 as described by Lemmas \ref{3} and \ref{4}, corresponding to $k_1=2.2393$ and $k_2=2.8005$. The two dashed 
closed curves correspond to Lemma \ref{3}, while the rest to Lemma \ref{4}. The bottom left shows a zoom around the origin. The top right figure illustrates the 
spectrum boundaries for the case $\zeta\rightarrow -\infty$ as described by  equations \eqref{curve++}.
The dashed curve corresponds to the first equation in \eqref{curve++}, while the solid curve corresponds to the second equation.
The bottom right shows a zoom around the origin.
}
\end{figure}

In the case where $\epsilon=0$ and $\delta\neq 0$, the boundary of the essential spectrum as $\zeta\rightarrow -\infty$ is also described by lemmas \ref{3} and \ref{4} and, qualitatively, looks the same. In the case where $\zeta\rightarrow\infty$, the boundary is still given by \eqref{curve++} but the first equation now describes a vertical line on the left side of the complex plane rather than a parabola. Figure \ref{Ess2} shows the boundary of the essential spectrum for the same parameter values as in Figure \ref{Ess1}, except that $\epsilon$ is set to 0.

In  the singular limits $\epsilon \to 0$ and  $\delta \to 0$ ($\epsilon\ll \delta $) we obtain  corresponding to the Fisher-KPP scalar equation
\eqref{reductionKPP}. In that case, the eigenvalue problem is given in \eqref{KPPlin}. 
 The essential spectrum is obtained by computing the spectra of the constant-coefficient problems obtained by applying the 
 limits $\zeta\rightarrow \pm \infty$ to \eqref{KPPlin}, using the fact that $w_f\rightarrow 0$ (resp.~ $w_f\rightarrow 1-\alpha^2$) as $\zeta\rightarrow \infty$ (resp.~$\zeta\rightarrow -\infty$). Using Fourier analysis, one finds that the essential spectrum is defined by the curves
 \eq{
\lambda&= -k^2 +cik + \frac{ 1-\alpha}{\eta+1},\\
\lambda&= -k^2 +cik - \frac{ 1-\alpha^2}{2\alpha(\alpha+\eta)}.
}{curvesKPP}
To stabilize the first curve of \eqref{curvesKPP}, we use a weight as in \eqref{disW} and obtained the 
second equation of \eqref{curve+s}. This curve is restricted to the left side of the complex plane if inequality \eqref{KPPE} is satisfied.
As mentioned in Section \ref{ESSECT}, it is known from the work of Sattinger \cite[Theorem 6.3]{Sattinger76} that the spectrum of the operator arising from the linearization of the Fisher-KPP
equation about the ``fast'' fronts (with $c$ satisfying the condition of Theorem \ref{T:1}) can be moved to the left side of the complex plane through the use {of a carefully chosen}  exponential weight.

\section{Bounds on the point spectrum}
\label{ESSECT}

The goal of that section is to obtain a bound on $|\lambda|$, for any eigenvalue in the point spectrum on the right side of the complex plane, i.e. such that $\realpart{\lambda}>0$. The point spectrum includes the discrete eigenvalues that are the zeros of the Evans function that we define in a subsequent section.   To obtain such a bound on unstable eigenvalues (if they exist), we consider the
 eigenvalue problem \eqref{evdelta2}, which we rewrite here:
\eq{
\lambda U&=\epsilon  U_{\zeta\zeta} +cU_{\zeta}+\frac{1}{\delta }f_u(u_f,w_f)U + \frac{1}{\delta }f_w(u_f,w_f)W,\\
 \lambda  W&= W_{\zeta\zeta} +cW_{\zeta}+ g_u(u_f,w_f)U + g_w(u_f,w_f)W.
}{evdelta}
We multiply the first equation by $\overline U$ and the second by $\overline W$ and integrate to obtain
\eq{
\lambda \intr |U|^2 d \zeta &=\epsilon  \intr U_{\zeta\zeta} \overline U d \zeta+c\intr U_{\zeta} \overline U  d \zeta+\frac{1}{\delta }\intr f_u(u_f,w_f)|U|^2 d \zeta+ \frac{1}{\delta }\intr f_w(u_f,w_f)W\overline U d \zeta,\\
 \lambda  \intr |W|^2 d \zeta &= \intr W_{\zeta\zeta} \overline W d \zeta+c \intr W_{\zeta}\overline W d \zeta+ \intr g_u(u_f,w_f)U \overline W d \zeta+  \intr g_w(u_f,w_f)|W|^2 d \zeta.
}{evdeltaint}
Since we assume $\lambda$ to be an eigenvalue, $U$ and $W$ are exponentially localized so we can integrate by parts. In particular, we have that
$$ \intr U_{\zeta\zeta} \overline U d \zeta = -\intr |U_{\zeta} |^2 d \zeta, \quad  \intr W_{\zeta\zeta} \overline W d \zeta = -\intr |W_{\zeta} |^2 d \zeta.$$
Another useful fact is that 
$$\intr U_{\zeta} \overline U d \zeta = - \intr \overline U_{\zeta}  U  d \zeta.$$
Therefore, $\intr U_{\zeta} \overline U  d \zeta$ and similarly $\intr W_{\zeta} \overline W  d \zeta$ are purely imaginary. 
It is then easy to see that, after taking the real parts of the \eqref{evdeltaint}, we obtain
\eq{
\rm{Re}(\lambda) \intr |U|^2 d \zeta \leq &  \frac{1}{\delta}\sup\{|f_u(u_f,w_f)|\}
 \intr |U|^2  d \zeta + \frac{1}{\delta}\sup\{|f_w(u_f,w_f)|\}\intr |\rm{Re}(W\overline U)| d \zeta\\&-\epsilon \intr |U_{\zeta} |^2 d \zeta,\\
 \rm{Re}(\lambda)  \intr |W|^2 d \zeta \leq & \sup\{|g_w(u_f,w_f)|\}
  \intr |W|^2  d \zeta+ \sup\{|g_u(u_f,w_f)|\}
  \intr|\rm{Re}( U \overline W)| d \zeta\\&
  -\intr |W_{\zeta} |^2 d \zeta.}{firstin}
We {{consider the following version of Young's inequality:
\eq{
ab\leq  a^2\beta_*+\frac{b^2}{4\beta_*},\;\;{\mbox{for any}}\;\;a,b\in\R,\text{ and }\beta_*>0.
}{Young}
We apply \eqref{Young} to the product of the mixed terms in \eqref{firstin}  and obtain
$$
|\rm{Re}( U \overline W)|=|\rm{Re}( W \overline U)|\leq |W||U|\leq
  |W|^2\beta_*+  \frac{|U|^2}{4\beta_*},\;\;{\mbox{for any}}\;\;\beta_*>0.
$$
We use the inequality above with $\beta_*=\beta_1$ and $\beta_*=\beta_2$, and}} rewrite \eqref{firstin} as follows 
\eq{
\rm{Re}(\lambda) \intr |U|^2 d \zeta \leq &\frac{1}{\delta}\sup\{|f_w(u_f,w_f)|\}\intr \left(\beta_1 |W|^2 +\frac{1}{4\beta_1}|U|^2\right) d \zeta\\&  +\frac{1}{\delta}\sup\{|f_u(u_f,w_f)|\}  \intr |U|^2  d \zeta  -\epsilon \intr |U_{\zeta} |^2 d \zeta ,\\ 
\rm{Re}(\lambda)  \intr |W|^2 d \zeta \leq &   \sup\{|g_u(u_f,w_f)|\} \intr \left(\beta_2 |U|^2 +\frac{1}{4\beta_2}|W|^2\right) d \zeta\\&+\sup\{|g_w(u_f,w_f)|\} \intr |W|^2  d \zeta -\intr |W_{\zeta} |^2 d \zeta. }{re}
After adding these two inequalities together  and ignoring the negative terms with the first order derivatives, we obtain 
\begin{eqnarray*}&& \rm{Re}(\lambda) \intr ( |U|^2 +|W|^2 ) d \zeta \\
&&\quad \leq\left( \frac{1}{\delta}\sup\{|f_u(u_f,w_f)|\}+  \beta_2 \sup\{|g_u(u_f,w_f)|\}  + \frac{1}{4\beta_1} 
\frac{1}{\delta}\sup\{|f_w(u_f,w_f)|\}  \right)\intr |U|^2  d \zeta \\ &&\qquad+ \left( \sup\{|g_w(u_f,w_f)|\}+\frac{\beta_1}{\delta}\sup\{|f_w(u_f,w_f)|\} + \frac{1}{4\beta_2 }\sup\{|g_u(u_f,w_f)|\}     \right) \intr  |W|^2 d \zeta,
\end{eqnarray*}
where   $\beta_1$ and $\beta_2$  are any positive numbers. The latter inequality implies  a bound on the real part of $\lambda$ if that real part is  positive, 
\eq{\rm{Re}(\lambda) &\leq  \max\left\{\frac{1}{\delta}\sup\{|f_u(u_f,w_f)|\} +  \beta_2 \sup\{|g_u(u_f,w_f)|\}+ \frac{1}{4\beta_1}\frac{1}{\delta}\sup\{|f_w(u_f,w_f)|\}  ,   \right . 
\\ & \qquad \qquad\left.
 \sup\{|g_w(u_f,w_f)|\}+\frac{\beta_1}{\delta}\sup\{|f_w(u_f,w_f)|\} + \frac{1}{4\beta_2 }\sup\{|g_u(u_f,w_f)|\}    \right \}.}{reall}
Instead of finding a  bound on $\rm{Im}(\lambda)$ for a $\lambda$ with $\rm{Re}(\lambda)>0$, if it exists, we plan to find an estimate for $\rm{Re}(\lambda)+ |\rm{Im}(\lambda)| $. By taking the imaginary part of \eqref{evdeltaint}, we obtain the following inequalities
\begin{eqnarray*}
|\rm{Im}(\lambda)| \intr |U|^2 d \zeta &\leq & c\intr |U_{\zeta}|| \overline U|  d \zeta+ \frac{1}{\delta}\sup\{|f_w(u_f,w_f)|\}\intr |W||U| d \zeta,\notag\\
|\rm{Im}(\lambda)|  \intr |W|^2 d \zeta &\leq &c \intr |W_{\zeta}|| W| d \zeta+ \sup\{|g_u(u_f,w_f)|\} \intr |U| |W| d \zeta,
\end{eqnarray*}
where we used the fact that $|\rm{Im}( W \overline U)|\leq |U||W|$.
We again use Young's inequality \eqref{Young}  to the terms $|W||U|$ and obtain
\eq{
|\rm{Im}(\lambda)| \intr |U|^2 d \zeta &\leq  c\intr |U_{\zeta}|| \overline U|  d \zeta+  \frac{1}{\delta}\sup\{|f_w(u_f,w_f)|\}  \intr \left(\beta_3 |W|^2+\frac{1}{4\beta_3}|U|^2\right)
 d \zeta,\\ 
|\rm{Im}(\lambda)|  \intr |W|^2 d \zeta &\leq c \intr |W_{\zeta}|| W| d \zeta+\sup\{|g_u(u_f,w_f)|\}   \intr  \left(\beta_4 |W|^2+\frac{1}{4\beta_4}|U|^2\right)  d \zeta.
}{im}
At this point we add the inequalities \eqref{re} and \eqref{im} together, 
\eqnn{
 &\left(\rm{Re}(\lambda)+ |\rm{Im}(\lambda)|\right) \intr (|U|^2  +  |W|^2  ) d \zeta \\
 & \quad  \leq
 \left[ \frac{1}{\delta}\sup\{|f_u(u_f,w_f)|\} +   \beta_2 \sup\{|g_u(u_f,w_f)|\} + \frac{1}{4\delta \beta_1}\sup\{|f_w(u_f,w_f)|\} \right.
  + \frac{1}{4\delta\beta_3}\sup\{|f_w(u_f,w_f)|\}\\
  &\qquad\left.   +  \frac{1}{4\beta_4}\sup\{|g_u(u_f,w_f)|\} \right ] \intr |U|^2  d \zeta
  \\&\qquad   + \left[ \sup\{|g_w(u_f,w_f)|\}+\frac{\beta_1}{\delta}\sup\{|f_w(u_f,w_f)|\} + \frac{1}{4\beta_2 }\sup\{|g_u(u_f,w_f)|\} \right. +  \frac{\beta_3}{\delta}\sup\{|f_w(u_f,w_f)|\} \\ 
 & \qquad    + \beta_4\sup\{|g_u(u_f,w_f)|\}    \bigg]\intr |W|^2 d \zeta 
 -\epsilon \intr |U_{\zeta} |^2 d \zeta  -\intr |W_{\zeta} |^2 d \zeta +  c\intr |U_{\zeta}|| \overline U|  d \zeta  +  c \intr |W_{\zeta}|| W| d \zeta.
}
We apply Young's inequality \eqref{Young} to obtain: 
\eq{
  c\intr |U_{\zeta}|| \overline U|  d \zeta &\leq \intr \left(\epsilon |U_{\zeta}|^2 +\frac{c^2}{4 \epsilon}| \overline U|^2\right)  d \zeta, \\ 
    c \intr |W_{\zeta}|| W| d \zeta &\leq   \intr \left(|W_{\zeta}|^2+\frac{c^2}{4}| W|^2\right) d \zeta,
}{Y1}
to modify the estimate above as 
\eqnn{
 &\left(\rm{Re}(\lambda)+ |\rm{Im}(\lambda)|\right) \intr (|U|^2  +  |W|^2  ) d \zeta\notag \\
 & \quad  \leq
 \left[\frac{c^4}{4\epsilon}+ \frac{1}{\delta}\sup\{|f_u(u_f,w_f)|\} +  \beta_2 \sup\{|g_u(u_f,w_f)|\}  + \frac{1}{4\delta\beta_1}\sup\{|f_w(u_f,w_f)|\}\right. \\ 
 & \qquad \left.  + \frac{1}{4\delta\beta_3}\sup\{|f_w(u_f,w_f)|\}       +\frac{1}{4\beta_4}\sup\{|g_u(u_f,w_f)|\} \right] \intr |U|^2  d \zeta  \notag \\ 
 & \qquad + \left[ \frac{c^4}{4}+ \sup\{|g_w(u_f,w_f)|\}+\frac{\beta_1}{\delta} \sup\{|f_w(u_f,w_f)|\} + \frac{1}{4\beta_2 }\sup\{|g_u(u_f,w_f)|\}  \right.\\ 
 & \qquad+  \frac{\beta_3}{\delta}\sup\{|f_w(u_f,w_f)|\}    + \beta_4\sup\{|g_u(u_f,w_f)|\}   \bigg]\intr |W|^2 d \zeta.\notag\\ 
}
We  notice from the expressions given in \eqref{fderiv} that 
\eq{\sup\{|f_w(u_f,w_f)|\} =\frac{1}{2}\text{ and }\sup\{|g_w(u_f,w_f)|\} =\frac{1-\alpha}{\eta+1}}{supexpl} 
and we state our result in the form of a lemma. 
\begin{Lemma} \label{L:2} For any eigenvalue $\lambda$  of the eigenvalue problem \eqref{evdelta} with $\rm{Re}(\lambda)>0$, the following estimate holds 
\begin{equation}\label{reim}
|\lambda|\leq \max\{M_1, M_2\} , \end{equation}
where 
\eq{
M_1&= \frac{c^4}{4\epsilon}+ \frac{1}{\delta}\left( \sup\{|f_u(u_f,w_f)|\}+\frac{1}{8\beta_1} + \frac{1}{8\beta_3}\right) + \sup\{|g_u(u_f,w_f)|\}\left( \beta_2  +  \frac{1}{4\beta_4}\right) 
   ,   \\
M_2&= \frac{c^4}{4}+ \frac{1-\alpha}{\eta+1}
+ \left(\frac{1}{4\beta_2}   + \beta_4 \right) \sup\{|g_u(u_f,w_f)|\}  +  \frac{\beta_1+\beta_3}{2\delta}
, 
 }{M1M2}
and  $\beta_1$, $\beta_2$, $\beta_3$, and $\beta_4$ are arbitrary positive constants. 
\end{Lemma}

\begin{Remark}
The bound given in Lemma \ref{L:2}, for $\epsilon$ small enough, is given by $M_1$ since $M_1={\mathcal{O}}\left(\frac{1}{\epsilon}\right)$ as $\epsilon\to 0$. Furthermore, 
Expression \eqref{reall} shows that if $\lambda$ is an eigenvalue on the right side of the complex plane, $\realpart{\lambda}$ is uniformly bounded in $\epsilon$. Indeed,  given the expressions for the derivatives of $f$ given 
in \eqref{fderiv}, and given the bounds on the components \eqref{frB}, and also given \eqref{supexpl} and \eqref{fuB}, one concludes that $M_2$ defined in \eqref{M1M2} and the RHS of \eqref{reall} can be bounded independently of 
$\epsilon$, for $\epsilon$ small enough.
\end{Remark}

We now obtain estimates on  $\lambda$ in the eigenvalue problem \eqref{evdelta} in the case when $\epsilon=0$. Clearly, the process that we used to obtain Lemma~\ref{L:2} does not apply in its completeness to this case since the first equation in \eqref{evdelta} now is of first order. Our starting point is the system: 
\eq{
\lambda U&=cU_{\zeta}+\frac{1}{\delta }f_u(u_f,w_f)U + \frac{1}{\delta }f_w(u_f,w_f)W,\\
 \lambda  W&= W_{\zeta\zeta} +cW_{\zeta}+ g_u(u_f,w_f)U + g_w(u_f,w_f)W.
}{evdelta0}
Since system \eqref{evdelta0} is obtained from \eqref{evdelta} by setting $\epsilon=0$, the  inequalities \eqref{re} hold with $\epsilon=0$, which in turn gives rise to 
\eqref{reall}. We thus have the following inequality in the case where $\lambda$ is an eigenvalue of the problem \eqref{evdelta0} with positive real part: 
\eq{\rm{Re}(\lambda)\leq  \max&\left\{\frac{1}{\delta}\sup\{|f_u(u_f,w_f)|\}   +  \beta_2 \sup\{|g_u(u_f,w_f)|\} + \frac{1}{8\beta_1}\frac{1}{\delta},\right.
\\  &
\quad \left.\frac{1-\alpha}{\eta+1}  +\frac{\beta_1}{2\delta}  + \frac{1}{4\beta_2 }\sup\{|g_u(u_f,w_f)|\}      \right\}
 }{reallc}
where $\beta_1$ and $\beta_2$ are arbitrary positive real numbers. Note that we have used the expressions for the supremums given in \eqref{supexpl} 
to write \eqref{reallc} from \eqref{reall}.

From \eqref{evdeltaint} with $\epsilon=0$, we have
\eq{| \impart{\lambda}|  \intr |U|^2 d \zeta \leq c \intr |U_{\zeta}||U| d \zeta+ \frac{1}{\delta} \intr f_w(u_f,w_f)|\impart{W \overline U}| d \zeta.
}{fisrt1}

We can improve the results  of Lemma~\ref{L:2} by showing that the bound on the  eigenvalues with positive real parts,  if there are any,  does not increase as $\delta$ decreases, by means of the argument  described below.

We write the solution of  the following ODE of the first order 
\eq{U_{\zeta} -\frac{1}{c} \left(\lambda-\frac{1}{\delta }f_u(u_f,w_f)\right)U=  - \frac{1}{c \delta }f_w(u_f,w_f)W,}{first}
 in terms of $W$, using the integrating factor
 \eq{\mu(\zeta)=\rm{exp}\left\{  
 \int^{\zeta}_{0}  \left(\frac{1}{c\delta} f_u(u_f(s),w_f(s)) -\frac{1}{c}\lambda  \right) d s   \right\}.
}{mudef}
Since, from \eqref{fderivlim}, the limits of $f_u(u_f(s),w_f(s))$ as $s\to\pm\infty$ are negative, and since $\realpart{\lambda}>0$, we have that $\mu\to 0$ as $\zeta\to + \infty$ and to $\mu\to\infty$ as $\zeta\to - \infty$. After  multiplying the equation \eqref{first}  by $\mu$ we obtain
  $$(U(\zeta) \mu(\zeta) )^{\prime}= -\frac{1}{ c \delta}\mu(\zeta) f_w(u_f(\zeta),w_f(\zeta))W(\zeta),$$
  which we integrate from $\zeta$ to $\infty$
  to obtain the unique solution of \eqref{first} that converges at $\pm \infty$ given by 
   $$U(\zeta) =  \frac{1}{ c \delta} \int_{\zeta}^{+\infty}  \frac{\mu(r)  }{  \mu(\zeta)}f_w(u_f(r),w_f(r))W(r) dr.
    $$  
We then have the following inequality 
     \begin{equation*}|U(\zeta)| \leq   \frac{1}{ c \delta} \int_{\zeta}^{+\infty} \left| \frac{\mu(r)  }{  \mu(\zeta)}\right| |f_w(u_f(r),w_f(r))||W(r)| dr\leq   \frac{1}{ c \delta} \sup\{| f_w(u_f,w_f)| \}   \int_{\zeta}^{+\infty} \left| \frac{\mu(r)  }{  \mu(\zeta)}\right| |W(r)| dr.
    \end{equation*}
   Using \eqref{mudef}, we have
    \begin{eqnarray*}\left|\frac{ \mu(r)}{\mu(\zeta)} \right|  &\leq&  \rm{exp}\left\{  
 \int^{r}_{\zeta}  \left(\frac{1}{c\delta} \sup\{f_u(u_f,w_f)\} -\frac{1}{c}Re\lambda  \right) d s   \right\} \notag
 \\&=&   \rm{exp}\left\{   \left( \frac{1}{c\delta} \sup\{f_u(u_f,w_f)\} -\frac{1}{c}Re\lambda  \right)( r -\zeta )   \right\}.
 \end{eqnarray*}
 We define
\eqnn{{{T}}(r-\zeta)\equiv
\left\{ 
\begin{array}{cl}
\displaystyle{\rm{exp}\left\{  
 \left( \frac{1}{c\delta} \sup\{f_u(u_f,w_f)\} -\frac{1}{c}Re\lambda  \right)( r -\zeta )\right\}},   &\text{ when }r>\zeta,\\
 \displaystyle{0},&\text{ otherwise}.
 \end{array}\right.
 }
 Then 
  \begin{equation}|U(\zeta)| \leq   \frac{1}{ c \delta} \sup\{| f_w(u_f,w_f)| \}   \int_{-\infty}^{+\infty} {\red{T}}( r -\zeta ) |W(r)| dr.
    \end{equation}
   Then by the Young's inequality for convolutions  $\|F*{{T}}\|_{s}\leq \|F\|_{p}\|{{T}}\|_{q}$ when $\frac {1}{p}+{\frac {1}{q}}={\frac {1}{s}}+1$  for $s=2$, $p=2$, and $q=1$,  we have 
   $$\|U \|_2=  \sup\{| f_w(u_f,w_f)| \}  \|W\|_2  \frac{1}{ c \delta} \int_{-\infty}^{+\infty} {{T}}( r -\zeta )  dr. $$  
We compute
\begin{equation} \int_{-\infty}^{+\infty} {{T}}( r -\zeta ) d r
 = \frac{c\delta }{-\sup\{f_u(u_f,w_f)\} +\delta\realpart{\lambda}}\leq  -\frac{c\delta }{\sup\{f_u(u_f,w_f)\} },
 \end{equation}
 where we use the fact that $\sup\{f_u(u_f,w_f)\} <0$ from \eqref{supexpl}.
 So, 
 \begin{equation}\|U \|_2=  \frac{\sup\{| f_w(u_f,w_f)| \}}{-\sup\{f_u(u_f,w_f)\} } \|W\|_2, 
 \end{equation} where we note that 
 $$I= \frac{\sup\{| f_w(u_f,w_f)| \}}{-\sup\{f_u(u_f,w_f)\} } >0,$$
 so 
 \begin{equation}\|U \|_2=  I \|W\|_2.
 \end{equation}
We estimate, after using Holder's inequality
\begin{eqnarray*}\left |\intr g_u(u_f,w_f)U \overline W d \zeta\right |\leq 
\sup\{g_u(u_f,w_f)\}  \sqrt{\intr  |U|^2 d\zeta}    \sqrt{\intr    |W|^2     d \zeta }\\
 \leq  \sup\{g_u(u_f,w_f)\}   I  \sqrt{\intr    |W|^2     d \zeta }\sqrt{\intr    |W|^2     d \zeta }\leq \sup\{g_u(u_f,w_f)\}I  \intr    |W|^2     d \zeta   \notag \end{eqnarray*}
 and 
  $$ c\intr |W_{\zeta}||\overline W| d \zeta \leq  \intr |W_{\zeta}|^2 d \zeta +\frac{c^2}{4} \intr  | W|^2 d \zeta,$$
   thus, 
\begin{eqnarray*}   \lambda  \intr |W|^2 d \zeta &=& \intr W_{\zeta\zeta} \overline W d \zeta+c \intr W_{\zeta}\overline W d \zeta+ \intr g_u(u_f,w_f)U \overline W d \zeta+  \intr g_w(u_f,w_f)|W|^2 d \zeta \notag \\ &= & - \intr |W_{\zeta}|^2d \zeta+c \intr W_{\zeta}\overline W d \zeta+ \intr g_u(u_f,w_f)U \overline W d \zeta+  \intr g_w(u_f,w_f)|W|^2 d \zeta
\end{eqnarray*}
and
\begin{eqnarray*} && | \lambda|  \intr |W|^2 d \zeta  \\
&&\quad \leq   - \intr |W_{\zeta}|^2d \zeta+\intr (|W_{\zeta}|^2+\frac{c^2}{4}|\overline W|^2) d \zeta + \intr |g_u(u_f,w_f)U \overline W d \zeta|+  
\sup\{|g_w(u_f,w_f)|\} \intr |W|^2 d \zeta
\notag \\ &&\quad \leq   \frac{c^2}{4}\intr  |W|^2 d \zeta  + \sup\{g_u(u_f,w_f)\}I  \intr    |W|^2     d \zeta +  \sup\{|g_w(u_f,w_f)|\} \intr |W|^2 d \zeta
\end{eqnarray*}
Therefore, we have  proved the following bound on the eigenvalues
\begin{equation} \label{delta0W} |\lambda| \leq  \frac{c^2}{4}   + \sup\{g_u(u_f,w_f)\}\frac{\sup\{| f_w(u_f,w_f)| \}}{\sup\{|f_u(u_f(s),w_f(s))|\} }  + 
 \sup\{|g_w(u_f,w_f)|\}.
\end{equation}
From \eqref{delta0W}, we use the expressions given in \eqref{supexpl}  for the supremums of $| f_w(u_f,w_f)|$  and $|g_w(u_f,w_f)|$ and state the following lemma.
\begin{Lemma} \label{L:3} Consider  the eigenvalue problem \eqref{evdelta}  with $\epsilon=0$. For any eigenvalue $\lambda$   with $\rm{Re}(\lambda)>0$ the following estimate holds 
  \begin{equation} \label{delta0Ws} |\lambda| \leq  \frac{c^2}{4}   + \frac{ \sup\{|g_u(u_f,w_f)|\}}{2\sup\{|f_u(u_f(s),w_f(s))|\} }  +\frac{1-\alpha}{\eta+\alpha} .
\end{equation}
\end{Lemma}

Taking  limits $\epsilon\to 0$ and $\delta \to 0$,  while $\epsilon\ll \delta$, in  an appropriate scaling of the independent variable,  the system can be reduced the Fisher-KPP scalar equation
\eqref{reductionKPP}, as shown in \cite{CGM}. In that case, the eigenvalue problem arising from the linearization about the front solutions takes the form
\eq{
\lambda W&=W_{\zeta\zeta}+cW_{\zeta}+\tilde{f}_w(w_f)W ,
}{KPPlin}
where
$$
\tilde{f}(w)\equiv  \frac{w\left(\sqrt{1-w}-\alpha\right)}{\eta+\sqrt{1-w}}.
$$
It follows from results obtained by Sattinger (see \cite[Theorem 6.3]{Sattinger76} and also the discussion in \cite[Section 2.2]{Xin00}) that the spectrum  of the operator defined by the right hand side of \eqref{KPPlin} is restricted to the left side of the complex plane $\realpart{\lambda}<0$, if considered on $L^2(\R)\cap L^2(\R,e^{-c\zeta/2})$. By $L^2(\R,e^{-c\zeta/2})$, we mean $L^2(\R)$ equipped with the weight  $e^{-c\zeta/2}$. We nevertheless present below how the point spectrum can be bounded on $L^2(\R)$. We also  present in Section \ref{EssSect} how the essential spectrum can be moved to the left side of the complex plane using an exponential weight $e^{\sigma \zeta}$, with $\sigma$ satisfying \eqref{KPPE}. Note that $\sigma=-c/2$ does fall in the range defined by \eqref{KPPE}.

\begin{Lemma} \label{L:4}  Consider  the eigenvalue problem \eqref{KPPlin}. For any eigenvalue $\lambda$   with $\rm{Re}(\lambda)>0$, the following estimate holds 
 \begin{equation} \label{delta00W} |\lambda| \leq  \frac{c^2}{4}   + \sup\{|\tilde{f}_w(w_f)|\}.
\end{equation}
\end{Lemma}
\begin{Proof}
The estimate \eqref{delta00W} is obtained in a very similar way as the previous ones. We only provide a general description of the process, leaving the details to the reader.
We first multiple  \eqref{KPPlin} $\overline{W}$ and integrate over $\R$. Then, one
takes the real and imaginary parts of the resulting equation. Inequalities are obtained by taking the absolute values of both sides of those equations and using the triangular inequality. We apply Young's inequality  as in the second part of \eqref{Y1}.
We add the two obtained inequalities and use the fact that $|\lambda|\leq \realpart{\lambda}+|\impart{\lambda}|$ to obtain \eqref{delta00W}.
\end{Proof}

  \section{Uniform in $\delta$ and $\epsilon$  bounds on the unstable eigenvalues \label{uniform}}

 The estimates \eqref{reim} obtained for the eigenvalue problem \eqref{evdelta} may be used for finding a specific area in the complex plane outside of which there are no unstable eigenvalues. We need to check for the presence of the eigenvalues numerically. Since we assume that $0< \epsilon\ll\delta$, we notice, from the inequality \eqref{reim},  that the bound on $|\lambda|$ contains terms of order  $O(1/\epsilon)$ and $O(1/\delta)$.   Decreasing either of these parameters should  cause an increase in the bound. On the other hand, the estimate \eqref{delta0Ws} obtained for the eigenvalue problem \eqref{evdelta0} which is  the system \eqref{evdelta} with $\epsilon=0$ implicitly depends on $\delta$. 
 Below we present  analytical arguments  showing that   there exists  a bound on the absolute value of the unstable eigenvalues of  \eqref{evdelta}  which is actually uniform in $\epsilon$ and $\delta$ as well as a bound on the unstable eigenvalues of  \eqref{evdelta0}   which  is uniform in $\delta$.

 Our first goal is to show that in the case of sufficiently small $0<\epsilon\ll\delta$, there are no unstable eigenvalues of \eqref{evdelta} outside of some semicircle with a radius  independent of parameters $\epsilon$ and $\delta$ provided these are chosen sufficiently small. This is done in Proposition \ref{Prop1}. Our second goal is to consider the case $\epsilon=0$ and  prove that in the case of sufficiently small $\delta>0$, there are no unstable eigenvalues of \eqref{evdelta0} outside of some semicircle with a radius independent of the parameter  $\delta$. This is done in Proposition \ref{Prop2}.
 
 \begin{Proposition}\label{Prop1} There exist $\epsilon_0$,  $\delta_0>0$ and  a constant  $r_1=r_1(\delta_0, \epsilon_0)$ so that for  all  $0<\epsilon\le \epsilon_0$  and $0<\delta \leq \delta_0$,  the eigenvalue problem \eqref{evdelta} has no nontrivial $L^2(\R)$-solutions for  values  $\lambda$ with $|\lambda|>r_1$ and $\argumc{\lambda}\in[-\pi/2,\pi/2]$. \end{Proposition}
To prove this claim  we rewrite \eqref{evdelta}, that is, the system
\eq{
\lambda U&=\epsilon  U_{\zeta\zeta} +cU_{\zeta}+\frac{1}{\delta }f_u(u_f,w_f)U + \frac{1}{\delta }f_w(u_f,w_f)W,\\
 \lambda  W&= W_{\zeta\zeta} +cW_{\zeta}+ g_u(u_f,w_f)U + g_w(u_f,w_f)W,
}{evdeltas}
in terms of a new variable $y=\zeta |\lambda|^{1/2}$
\eq{
\lambda U&=\epsilon  U_{yy} |\lambda|  +cU_{y} |\lambda|^{1/2}+\frac{1}{\delta }f_u(u_f,w_f)U + \frac{1}{\delta }f_w(u_f,w_f)W,\\
 \lambda  W&= W_{yy}|\lambda| +cW_{\zeta}|\lambda|^{1/2}+ g_u(u_f,w_f)U + g_w(u_f,w_f)W.
}{evdeltay}

The key idea  in the proof of this proposition is  the persistence   of exponential dichotomies under small perturbations as described in  \cite[Theorem 3.1]{Sandstede} or the classical source \cite [Chapter 4, Prop 1, page 34]{Coppel}.  The plan is to show that: (i) the first order system corresponding to \eqref{evdeltay} is a small perturbation of some limiting system, and that (ii) the limiting system has exponential dichotomies.

Let us rewrite \eqref{evdeltay} as a system of the first order equations as follows
\eq{
  \frac{dU_1}{dy}&=U_2,\\
\epsilon  \frac{dU_2}{dy}& = \left(e^{i \argum{\lambda}}  - \frac{1}{  \delta |\lambda| } f_u(u_f,w_f) \right)U_1 -\frac{c}{ |\lambda|^{1/2}} U_{2}-  \frac{1}{  \delta |\lambda| }f_w(u_f,w_f)W_1,\\
 \frac{dW_1}{dy}&=W_2,\\
 \frac{dW_2}{dy}&=  - \frac{1}{ |\lambda|} g_u(u_f,w_f)U_1+\left(e^{i \argum{\lambda}} - \frac{1}{ |\lambda|} g_w(u_f,w_f)\right)W_1-\frac{c}{ |\lambda|^{1/2}}  W_{2} .
}{y1}

We use the idea of \cite{AGJ} to freeze the parameters in the second equation in system \eqref{y1}.  We point out that the situation here is different from one considered in \cite{AGJ} because  \cite[Proposition 2.2]{AGJ} does not directly apply  since that equation in our case is not autonomous. 
We introduce a parameter $\gamma =(\epsilon \delta |\lambda|)^{-1}$, and rewrite \eqref{y1} as
\eq{
  \frac{dU_1}{dy}&=U_2,\\
  \frac{dU_2}{dy}& = \left(\frac{1}{\epsilon}e^{i \argum{\lambda}}  - \gamma f_u(u_f,w_f) \right)U_1-\frac{c\delta^{1/2} \gamma^{1/2} }{  \epsilon^{1/2}} U_{2} - \gamma f_w(u_f,w_f)W_1,\\
 \frac{dW_1}{dy}&=W_2,\\
 \frac{dW_2}{dy}&=   - \frac{1}{ |\lambda|} g_u(u_f,w_f)U_1+\left(e^{i \argum{\lambda}}- \frac{1}{ |\lambda|} g_w(u_f,w_f)\right)W_1-\frac{c}{ |\lambda|^{1/2}}  W_{2}.
}{y1eta}
Here, we abused the notations in the right hand side of the \eqref{y1eta} in the sense that, at the beginning,  $u_f$ and $w_f$  in \eqref{evdeltas} were, in fact, functions of the variable $\zeta$. In \eqref{evdeltay} their argument was changed to  $y / |\lambda|^{1/2}$, and so $u_f$ and $w_f$ became functions of the variable $y / |\lambda|^{1/2}$.
Then the parameter $\gamma$ was introduced, and,
 in terms of $\gamma$, the functions $u_f$ and $w_f$ in \eqref{y1eta} depend on the variable  $(\gamma \epsilon\delta )^{1/2}y $. However, the sup-norm of the functions $u_f$ and $w_f$ is independent of the pre-factor $(\gamma \epsilon\delta )^{1/2}$. We stress that passing to the autonomous asymptotic system in \eqref{y1eta} by taking the limit as $y\to\pm\infty$ should be done carefully as the limit is not uniform in $\gamma$. For instance, the norm of the difference of the coefficients of \eqref{y1eta} and the respective asymptotic systems  will become small for all $|y|\ge L(\gamma)$ with some $L(\gamma)>0$ that depends on $\gamma$. As we will see, this fact however will not affect our conclusions.

 We now  take in \eqref{y1eta} the limit as $|\lambda|\to \infty$, independently of the rest of the parameters, including  $\gamma$ and $\argumc{\lambda}$, and get
\eq{
  \frac{dU_1}{dy}&=U_2,\\
  \frac{dU_2}{dy}& =  \left(\frac{1}{\epsilon}e^{i \argum{\lambda}}- \gamma f_u(u_f,w_f)\right )U_1 -\frac{c\delta^{1/2} \gamma^{1/2 } }{ \epsilon^{1/2}} U_2- \gamma f_w(u_f,w_f)W_1,\\
 \frac{dW_1}{dy}&=W_2,\\
 \frac{dW_2}{dy}&=
 e^{i \argum{\lambda}}W_1.
}{y1etalim}

Let us denote   the coefficient matrix on the right hand side of   \eqref{y1etalim}  by $\mathcal{A}(y)$,   and  the coefficient associated with the right hand side of \eqref{y1eta} by $\mathcal{A}(y)+\mathcal{B}(y)$. More precisely, 
\begin{equation}\label{B}
\mathcal{B}(y)= \begin{pmatrix} 0&0&0&0\\0&0&0&0\\0&0&0&0\\
 - \frac{1}{ |\lambda|} g_u(u_f,w_f)&0&
  - \frac{1}{ |\lambda|} 
  g_w(u_f,w_f)) & -\frac{c}{ |\lambda|^{1/2}} 
\end{pmatrix}.
\end{equation}
The formulas  for $g_u$ and $g_w$ are given in \eqref{fderiv}. It is  easy to see that Lemma~\ref{L:de} implies that for sufficiently small  $\delta$ and $\epsilon>0$, there exists bounds on  $\|g_u\|_\infty$ and $\|g_w\|_\infty$ uniform in   $\delta$ and $\epsilon$.  
Therefore,
 \begin{equation}\label{Bnorm}
\|\mathcal{B}(y)\|=\mathcal{O}\left(\frac{1}{|\lambda|^{1/2}}\right) \text{ as $|\lambda|\to\infty$,}
\end{equation}
 that is, the system  \eqref{y1eta} is a perturbation of   \eqref{y1etalim} by terms that are small for $|\lambda|\to \infty$ uniformly with respect to $y$, $\epsilon$, $\delta$ and $\gamma$.

Next we want to prove that \eqref{y1etalim}  has exponential dichotomies on $\mathbb R^+$ and $\mathbb R^-$ with the rates of decay and growth of the solutions at infinities which is uniform in $\epsilon$, $\delta$ and $\gamma$. We formulate this statement as a lemma. 

\begin{Lemma} \label{L2} There exist $\epsilon_0$ and  $\delta_0>0$ such that for  any $\argumc{\lambda}\in[-\pi/2,\pi/2]$, $0<\epsilon \leq \epsilon_0$, $0<\delta\leq \delta_0$ and $\gamma >0$ the spectra of the matrices  $\mathcal{A}_{\pm}:=\lim_{y\to\pm\infty} \mathcal{A}(y) $  are separated from the imaginary axis uniformly in  $\argumc{\lambda}$, $\epsilon$, $\delta$, $\gamma$. As a result, system \eqref{y1etalim} has exponential dichotomies on $\R^-$ and on $\R^+$ with  the rates of decay and growth of the solutions at infinity uniform in $\argumc{\lambda}$, $\epsilon$, $\delta$, $\gamma$.
\end{Lemma}
\begin{Proof}
First, we consider system  \eqref{y1etalim} at the limit $y \to +\infty$ as we are interested in the solutions to \eqref{y1etalim} that are exponentially decaying and growing at $+\infty$. Using that $f_u(u_f,w_f) \to f_u(0,1) =-1$ and $f_w(u_f,w_f) \to f_w(0,1) =-1/2$, we have that under the limit $y \to +\infty$, system  \eqref{y1etalim} becomes
\eq{  \frac{dU_1}{dy}&=U_2,\\
  \frac{dU_2}{dy}& = \left(\frac{1}{\epsilon}e^{i \argum{\lambda}}  + \gamma \right)U_1 -\frac{c\delta^{1/2} \gamma^{1/2 } }{ \epsilon^{1/2}} U_{2}+ \frac{\gamma}{ 2 }W_1,\\
 \frac{dW_1}{dy}&=W_2,\\
 \frac{dW_2}{dy}&=
 e^{i \argum{\lambda}}W_1.
}{y1etalim0}
The spatial eigenvalues $\kappa$ of  this system are given by the solutions of the equation
\begin{equation}\label{m1}
\det \begin{pmatrix}
-\kappa&1&0&0\\
 \frac{1}{\epsilon}e^{i \argum{\lambda}}  +\gamma &-\frac{c\delta^{1/2} \gamma^{1/2 }}{  \epsilon^{1/2}}-\kappa &\frac{\gamma}{ 2 }&0\\
 0&0&-\kappa&1\\
 0&0&e^{i \argum{\lambda}} &-\kappa
\end{pmatrix}=0.
\end{equation}
The spatial eigenvalues are found to be  
\begin{equation}\label{97bis}
\kappa=\pm e^{i \argum{\lambda}/2}\end{equation} and the solutions of the equation
\begin{equation} \kappa^2 + \frac{c\delta^{1/2} \gamma^{1/2 }}{  \epsilon^{1/2}} \kappa -  \frac{1}{\epsilon}(e^{i \argum{\lambda}}  + \gamma)=0.
\end{equation}
The latter are
\begin{equation}\label{eqn71}
\kappa_{\pm}=\frac{1}{  2\epsilon^{1/2}}\left(-c\delta^{1/2} \gamma^{1/2} \pm  \sqrt{c^2\delta \gamma  + 4  (e^{i \argum{\lambda}}  + \gamma)}\right).
\end{equation}
Since $\epsilon$ in our model is a small parameter,  it is easy to see  that the $|\realpart{ \kappa_{\pm}}|$ have  positive lower bounds uniform in 
 $\epsilon\in (0,\epsilon_0]$ for all small $\epsilon_0$. 
 We next show that the real parts of $\kappa_{\pm}$ are bounded away from zero uniformly in  $\delta$  as long as it is smaller than some carefully chosen  $\delta_0$, for all $\gamma>0$ without any assumption on  the smallness of $\gamma$.  For the  real parts of $\kappa_{\pm}$ we have 
 \begin{equation}\epsilon^{1/2}\realpart{ \kappa_{\pm}}= -\frac{c\delta^{1/2} \gamma^{1/2} }{2} {\pm}\realpart{\sqrt{\left(\frac{c\delta^{1/2} \gamma^{1/2} }{2 }\right)^{2} +e^{i \argum{\lambda}}  + \gamma}}.\label{eqn72}
 \end{equation}
 Let us  introduce notation $\kappa_\pm(\gamma)$ for $\kappa_\pm$ in \eqref{eqn72} so that
  \begin{equation}
  \epsilon^{1/2} \realpart{ \kappa_{\pm}(\gamma)}
  = - \frac{1}{2} c\delta^{1/2} \gamma^{1/2}  {\pm} \left(\left(\frac{c^2\delta \gamma }{4 } + \gamma+\cos(\argumc{\lambda}) \right)^2 + \sin^2 (\argumc{\lambda}) \right)^{1/4} \cos(\phi),\label{eqn73}\end{equation} 
  where 
  \begin{equation}\phi (\gamma)= \frac{1}{2} \arctan\left(\frac{\sin (\argumc{\lambda})}{\gamma+\cos(\argumc{\lambda})+\frac{c^2\delta \gamma }{4 }   }\right).
  \end{equation}
 Here and in what follows we assume, without loss of generality, that  $\lambda$ is located in the first quadrant, i.e. $ 0\leq \argum{\lambda}\leq \pi/2$, so that both $\sin(\argumc{\lambda})$ and $\cos(\argumc{\lambda})$ are non-negative. We can make this assumption because the eigenvalue problem \eqref{evdeltas} has the symmetry that if $\lambda=\lambda_0$ is an eigenvalue, then so is $\lambda=\overline{\lambda}_0$. We  will show that $\cos \phi (\gamma)$ is  bounded away from $0$ uniformly in $\delta$ and $\gamma$.  Indeed, 
   \begin{equation}0\leq \frac {\sin (\argumc{\lambda})}{\gamma+\cos(\argumc{\lambda})+\frac{c^2\delta \gamma}{4 }  } \leq \frac{\sin (\argumc{\lambda})}{\cos(\argumc{\lambda}) }. 
  \end{equation}
By monotonicity then
  \begin{equation} 0\leq \phi (\gamma)= \frac{1}{2} \arctan\left(\frac{\sin (\argumc{\lambda})}{\gamma+\cos(\argumc{\lambda})+\frac{c^2\delta \gamma }{4 }   } \right)\leq \frac{1}{2} \arctan\left(\frac{\sin (\argumc{\lambda})}{\cos(\argumc{\lambda}) } \right) \leq \frac{1}{2}\argumc{\lambda} \leq \frac{\pi}{4}.
  \end{equation}
  Therefore   $ \sqrt{2}/2\leq  \cos (\phi (\gamma))   \leq 1$ for all $\delta >0$ and $\gamma >0$.
  
We focus on $\kappa_{+}$ first. From  \eqref{eqn73} it follows that
    \begin{equation}\label{77} \begin{split}
    \epsilon^{1/2} \realpart{\kappa_{+} (\gamma)} 
 & \geq  -\frac{c\delta^{1/2} \gamma^{1/2} }{2}   + \left(\left( \gamma \left(1+\frac{c^2\delta}{4}\right)+\cos(\argumc{\lambda}) \right)^2 + \sin^2 (\argumc{\lambda}) \right)^{1/4} \frac{\sqrt{2}}{2} \\
 &\geq -\frac{c\delta^{1/2} \gamma^{1/2} }{2}   +  \left(\gamma^2 \left(1+\frac{c^2\delta}{4}\right)^2+1  \right)^{1/4} \frac{\sqrt{2}}{2}.\end{split}
\end{equation}
We view the last bound as a function of $w=\gamma^{1/2}$ defined by
 $$b(\delta, w)= -\frac{c\delta^{1/2} }{2}  w+ \frac{\sqrt{2}}{2}\left( \left(1+\frac{c^2\delta}{4}\right)^2 w^4  +1\right)^{1/4}.$$
The function $w\mapsto b(\delta,w)$ has a critical point $w_0=w_0(\delta)$ that satisfies 
the equation

$$   \left( \frac{\sqrt{2}}{c\delta^{1/2} } \left(1+\frac{c^2\delta}{4} \right) \right)^{4/3} - \left(1+\frac{c^2\delta}{4}\right)^2w^4_0  = 1.  $$
The latter  equation  has a unique  positive solution  $w_0=w_0(\delta)$ if $0<\delta\leq \delta_0$  with $\delta_0>0$ so small that $c\delta_0^{1/2}(1+c^2\delta_0/4)^{2/3}<\sqrt2$. Moreover, 
it is easy to see that $w_0(\delta) \to 0 $ when $\delta \to 0$, but $b(0,0)=\sqrt{2}/{2}$,  so we can choose $\delta_0$ so small that $b(\delta, w_0(\delta)) \geq \sqrt{2}/4$ provided $\delta\in(0,\delta_0]$. 
Thus, the right hand side of  \eqref{77}  is greater than or equal to $\sqrt{2}/4$ for all $\gamma >0$ and $\delta\in (0,\delta_0]$.  
A simpler analogous argument shows that $ \epsilon^{1/2} \realpart{\kappa_{-} (\gamma) }< -1$    for all $\delta$ and  $\gamma$. 

 Thus, there is a small $\delta_0>0$ and $\epsilon_0>0$ such that  system \eqref{y1etalim0}  has two eigenvalues with positive and two with negative real parts. The  eigenvalues are bounded away from the imaginary axis uniformly in $\argumc{\lambda}\in[-\pi/2,\pi/2]$, 
 $\epsilon\in(0,\epsilon_0]$, $\delta\in(0,\delta_0]$ and $\gamma>0$. Since \eqref{y1etalim0} is the limiting system for \eqref{y1etalim},  system \eqref{y1etalim} has exponential dichotomy on $\R^+$ with the rates of decay and growth of the solutions at infinity  which are uniform in $\argumc{\lambda}$, $\epsilon$, $\delta$ and $\gamma$ provided $\argumc{\lambda}\in[-\pi/2,\pi/2]$, $\epsilon \in (0,\epsilon_0]$, $\delta \in (0, \delta_0]$ and $\gamma >0$. Indeed, we stress that the $y$-dependent coefficients in \eqref{y1etalim} actually depend on the variable $(\gamma\epsilon\delta)^{1/2}y$. As a result, the difference of $\mathcal{A}(y)$ and $\mathcal{A}_+$ is small for all $y>L(\gamma)$ for some large constant $L(\gamma)>0$, and by roughness of exponential dichotomy on the half-line $[L(\gamma),\infty)$, equation \eqref{y1etalim} will have exponential dichotomy on $[L(\gamma),\infty)$ because the respective asymptotic equation with the coefficient $\mathcal{A}_+$ does have exponential dichotomy on that half-line. The rate of exponential growth and decay of the solutions of  \eqref{y1etalim}  is controlled by the rate of exponential growth and decay of the solutions of the asymptotic equation, and so the rate does {\em not} depend on $\gamma$. As a result, equation \eqref{y1etalim} will have exponential dichotomy on $\mathbb{R}^+$. 
 While the constants that control  the actual growth or decay of the solutions of \eqref{y1etalim} on $\mathbb{R}^+$ {\em will} depend on $\gamma$, the exponential {\em rate } of growth and decay of the solutions is independent of $L(\gamma)$ and so of $\gamma$, as claimed in the lemma.

We will now discuss dichotomy on $\R^-$.
At the limit $y\to -\infty$,  
$$f_u(u_f,w_f) \to f_u(\alpha,1-\alpha^2 ) =-\frac{2\alpha^2}{1+\alpha}$$ and 
$$f_w(u_f,w_f) \to f_w(\alpha,1-\alpha^2) =-\frac{\alpha}{1+\alpha},$$ 
so the  limit of the  system  \eqref{y1etalim}  reads 
\eq{
  \frac{dU_1}{dy}&=U_2,\\
  \frac{dU_2}{dy}& = \left(\frac{1}{\epsilon}e^{i \argum{\lambda}}  +\gamma \frac{2\alpha^2}{1+\alpha}\right)U_1-\frac{c\delta^{1/2} \gamma^{1/2} }{  \epsilon^{1/2}} U_{2} + \gamma \frac{\alpha}{1+\alpha} W_1,\\
 \frac{dW_1}{dy}&=W_2,\\
 \frac{dW_2}{dy}&=
 e^{i \argum{\lambda}}W_1.
}{y1etalim12}
The spatial eigenvalues $\mu$ of  this system are given by the solutions of
\begin{equation}\label{m11}
\det \begin{pmatrix}
-\mu&1&0&0\\
 \frac{1}{\epsilon}e^{i \argum{\lambda}}  + \gamma\frac{2\alpha^2}{1+\alpha}&-\frac{c\delta^{1/2} \gamma^{1/2 }}{  \epsilon^{1/2}}-\mu &\gamma \frac{2\alpha^2}{1+\alpha} &0\\
 0&0&-\mu&1\\
 0&0&e^{i \argum{\lambda}} &-\mu
\end{pmatrix}=0.
\end{equation}
The eigenvalues are  $\mu=\pm e^{i \argum{\lambda}/2}$ and the solutions of the equation
\begin{equation}\mu^2 + \frac{c\delta^{1/2} \gamma^{1/2} }{  \epsilon^{1/2}} \mu -  \frac{1}{\epsilon}e^{i \argum{\lambda}}  - \gamma\frac{2\alpha^2}{1+\alpha}=0,
\end{equation}
which are
\begin{equation}
\mu_{\pm}= \frac{1}{\epsilon^{1/2}}\left(-\frac{c\delta^{1/2} \gamma^{1/2} }{  2}  \pm  \sqrt{\left(\frac{c\delta^{1/2} \gamma^{1/2} }{  2} \right)^2 +  e^{i \argum{\lambda}}  + \gamma\frac{2\alpha^2}{1+\alpha}}\right).
\end{equation}
Analysis similar to the case of $y\to +\infty$, with $\gamma$ replaced by $\frac{2\alpha^2}{1+\alpha}\gamma$,   shows that  there exist positive  $\delta_0$  such that for any $\delta <\delta_0$ and $\gamma>0$, 
the quantities $\epsilon^{1/2}\Re\mu_{\pm} $  are bounded away from zero uniformly in $\argumc{\lambda}$, $\delta $ and $\gamma$. System \eqref{y1etalim12}  has two eigenvalues with positive and two with negative real parts  which are  bounded away from the imaginary axis uniformly in $\epsilon$, $\delta$ and $\gamma$.   Therefore, the system \eqref{y1etalim} has exponential dichotomy at $-\infty$ and the rate of decay and growth of solutions  is uniform in $\epsilon$, $\delta$ and $\gamma$, as desired. 
\end{Proof}

We further notice that  in system \eqref{y1etalim} the term $f_u(u_f,w_f)$  is strictly negative for all $y\in\R$. As we will see in a moment, the negative sign guarantees the existence of the exponential dichotomy not only on both semi-axes but also on all of $\R$. Indeed, to show the existence of dichotomy on $\R$, it is enough to prove that the respective stable dichotomy subspaces coming from the dichotomies on the semi-axes have zero intersection at $y=0$. 
This is of course equivalent to the fact that system \eqref{y1etalim}  does not have a nontrivial bounded (equivalently, $L^2$-) solutions on $\R$. We show this fact next.
 
\begin{Lemma}\label{L3}The nonautonomous equation \eqref{y1etalim} has  no nontrivial $L^2(\R)$-solutions.
\end{Lemma}
\begin{Proof}
The system \eqref{y1etalim} has a 2-dimensional unstable manifold at $-\infty$ and a 2-dimensional stable manifold at $+\infty$. We next show that these manifolds do not intersect. To begin, we note that  the equations for  $W_1$ and $W_2$ in \eqref{y1etalim} are decoupled  from the equations for $U_1$ and $U_2$. The autonomous system for $W_1$ and $W_2$ has exponential dichotomy on $\R$ because $|\realpart{\kappa}|\ge\sqrt2/2$ for the real parts of the eigenvalues \eqref{97bis} of the respective matrix as $\argumc{\lambda}\in[-\pi/2,\pi/2]$. Thus, 
these autonomous  equations do not have  nontrivial  bounded or $L^2$-solutions. Furthermore, $U_1$ and $U_2$ then must satisfy the system
\eq{
  \frac{dU_1}{dy}&=U_2,\\
  \frac{dU_2}{dy}& = \left(\frac{1}{\epsilon}e^{i \argum{\lambda}}  - \gamma f_u(u_f,w_f) 
  \right)U_1-\frac{c\delta^{1/2} \gamma^{1/2} }{  \epsilon} U_{2}. 
}{y1etalim1233}
The $U_1$ component then satisfies the equation
\begin{equation}\label{y1etalim1234}
  \frac{d^2U_1}{dy^2}  +\frac{c\delta^{1/2} \gamma^{1/2} }{  \epsilon^{1/2}}\frac{dU_1}{dy}-\left(\frac{1}{\epsilon}e^{i \argum{\lambda}} -\gamma f_u(u_f,w_f)  \right)U_1 =0.  \end{equation}
Let us assume that there is  a nontrivial $U_1\in L^2(\R)$ that satisfies \eqref{y1etalim1234}. 
We multiply \eqref{y1etalim1234}  by $\overline{U}_1$ and use the fact that $U_1$ is in $L^2(\R)$ to perform the following integration
\begin{equation}\label{y1etalim0est}
 \intr  \frac{d^2U_1}{dy^2}\overline{U}_1 dy + \frac{c\delta^{1/2} \gamma^{1/2} }{  \epsilon^{1/2}}  \intr   \frac{dU_1}{dy}\overline{U}_1 dy -\intr \left(\frac{1}{\epsilon}e^{i \argum{\lambda}}  - \gamma f_u(u_f,w_f) \right) U_1\overline{U}_1dy=0.  \end{equation} 
The real part of this expression is 
 \begin{equation}\label{y1etalim0real}
 \intr  \frac{d^2U_1}{dy^2}\overline{U}_1dy  -\intr \left(\frac{1}{\epsilon}\cos (\argumc{\!\lambda}) - \gamma f_u(u_f,w_f) \right) U_1\overline{U}_1 dy=0.  \end{equation}
 We apply integration by parts to the first integral and obtain:
  \begin{equation}\label{y1etalim0real1}
- \intr  \left |\frac{dU_1}{dy}\right | ^2dy-\intr \left(\frac{1}{\epsilon}\cos (\argumc{\lambda})  - \gamma f_u(u_f,w_f) \right)| U_1|^2\ dy=0.  \end{equation}
For $\lambda$ in the closed half-plane one has \begin{equation}
\frac{1}{\epsilon}\cos (\argumc{\lambda})  - \gamma f_u(u_f,w_f)  >0,\end{equation}
due to Lemma~\ref{negativeBound2}. 
Therefore, \eqref{y1etalim0real1} may hold if only if $U_1=0$, but then $U_2 =0$ as well from \eqref{y1etalim1233} if $U_1=0$.  So,   the limiting system \eqref{y1etalim}   has no nontrivial bounded solutions. \end{Proof}

We are now ready to complete the proof of Proposition~\ref{Prop1}.  
As we have mentioned earlier, Lemmas~\ref{L2} and \ref{L3} imply that equation \eqref{y1etalim} has exponential dichotomy on $\mathbb R$. Moreover, Lemma ~\ref{L2} implies that there exist $\epsilon_0$ and $\delta_0$ such that for $0<\epsilon \leq \epsilon_0$, $0<\delta\leq \delta_0$ and for all $\gamma>0$,   the rates of exponential decay and growth of the solutions of \eqref{y1etalim} at $\pm\infty$ are uniform in $\epsilon$, $\delta$, and $\gamma$, as the rates are determined by the spectral gaps of the asymptotic matrices $\mathcal{A}_{\pm}$. We next choose $r_1=r_1(\epsilon_0,\delta_0)$
 such that for any $\lambda$ with $\realpart{\lambda} \geq 0$ and such that   $|\lambda|\geq r_1$ the $L_{\infty}$-norm of  the perturbation $\mathcal{B}(\cdot)$ defined by \eqref{B} is small compared to the spectral gap for  $\mathcal{A}_{\pm}$ for all $0<\epsilon \leq \epsilon_0$, $0<\delta\leq \delta_0$ and all $\gamma >0$. Then, by  \cite[Theorem 3.1]{Sandstede}, see also \cite{ChiconeLatushkin}, equation \eqref{y1eta}, which is a small perturbation of \eqref{y1etalim} for large $|\lambda|>r_1$  retains the exponential dichotomy on all of $\mathbb{R}$  for these values of the parameters. Therefore, equation   \eqref{y1} and thus \eqref{evdeltay} has no nontrivial $L^2(\R)$-solutions, as needed in Proposition~\ref{Prop1}.

This completes the analysis for $\epsilon>0$. Our next goal is to show that  when $\epsilon=0$,  and for sufficiently small $\delta>0$,  the eigenvalue problem \eqref{evdelta0} has no unstable eigenvalues outside of some semi-disk with radius  independent of  $\delta$. 
More precisely, we  prove the following statement.

\begin{Proposition}\label{Prop2} There exist  $\delta_0>0$ and  a constant  $r_1=r_1(\delta_0)$ so that for  all $0<\delta \leq \delta_0$  the eigenvalue problem \eqref{evdelta0} has no nontrivial $L^2(\R)$-solutions for  values  $\lambda$ with $|\lambda|>r_1$ and $\argumc{\lambda}\in[-\pi/2,\pi/2]$. \end{Proposition}

We fix $\delta_0>0$ sufficiently small as specified in the earlier parts of the paper, cf., e.g., Lemma \ref{L:d}.
 Let us introduce the  scalar operators 
\eq{  A_{11}&= A_{11}(\lambda)= c\frac{d}{d \zeta} +  \frac{1}{\delta }f_u(u_f,w_f) -\lambda, \\ A_{12}&=\frac{1}{\delta }f_w(u_f,w_f),\\ A_{21}&=g_u(u_f,w_f), \\  A_{22}&=A_{22}(\lambda) = \frac{d^2}{d  \zeta^2}+ c\frac{d}{d  \zeta} + g_w(u_f,w_f)-\lambda, 
 }{AA}
so that the eigenvalue problem  \eqref{evdelta0} can be rewritten as  
\begin{equation}\label{eqn93}
A\begin{pmatrix}U\\W\end{pmatrix} =0, \text{ where }  A= \begin{pmatrix}A_{11}&A_{12}\\A_{21}&A_{22}\end{pmatrix}.
\end{equation}

Note that  the estimate \eqref{delta0Ws} depends on $\delta\le\delta_0$ through the dependence on $\delta$ of $u_f$, $w_f$ and so the coefficient $g_w$ in $A_{22}$ depends on $\delta$ through the  component $u_f$  (it follows from \eqref{fderiv} that $g_w$ does not depend on $w_f$).
However, it follows from \eqref{fderiv} and \eqref{fderivlim} that the $L^\infty$-norm of $g_w$ has a bound which is uniform in $\delta\le\delta_0$  and depends on $\delta_0$.
The operator given by the expression  $ \frac{d^2}{d  \zeta^2}+ c\frac{d}{d  \zeta} + g_w(u_f,w_f)$ is sectorial \cite{Henry81}. The size of the sector for the spectrum of the operator and the estimate for its resolvent are uniform in $\delta\le\delta_0$ by the bound on $g_w$. So, there is then a $x_0=x_0(\delta_0)>0$ and $\theta=\theta(\delta_0)$ such that  the part of the spectrum of this operator in $\{\lambda\in\mathbb {C} :  \realpart{\lambda} \geq 0\}$ belongs to the triangle  
$$\Sigma_+=\left\{\lambda\in\mathbb {C}: \, \realpart{\lambda} \geq 0,  \,\, \realpart{\lambda} \leq x_0, \, | \mathrm{Im} \lambda |\leq \tan(\theta) (x_0-\realpart{\lambda} )\right\},$$ 
and  there exists a $C=C(\delta_0)>0$ such that if $\lambda \notin \Sigma_+$ and $ \realpart{\lambda}\ge0$,  then $A_{22}$ given in \eqref{AA} is invertible in $L^2(\mathbb R)$ with 
\begin{equation}\label{sect}
\|A_{22}^{-1}\|\leq \frac{C}{|\lambda|}.\end{equation}
For these values of $\lambda$, we may use the Schur complement  
\eq{B(\lambda) = A_{11}(\lambda)-A_{12}A_{22}^{-1}(\lambda)A_{21}}{Bdef}
and  write the matrix operator $A$  as
 \begin{equation}  A(\lambda)= \begin{pmatrix}A_{11}(\lambda)&A_{12}\\A_{21}&A_{22}(\lambda)\end{pmatrix} = \begin{pmatrix}B(\lambda)&A_{12}A_{22}^{-1}(\lambda)  \\0&I \end{pmatrix} \begin{pmatrix}I&0\\A_{21}&A_{22}(\lambda)\end{pmatrix} .
\end{equation}
 It follows that for $\lambda \notin \Sigma_+$ and $ \realpart{\lambda}\ge 0$ the operator $A(\lambda)$ is invertible in the space $L^2(\R)$ of $2$-dimensional vector valued functions if and only if  the operator $B(\lambda)$ is invertible in $L^2(\R)$.

Using this we conclude that the following implications hold:
(a) if $B(\lambda)$ is invertible then $A(\lambda)$ is invertible; 
and (b) if $A(\lambda)$ is invertible then 
the eigenvalue problem \eqref{eqn93} or, which is the same, the eigenvalue problem \eqref{evdelta0}, has no nontrivial solutions in $L^2(\R)$. Thus, to establish assertions required in Proposition \ref{Prop2}, it suffices 
to prove the following lemma.

\begin{Lemma}\label{lem14} There exists $\delta_0>0$ and a constant $r_1=r_1(\delta_0)$ 
such that for any $0<\delta\le\delta_0$ and  $\lambda$  with $\realpart{\lambda} \ge0$ and  $| \lambda|\ge r_1$
the operator  $B(\lambda)$ is invertible. \end{Lemma}
\begin{Proof}
Note that $A_{11} $ defined in \eqref{AA} is a scalar differential  operator of first order   with an asymptotically constant coefficient. It is invertible because, as proved earlier (see lemmas \ref{negativeBound1} and \ref{negativeBound2}), for $\realpart{\lambda}\ge0$ we have
\begin{equation}  \frac{1}{\delta }f_u(u_f(\pm\infty),w_f(\pm\infty))- \realpart{\lambda}<0,\end{equation}
 where $f_u(u_f(\pm\infty),w_f(\pm\infty))= \lim _{ \zeta\to\pm\infty}f_u(u_f( \zeta),w_f( \zeta))$.
  Then $B(\lambda)$ from \eqref{Bdef} can be rewritten as
\begin{equation}B(\lambda) =A_{11}(\lambda)(I-A_{11}^{-1} (\lambda)A_{12}A_{22}^{-1}(\lambda)A_{21}).
\end{equation}
To prove that $B(\lambda)$ is invertible, in the remaining part of the argument, we show that 
\begin{equation}\label{Bdecay}
A_{11}^{-1}(\lambda) A_{12}A_{22}^{-1}(\lambda)A_{21}\to 0 \text{ as } |\lambda|  \to \infty,\end{equation}
 uniformly in  $\delta$.

We first show that   the inverse of  $A_{11}(\lambda)$  is given by 
\begin{equation} \label{A11inv} v( \zeta)=(A_{11}^{-1} u)( \zeta)=-\frac{1}{c}\int_{ \zeta}^{+\infty} \exp\left({\int^{\tilde  \zeta}_ \zeta \frac{1}{c} \left( \frac{1}{\delta }f_u(u_f(s),w_f(s)) -\lambda\right)\mathrm d s}\right)u(\tilde  \zeta)\mathrm d \tilde{ \zeta}.
\end{equation}
This can be checked directly, by calculating $A_{11} v=c\frac{dv}{d  \zeta} +  \left(\frac{1}{\delta }f_u(u_f,w_f) -\lambda\right) v$, with $v$ given above. Indeed,  
\begin{eqnarray*}
c\frac{dv}{d  \zeta}  
&=&  u(  \zeta) + \left( \frac{1}{\delta }f_u(u_f( \zeta),w_f( \zeta)) -\lambda\right)  \frac{1}{c} \int_{ \zeta}^{+\infty} \exp\left({\int_{\tilde  \zeta}^ \zeta   \frac{1}{c}\left(  \frac{1}{\delta }f_u(u_f(s),w_f(s)) -\lambda\right)\mathrm d s}\right)u(\tilde  \zeta)\mathrm d \tilde{ \zeta} \\
&=& u(\zeta)- \left( \frac{1}{\delta }f_u(u_f( \zeta),w_f( \zeta)) -\lambda\right) v(\zeta).
\end{eqnarray*}
It follows that
$$
A_{11}v=c\frac{dv}{d  \zeta} +  \left(\frac{1}{\delta }f_u(u_f,w_f) -\lambda\right) v = u,$$
thus proving \eqref{A11inv}.

To simplify notation, we {{temporarily}} denote:
\eq{
\tilde{f}(s)&=f_u(u_f(s),w_f(s)),\quad \nu=-\max\{\tilde{f}(s): s\in\mathbb{R}\}, 
\\ {{D}}(s)&=\begin{cases} 0, & s>0,\\ \frac{1}{c}\exp\big(\frac{1}{c}(\frac{\nu}{\delta}+{\rm Re } \lambda)s\big), & s\le 0,\end{cases}}{101}

and remark that $\nu>0$ by Lemmas \ref{negativeBound1} and \ref{negativeBound2}. Using \eqref{A11inv}, we conclude:
\begin{eqnarray*}
\| A_{11}^{-1}u\|_{L^2}&=&\Big(\int_{-\infty}^\infty \left|\frac{1}{c}\int_ \zeta^\infty\exp\left(\int_ \zeta^{\tilde{ \zeta}}\frac{1}{c}\big(\frac{1}{\delta}\tilde{f}(s)-\lambda\big)ds\right)u(\tilde{ \zeta})\,d\tilde{ \zeta}\right|^2\, d \zeta\Big)^{1/2}\\
&\le&\Big(\int_{-\infty}^\infty \Big(\frac{1}{c}\int_ \zeta^\infty\exp\left(\frac{1}{c}\big(-\frac{\nu}{\delta}-\realpart{\lambda}\big)( \zeta-\tilde{ \zeta})\right)\big|u(\tilde{ \zeta})\big|\,d\tilde{ \zeta}\Big)^2\, d \zeta\Big)^{1/2}\\
&=&\Big(\int_{-\infty}^\infty \Big(\int_{-\infty}^\infty {{D}}( \zeta-\tilde{ \zeta})\big|u(\tilde{ \zeta})\big|\,d\tilde{ \zeta}\Big)^2\, d \zeta\Big)^{1/2}\\
&=&\| {{D}}*\big| u\big|\|_{L^2}\le\|{{D}}\|_{L^1}\|u\|_{L^2},
\end{eqnarray*}
where, in the last line, we used Young's inequality for convolutions.
Computing $\|{{D}}\|_{L^1}=\delta/(\nu+\delta \realpart{\lambda})$ and using 
$\realpart{\lambda}\ge0$, we arrive at the inequality $\|A_{11}^{-1}\|\le\delta/\nu$.
Therefore, there is a $\delta$-independent constant $\mathcal {C}$ such that
\begin{equation}\|A_{11}^{-1} A_{12}A_{22}^{-1}A_{21}\| \leq  \mathcal {C}\delta\cdot\frac{1}{\delta} \frac{1}{|\lambda|} =\frac{\mathcal {C}}{|\lambda|},
\end{equation}
where we also used the definitions of $A_{12}$ and $A_{21}$ in \eqref{AA} and \eqref{sect}, thus yielding \eqref{Bdecay}. This concludes the proof of Lemma~\ref{lem14} and therefore Proposition~\ref{Prop2}.
\end{Proof}

\bigskip

\section{Discrete spectrum}\label{DS}
{{The aim of this section is to relate the discrete spectrum to the spectrum associated to the KPP equation \eqref{reductionKPP}, complete the proof of Theorem \ref{T:2}, and conclude the asymptotic stability of the waves.}}
\subsection{Goals and strategy}
Our goal is to study how the discrete spectrum of \eqref{evdelta2} with a nonzero but small $\epsilon$ is related to the spectrum of   \eqref{evdelta2}  with $\epsilon=0$. 
The Evans function is a convenient tool to study the discrete eigenvalues but it does not survive taking a limit $\epsilon\to 0$, because it  effectively counts the  $\lambda$-dependent dimension of the manifold constructed using the unstable manifold of the trivial equilibrium of the equation of the linearization of the underlying PDE   about the traveling front  at $-\infty$ and the stable manifold of the trivial equilibrium of the equation of the linearization of the underlying PDE   about the traveling front  at $\infty$. Outside of the continuous spectrum,  the sum of the dimensions of these manifolds is equal to the dimension of the phase space. When $\lambda$ is an eigenvalue, these two manifolds intersect. 
In the limit $\epsilon \to 0$  at least one these manifolds collapses. The phase space for the eigenvalue problem with $\epsilon=0$ is strictly less than that of the problem with a nonzero $\epsilon$.  


To circumvent this issue, we use the concept of the unstable augmented bundles as introduced in \cite{AGJ} (see also \cite{GJ90}). The unstable augmented bundles are topological objects and as such they survive the limit operation. It is known that for the  multi-scale problems  the first Chern number, which is a characteristic of the bundle, is equal to the sum of the first Chern numbers  for the  bundles of the unstable augmented so-called ``fast'' and ``slow'' problems.  The Evans function plays an important role in the construction of the augmented unstable bundles. 

{ 
The geometric approach based on the construction of the unstable bundles was used in a number of papers. We refer the reader  to  \cite{RDR} and references therein.  In \cite{AGJ}, the authors applied the unstable bundle approach when investigating the stability of  traveling pulses in the FitzHugh–Nagumo system.   One of the interpretations of their result is that  the slow-fast structure in the  eigenvalue problem induces a factorization of the Evans function into a slow and a fast component. As each of the factors  is associated with an unstable subbundle 
 (slow or fast), the zeros of the factors are approximations  of the zeros of the  Evans functions in the limiting  slow and fast  systems. 
Similarly, the factorization of the  Evans function into a slow and a fast component may be concluded from the results of  \cite{GJPred}.  In this case, where the solutions are of a complex structure,  the analyses are again based on the unstable bundle decomposition on the slow and fast subbundles, where the elephant trunk lemma was used to develop geometric and topological arguments for establishing the count of eigenvalues.
 
 A completely different, analytic (as opposed to geometric)  approach was developed in \cite{RDR}. The key is to  use the Riccati transformation and exponential dichotomies to show that the Evans function may be approximated by a product of  Evans functions for the limiting slow
and  fast systems. In \cite{RDR}, this approach is used to study the  spectra of the linearization  about stationary, spatially periodic pulses in a class of reaction-diffusion systems.  The Riccati transformation leads to slow-fast factorization of the Evans function and provides a framework for passing  to the singular limit.  We note that the nonlinear stability is investigated in \cite{dR} for the same spatially periodic pulses.   
 
 Our system  shares some properties with the system considered in \cite{RDR}. For example, we were able to prove that the constants in the exponential dichotomies do not depend on the small parameters. However, we have chosen to use the geometric  technique, i.e., the unstable bundles to treat the spectrum of traveling fronts considered in this paper.}

Our strategy is as follows:
\begin{itemize}
\item For our analysis we need to know that in the complex plane there is a  simple closed contour that encircles all possible eigenvalues of \eqref{evdelta2} with a nonzero $\epsilon$.
There are two requirements for a such contour. Without loss of generality, we may imagine that it is  a boundary of a  right-hand side section of a circle with the center at the origin  cut slightly to the left of the imaginary axis.    
  We want the radius  of the circle not to expand  unboundedly to the right of the imaginary axis as $\epsilon\to 0$ and $\delta\to 0$.  In other words,  our analysis works if the radius of the circle  is uniformly bounded in $\epsilon$ and $\delta$ when  $\epsilon$  and $\delta$ are sufficiently small. Existence of such a contour follows from 
Proposition~\ref{Prop1}. 
 We also want the contour to slightly extend across and to the left  the imaginary axis to count for the unstable eigenvalues possibly embedded in the essential spectrum. We show in the next subsection  that this is possible to do despite the presence of the continuous spectrum.  
 We call a contour that satisfies these two conditions $K$.

\item Given a contour $K$, we can then  define Evans function as in \cite[Sect.~4]{AGJ}  and the augmented unstable bundle as in   \cite[Sect.~3]{AGJ}. The Evans function serves as ``intermediary'' between the count of the eigenvalues and the bundle, in the sense of the following  theorem, cf.\ \cite[Sect.~2]{AGJ}: {\it The following three numbers are equal:  (1) the winding number of the Evans function over the contour $K$, (2) the first Chern number of the augmented unstable bundle over $K$ and its interior, (3) the number of eigenvalues of  \eqref{evdelta2} inside contour $K$ counting their algebraic multiplicity.} We point out that if the analytic extension of the Evans function  exists then the augmented unstable bundle construction is still applicable and  the count of zeros of the Evans function  and the first Chern number of the augmented unstable bundle coincide, see \cite{AGJ, GJ90}. 
\item We then use multi-scale nature of the  eigenvalue  problem \eqref{evdelta2}  to show  that the for sufficiently small $\epsilon$ and, consequently, $\delta$ the eigenvalues of  \eqref{evdelta2}  are small perturbations of the eigenvalues of some limiting problem.  This is accomplished   by using the two important properties  of the bundle: 1) as a topological object it is robust under taking limits, including when the limit is singular; 2)  the first Chern number of a bundle  is additive  in the Whitney sum of bundles, in this particular case  the decomposition that we are interested in is related to the multi-scale nature of our system. The unstable bundle may be viewed as the Whitney sum of bundles associated to the fast and slow reduced systems of  \eqref{evdelta2}. We construct this decomposition using   the small parameters $\epsilon$ and $\delta$ consequently. We show that  the fast (with respect to both $\epsilon$ and $\delta$) bundles do not contribute to the count of the unstable eigenvalues, and therefore it is the  slow limiting system that is responsible for the stability of the wave. 
\item {{We apply the results from this and previous sections and also Sattinger's work \cite{Sattinger76} to complete the proof of Theorem \ref{T:2} and obtain asymptotic stability.}}  
\end{itemize}

\subsection{Definition of the Evans function} 
\label{EvansDef}
Recall that the Evans function is analytic outside of the essential spectrum with zeroes corresponding to the discrete eigenvalues of the operator of the linearization about the wave.  The plan for this  section is as follows.  We prepare our  eigenvalue problem for the bundle construction and build some machinery for the bundle's multiscale analysis   with respect to $\epsilon\ll\delta\ll1$. We deal with  these parameters sequentially, first with  $\epsilon$ and, later, with $\delta$. In this section we focus on $\epsilon$.   
  We then briefly recall how  Evans function $E(\lambda)$ was constructed in  \cite{AGJ}. We then argue that in our case  Evans function $E(\lambda)$  
  allows  an analytic  extension  $\tilde E(\lambda)$  across the boundary of the essential spectrum. 
  


To start, we  consider  the eigenvalue problem \eqref{evdelta2} as a dynamical system,
\eq{
  \frac{dU_1}{d\zeta}&=U_2,\\
 \frac{dU_2}{d\zeta}& =\frac{1}{\epsilon}\left( (\lambda - \frac{1}{\delta } f_u(u_f,w_f) )U_1 -cU_{2}-  \frac{1}{\delta }f_w(u_f,w_f)W_1\right),\\
 \frac{dW_1}{d\zeta}&=W_2,\\
 \frac{dW_2}{d\zeta}&=  -g_u(u_f,w_f)U_1+(\lambda - g_w(u_f,w_f))W_1-cW_{2} .
}{newslow}

We treat   \eqref{newslow} as a multi-scale system. We assume here that $0<\epsilon\ll 1$. 
 We refer to the original independent variable $\zeta$ as the slow variable and, therefore,   to system \eqref{newslow} as the slow system associated with  the eigenvalue problem \eqref{evdelta2} as opposed to the fast system that will be introduced below. 
We denote\footnote{Here and in what follows we suppress the dependence of $A$ from $\delta$} the coefficient matrix in  the slow system by \eqref{newslow}  as
\begin{equation} \label{Am} A_s(\zeta, \lambda, \epsilon) = \begin{pmatrix}
 0&1  &0&0\\
\frac{1}{\epsilon}(\lambda - \frac{1}{\delta } f_u(u_f,w_f) ) &-\frac{c}{\epsilon} &-  \frac{1}{\epsilon \delta }f_w(u_f,w_f)&0\\
 0&0& 0&1\\
-  g_u(u_f,w_f) &0&  \lambda - g_w(u_f,w_f)&- c
\end{pmatrix},\end{equation}
where the subscript $s$ stands for ``slow''.
  The slow system \eqref{newslow} then can be written as
  \eq{
\frac{d \vec{U}}{d\zeta}  = A_s(\zeta, \lambda, \epsilon)  \vec{U}.
}{evslowzeta}
We next  introduce a fast variable $\xi =\zeta/\epsilon$ and obtain the system that we call the {\em fast system} associated with \eqref{evdelta2},
\eq{
  \frac{dU_1}{d\xi}&=\epsilon U_2,\\
 \frac{dU_2}{d\xi}& = (\lambda - \frac{1}{\delta } f_u(u_f,w_f) )U_1 -cU_{2}-  \frac{1}{\delta }f_w(u_f,w_f)W_1,\\
 \frac{dW_1}{d\xi}&=\epsilon W_2,\\
 \frac{dW_2}{d\xi}&=  \epsilon \left(-g_u(u_f,w_f)U_1+(\lambda - g_w(u_f,w_f))W_1-cW_{2} \right).
}{fastzeta}
We point out that this system is written with an abuse of notation as we use the same notation for $U_{1,2}$ and $W_{1,2}$, as well as $u_f$ and $w_f$ regardless of whether the independent variable is $\xi$ or  $\zeta$.

We denote the  coefficient matrix in \eqref{fastzeta} by $A_f(\xi, \lambda, \epsilon)$. In this notation the subscript $f$ stands for ``fast''.  The slow and fast coefficient matrices are related by the formula
 \eq{ A_f(\xi, \lambda, \epsilon) =\epsilon A_s(\epsilon\xi, \lambda, \epsilon).}{fsrel}
 The fast  system \eqref{fastzeta} in  the matrix form is  given by
  \eq{
\frac{d \vec{U}}{d\xi}  = A_f(\xi, \lambda, \epsilon) \vec{U}.
}{evfastxi}

For $\lambda$ to be an eigenvalue of \eqref{evdelta2},  there should be a bounded (equivalently, $L^2(\R)$) solution to the eigenvalue problem. That implies that there should be a solution to \eqref{evfastxi} (or \eqref{evslowzeta}) that approaches the zero equilibrium of \eqref{evfastxi} (or \eqref{evslowzeta}) along the equilibrium's stable manifold at $+\infty$ and along the unstable manifold at $-\infty$.  So we consider the limits of  $A_f(\xi, \lambda, \epsilon)$  ($A_s(\zeta, \lambda, \epsilon)$) as $\xi$ (or $\zeta$) approaches $\pm\infty$.
When $\epsilon>0$, the coefficient matrices have limits when their corresponding variables approach $\pm\infty$. This is due to  the exponential convergence of  the components of  the traveling front to its rest states.  More specifically, these limits are 
\begin{equation} \label{Amf+} A_f^{+}  (\lambda, \epsilon) \equiv A_f(+\infty, \lambda, \epsilon) = \begin{pmatrix}
 0&\epsilon  &0&0\\
\lambda +\frac{1}{\delta }  &-c &  \frac{1}{2 \delta}&0\\
 0&0& 0&\epsilon\\
0 &0&  \epsilon(\lambda - \frac{1-\alpha}{1+\eta} )&- \epsilon c
\end{pmatrix},
\end{equation}
\begin{equation} \label{Amf-} A_f^{-}  (\lambda, \epsilon) \equiv A_f(-\infty, \lambda, \epsilon) = \begin{pmatrix}
 0&\epsilon  &0&0\\
\lambda +\frac{1}{\delta } \frac{2\alpha^2}{1+\alpha}  &-c &  \frac{1}{ \delta } \frac{\alpha}{1+\alpha}&0\\
 0&0& 0&\epsilon\\
-   \frac{\epsilon (1-\alpha^2)}{\eta+\alpha} &0&  \epsilon \lambda &- \epsilon c
\end{pmatrix},\end{equation}
and \eq{ A_s^{\pm}(\lambda, \epsilon)\equiv A_s(\pm\infty, \lambda, \epsilon) =\frac{1}{\epsilon} A_f^{\pm}( \lambda, \epsilon).}{fsrel2}

Due to the dynamics at $+\infty$,  when $\lambda$ is to the right of the rightmost boundary of the essential spectrum, The matrices $A_f^{\pm}( \lambda, \epsilon)$ have two eigenvalues with positive real parts and two eigenvalues with negative real parts. Indeed, the eigenvalues of $A_f^{+}  (\lambda, \epsilon)$ can be calculated in closed form; the eigenvalues with positive real parts are given by
\eq{\mu^+_{1}=\frac{1}{2}\left(-c+\sqrt{c^2+4\epsilon\left(\lambda+\frac{1}{\delta}\right)}\right),\quad  \mu^+_{3}=\frac{\epsilon}{2}\left(-c+\sqrt{c^2+4\left(\lambda-\frac{1-\alpha}{1+\eta}\right)}\right),}{APP}
 while the eigenvalues with negative real parts  are given by
\eq{\mu^+_{2}=\frac{1}{2}\left(-c-\sqrt{c^2+4\epsilon\left(\lambda+\frac{1}{\delta}\right)}\right),\quad  \mu^+_{4}=\frac{\epsilon}{2}\left(-c-\sqrt{c^2+4\left(\lambda-\frac{1-\alpha}{1+\eta}\right)}\right).}{APM}
The tangent subspace to the stable manifold at the equilibrium at $+\infty$ is spanned by the eigenvectors corresponding to $\mu^+_2$ and $\mu^+_4$.
We note that when $\realpart{\lambda} =-\frac{c^2}{4} +\frac{1-\alpha}{1+\eta}$ the two eigenvalues $\mu^+_3$ and $\mu^+_4$ coincide. Recall that we assume that   
$c > 2\sqrt{ \frac{1-\alpha}{1+\eta}}$ . 

The eigenvalues of $A_f^{-}  (\lambda, \epsilon)$  do not have convenient closed form expressions. However, 
the essential spectrum for the original linear system determined by the matrix $A_s$ given in \eqref{newslow} is studied in 
 Section~\ref{EssSect}. In that case, it is shown {{in Section \ref{S:Numerics}}} that to the right of the curves of the essential spectrum given by \eqref{curve-}, two of the  eigenvalues  have positive real parts and two of them have negative real part.  Both curves  \eqref{curve-}  are strictly to the left of the imaginary axis. Lemma \ref{5} states that the expression  \eqref{spgap} gives the gap between the right most curve and the imaginary axis  is independent of $\epsilon$.  Let  us denote by ${{S}}$ the length  of the  gap given in \eqref{spgap}. In view of the relation \eqref{fsrel} between $A_s$ and $A_f$, or more importantly here, in view of the relation \eqref{fsrel2} between the asymptotic matrices, to obtain the essential spectrum for the system specified by $A_f$,  one only needs to replace $k$ by $k/\epsilon$ in the expressions giving $\lambda$ as a function of $k$ obtained in Section~\ref{EssSect}. In particular 
 the gap for the system specified by $A_f$ does not change and is still given by ${{S}}$.
  It has a  nonzero  finite limit as $\epsilon\to 0$ and also as $\delta \to 0$.

Because of the  $2$--$2$ consistent splitting, following \cite{AGJ} we  consider \eqref{evfastxi} (or \eqref{evslowzeta}) in the framework of exterior powers 
$\Lambda ^2 (\mathbb {C}^4)$ (see \cite{[6], [20]}):
\eq{
\frac{d\hat{U}}{d\xi}= [A_f ]^{(2)}(\xi, \lambda, \epsilon)\hat{U}.}{lambda2fast}
Here $\left[A_f\right]^{(2)}$ is a $6 \times 6$ matrix 
generated by  $A_f$ on the wedge-space $\Lambda^2(\mathbb{C}^4)$. 
For every fixed $\lambda$, the asymptotic systems for  \eqref{lambda2fast} is given by
\eq{
\frac{d\hat{U} }{d\xi}= [A_f ^{\pm}]^{(2)}(\lambda, \epsilon)\hat{U}.}{lambda2fastasymp}
The matrices $[A_f ^{\pm}]^{(2)}$ have as their eigenvalues the sums of any two of the eigenvalues of $A_f ^{\pm}$.
 For example,  when $\lambda$ is to the right of the right most boundary of the essential spectrum, $A_f ^{+}$ has two eigenvalues $\mu^+_1$ and $\mu^+_3$ with positive real parts and two 
eigenvalues  $\mu^+_2$ and $\mu^+_4$ with  negative real parts.  Then  the eigenvalue  of  $[A_f ^{+}]^{(2)} $ with the ``most negative'' real part is $\mu_+(\lambda)= \mu^+_2 +\mu^+_4$. 
The eigenvalue  of  $[A_f ^{-}]^{(2)}$ with the largest positive real part  is  $\mu_-(\lambda)= \mu^-_1 +\mu^-_3$. 

As mentioned before, the asymptotic matrix $A_f^-$ has two of its eigenvalues with positive real parts and two of them with negative real part as long as $\realpart{\lambda}\geq -{{S}}$, where ${{S}}$ the length of the gap given in \eqref{spgap}, as stated in Lemma \ref{5}, between the boundary of the essential spectrum and the imaginary axis. The essential spectrum associated to the asymptotic matrix at $+\infty$, $A_f^+$, intersects with the right side of the complex place. However, it follows from the explicit expressions of the eigenvalues given in \eqref{APP} and \eqref{APM}, that the pair of eigenvalues given in \eqref{APM} remain the pair of eigenvalues of $A_f^+$ whose sum has the smallest real part if $\realpart{\lambda}\geq -\frac{c^2}{4}+\frac{1-\alpha}{1+\eta}$. Note that the RHS of that inequality is a negative expression given the condition on $c$ given in Theorem \ref{T:1}. 
Both eigenvalues $\mu_+(\lambda)$ and $\mu_-(\lambda)$ thus are well-defined and simple if $\realpart{\lambda}$ satisfies the two inequalities just mentioned. 
  As the eigenvalues are simple, the corresponding eigenvectors are analytic in $\lambda$. The Evans function $E(\lambda)$ now can then be constructed as in \cite[Sect.~4]{AGJ}
in the region $ \realpart{\lambda} >- \rho$, where  $\rho< \min\{\frac{c^2}{4}-\frac{1-\alpha}{1+\eta}, {{S}}\}$.


We now recall the definition of the Evans function that stems from the compactified version of the eigenvalue problem.  Indeed,
following \cite{AGJ},  the eigenvalue problem \eqref{newslow} is compactified by introducing  a new variable  $\tau$ as in $\zeta= \frac{1}{2k}\ln \left(\frac{ 1+\tau}{1-\tau}\right)$, where $k$ is a positive constant.  
We replace $\zeta$  by $\tau(\zeta)$ only in the expressions $(u_f,w_f)$, and keep the same notation for the unknown functions as functions of $\zeta$, thus obtaining from \eqref{newslow} the following autonomous nonlinear system,
\eq{
  \frac{dU_1}{d\zeta}&=U_2,\\
 \frac{dU_2}{d\zeta}& = \frac{1}{\epsilon}\left((\lambda - \frac{1}{\delta } f_u(u_f,w_f) )U_1 -cU_{2}-  \frac{1}{\delta }f_w(u_f,w_f)W_1\right),\\
 \frac{dW_1}{d\zeta}&=W_2,\\
 \frac{dW_2}{d\zeta}&=   -g_u(u_f,w_f)U_1+(\lambda - g_w(u_f,w_f))W_1-cW_{2},\\
  \frac{d \tau}{d\zeta} &=k(1-\tau^2).
}{slowdtau}
Still suppressing in the notations dependence on $\delta$, we denote 
\begin{equation}
A_s(\tau, \lambda, \epsilon) =\begin{cases} \,\,\,\,A_s(\zeta(\tau), \lambda, \epsilon), \, \text{ for } \,\,\tau \neq \pm 1, \\
\lim\limits_{\zeta \to \pm \infty}A_s(\zeta, \lambda, \epsilon), \, \text{ for } \,\, \tau  = \pm 1. \end{cases} \end{equation}
The system 
\eq{
\frac{d \vec{U}}{d\zeta}  &= A_s(\tau, \lambda, \epsilon) \vec{U},\\
  \frac{d \tau}{d\zeta} &=  k(1-\tau^2),
}{evatau}
where $ \vec{U} =(U_1,U_2,W_1, W_2)^T$, is  equivalent to \eqref{slowdtau}.

We consider   \eqref{slowdtau} as a multi-scale system. More precisely, we introduce a fast variable $\xi =\zeta/\epsilon$,  
\eq{
  \frac{dU_1}{d\xi}&=\epsilon U_2,\\
 \frac{dU_2}{d\xi}& = (\lambda - \frac{1}{\delta } f_u(u_f,w_f) )U_1 -cU_{2}-  \frac{1}{\delta }f_w(u_f,w_f)W_1,\\
 \frac{dW_1}{d\xi}&=\epsilon W_2,\\
 \frac{dW_2}{d\xi}&=  \epsilon \left(-g_u(u_f,w_f)U_1+(\lambda - g_w(u_f,w_f))W_1-cW_{2} \right),\\
 \frac{d \tau}{d\xi}  &= \epsilon k(1-\tau^2).
}{fastdtau} 
Recall that, with the same as above abuse of notation, the components of the  front $u_f$ and $w_f$  are  functions of $\tau$. In   the   matrix notation, this system  reads
\eq{
\frac{d \vec{U}}{d\xi}  &=\epsilon  A_s(\tau(\epsilon\xi), \lambda, \epsilon)\vec{U}, \\
\frac{d \tau}{d\xi} &= \epsilon k(1-\tau^2).}{evataufast}

The exponential  rates of convergence of the front to its rest states  have well defined limits bounded away from $0$ when  $\epsilon\to 0$ (as well as  when $\delta \to 0$). 
It is proven in  \cite[Lemma 3.1]{AGJ}  that for a sufficiently small $\epsilon$ it is possible to choose a small $k$ so that the right hand sides of \eqref{slowdtau} and  \eqref{fastdtau}  are $C^1$ on $\mathbb{C}^4 \times [-1,1]$.  For any $\tau\neq 1$ the smoothness is clear. The cases of $\tau=\pm 1$ are treated in \cite[Lemma 3.1]{AGJ} in for general case of semilinear parabolic systems. It is shown \cite[Lemma 3.1]{AGJ} that it is enough to take  $k$ less than a half of the smallest of the rates of convergence of the front to its rest states to guarantee smoothness at the endpoints $\tau=\pm 1$.  

The companion equations for \eqref{evatau} in $\Lambda^2(\C^4)$ are given by 
\eq{
\frac{d\hat{U}}{d\zeta}& = \left([A_s(\tau, \lambda, \epsilon)]^{(2)} - \mu_{\pm}(\lambda)I_{6\times6} \right) \hat{U},\\ 
  \frac{d \tau}{d\zeta} &=  k(1-\tau^2).
}{evataucent}
The RHS of the systems \eqref{evataucent} are also  $C^1$ on $\Lambda^2 (\mathbb{C}^4) \times [-1,1]$.  
The system \eq{
\frac{d\hat{U}}{d\zeta}& = \left([A_s(\tau, \lambda, \epsilon)]^{(2)} - \mu_{-}(\lambda)I_{6\times6} \right) \hat{U},\\ 
  \frac{d \tau}{d\zeta} &=  k(1-\tau^2)
}{evataucent_-}
has an invariant sets $\tau=- 1$, on which  zero  is the eigenvalue of $[A_s^-(\lambda,\epsilon)]^{(2)}- \mu_{-}(\lambda)I_{6\times6}$ with  the largest real part.  The associated  one-dimensional central manifold consists of the solutions $\hat{U}$ to 
\eq{\left([A_s ^{-}(\lambda, \epsilon)]^{(2)} - \mu_{-}(\lambda)I \right) \hat{U}=0.}{wc-}
In other words, the central manifold at $\tau=-1$ is the linear span of the eigenvector  $z_{-}(\lambda)$ corresponding to the eigenvalue $\mu_{-}(\lambda)$ of the matrix $[A_s^-(\lambda,\epsilon)]^{(2)}$.   The equilibrium $(z_{-}(\lambda),  - 1)$ of the nonlinear equation \eqref{evataucent_-} have a local one-dimensional unstable manifold given as the graph of the function $\tau\mapsto\hat Z_{-}(\tau, \lambda, \epsilon)$ with values in $\Lambda^2(\C^4)$ for $\tau$ near $-1$. The manifold can be spread for $\tau\in[-1,1)$ by applying the flow of  \eqref{evataucent_-}. More details of this construction can be found in \cite[p. 182]{AGJ}.

Similarly, for the system \eq{
\frac{d\hat{U}}{d\zeta}& = \left([A_s(\tau, \lambda, \epsilon)]^{(2)} - \mu_{+}(\lambda) I_{6\times6} \right) \hat{U},\\ 
  \frac{d \tau}{d\zeta} &=  k(1-\tau^2),
}{evataucent_+}
 the set $\tau= + 1$ is invariant  and  zero  is the eigenvalue of  $[A_s^+(\lambda, \epsilon)]^{(2)} - \mu_{+}(\lambda) I_{6\times6}$ with  the smallest real part which has the one-dimensional central manifold that consists of the solutions $\hat{U}$ to 
\eq{\left([A_s ^{+}(\lambda, \epsilon)]^{(2)} - \mu_{+}(\lambda) I \right) \hat{U}=0;}{wc+}
thus, the central manifold  is the span of  the eigenvector  $z_{+}(\lambda)$ of $[A_s ^{+}(\lambda, \epsilon)]^{(2)}$
corresponding to the eigenvalue $\mu_{+}(\lambda)$.
 The equilibrium $(z_{+}(\lambda), + 1)$  have a one-dimensional stable manifold given by $\tau\mapsto\hat Z_{+}(\tau, \lambda, \epsilon)$ for the equation \eqref{evataucent_+}.   
 
We define the Evans function $E(\lambda)$ as  
 \eq{ E(\lambda)= e^{-\int_0^\zeta {\rm Trace } \left(A_s(\tau(\zeta), \lambda, \epsilon)\right)} \left(\left(e^{\mu_{-}\zeta} \hat Z_{-}(\tau(\zeta), \lambda, \epsilon)\right)\wedge\left(e^{\mu_{+}\zeta} \hat Z_{+}(\tau(\zeta), \lambda, \epsilon)\right)\right).}{defev}
The wedge product  here gives a scalar and the expression in the right hand side  is shown in \cite{AGJ} to be  independent of $\zeta$,  and, therefore, can be calculated for a fixed value of $\zeta$.  For example, $\zeta=0$ implies $\tau=0$, and so 
\eqnn{  E(\lambda)= \hat Z_{-}(0, \lambda, \epsilon)\wedge \hat Z_{+}(0, \lambda, \epsilon).}
Detailed deduction and the proof of further  properties of the Evans function are given in \cite[Lemma 4.1]{AGJ}.
We point out that this definition stands not only for  $\lambda$ outside of the essential spectrum, but the formula \eqref {defev} is well defined for any $\lambda$ with   $ \realpart{\lambda} >- \rho$, where  $\rho< \min\{\frac{c^2}{4}-\frac{1-\alpha}{1+\eta}, G\}$.  Therefore, we have the following  Lemma: 

\begin{Lemma}\label{lemma8.1}
Given $\eta,\alpha$ as in Theorem \ref{T:1}, there  exists a constant $\rho=\rho(\eta, \alpha)>0$, independent of $\epsilon$ and $\delta$ and an analytic function  $\tilde E(\lambda)$ defined for all $\lambda$  with $\mathrm{Re}\lambda > -\rho$  such that $\tilde E(\lambda)$ is an analytic continuation of $E(\lambda)$. \end{Lemma}

 \subsection{Definition of the Augmented Unstable Bundle} \label{AUB} We fix positive $\epsilon$ and $\delta$, and also $\lambda$ so that $\Re\lambda>-\rho$. Following \cite{AGJ}, we will now define the fibers of the augmented unstable bundle. We stress that they depend on 
 $\epsilon$, $\delta$ and $\lambda$ but we will suppress this dependance in the notations as much as possible until further notice.
 
 The equilibrium $({\bf 0},-1)$ 
 of the system \eqref{evatau} has an  unstable manifold in 
 $\C^4\times[-1,1]$ which has two complex dimensions  plus one real dimension in the direction of $\tau$.  We denote by $U_-(\lambda)$ the (unstable, that is, corresponding to the eigenvalues of $A^-_f(\lambda,\epsilon)$ with positive real parts) linear subspace obtained by intersecting the unstable manifold $W^u$ of the 
equilibrium $({\bf 0},-1)$ of the system \eqref{evatau} with the set $\{\tau=-1\}$. Note that from here on, when writing ``$\{\tau=\tau_0\}$'', we mean the subset of  $\C^4\times[-1,1]$ whose elements are the ones whose last component is $\tau_0$.
A solution of equation \eqref{evslowzeta} approaches ${\bf 0}$ as $\zeta \to -\infty$ if and only if  $(\vec{U}(\zeta),\tau(\zeta)) \to ({\bf 0},-1)$ as $\zeta\to -\infty$. The equation \eqref{evslowzeta} (or \eqref{evatau}) is linear in $\vec{U}$. 
Therefore, for every fixed $\tau_0\in[-1,1)$ the intersection of the global unstable manifold $W^u=W^u(\lambda,\epsilon,\delta)$ of $({\bf 0},-1)$  of  equation \eqref{evatau} with the set $\{\tau=\tau_0\}$ is a two-dimensional linear subspace of $\mathbb C^4$, cf.\ \cite[Lemma 3.3]{AGJ}. 
We denote this subspace by
\begin{equation*}
\Phi_-(\lambda, \tau_0)= W^u(\lambda) \cap\{ \tau=\tau_0\}.\end{equation*}
When we vary $\tau_0$, the linear subspaces $\Phi_-(\lambda, \tau_0)$ form a bundle. The definition of the unstable bundle from \cite[Section 4D]{AGJ} adopted to the eigenvalue problem \eqref{evdelta} reads as follows.
\begin{Definition}   A two-dimensional complex subbundle of $\mathbb C^4\times [-1,1)$  
is defined by
\begin{equation*} 
\Phi_-(\lambda) = \cup_{\tau\in [-1,1)} \Phi_-(\lambda, \tau),\end{equation*}  
and is called the {\em augmented unstable bundle}. 
\end{Definition}
We let $G^2(\C^4)$ denote the Grassmanian which is a complex manifold consisting of two-dimensional subspaces in $\C^4$. The unstable bundle $ \Phi_-(\lambda)$ can be also viewed  as a map 
\begin{equation}\tilde \Phi_-(\lambda,\cdot): [-1,1) \to G^2(\mathbb C^4)
\end{equation}
which maps  each value of $\tau$ to  $\Phi_-(\lambda, \tau)$, the latter being a subspace of $\mathbb C^4$. Following the notation in  \cite[Proposition C1]{AGJ},   we denote  $\Phi_-(\lambda, \tau)$  by $\tilde \Phi_-(\lambda, \tau)$ when it is being considered  as a point in $G^2(\mathbb C^4)$. 

Analogous considerations work at $+\infty$, that is, at $\tau=+1$. Indeed, the equilibrium $(0,1)$ of the system \eqref{evatau} also has a stable manifold, denoted $W^s$, which has two complex dimensions and one real dimension in the direction of $\tau$.  We denote by  $U_+(\lambda)$  the unstable subspace (corresponding to the eigenvalues of $A^+_f(\lambda,\epsilon)$ with positive real parts) at $\tau=+1$ viewed as
a subspace of $\mathbb C^4$ and by $\tilde U_+(\lambda)$ the respective  point in  $G^2(\mathbb C^4)$.  We denote by $\Phi_+(\tau_0,\lambda)$ the subspace in $\C^4$ given by the intersection of $W^s$ and the set $\{\tau=\tau_0\}$ and by $\tilde\Phi_+(\tau_0,\lambda)$ the corresponding point in $G^2(\mathbb C^4)$, cf.\ \cite[Lemma 3.5]{AGJ}.

In \cite[Lemma 3.7]{AGJ} the following fact was established:  If $\lambda$ is not an eigenvalue of \eqref{evdelta}, then  $\tilde \Phi_-(\lambda,\tau) \to \tilde U_+(\lambda)$  as $\tau \to 1$ and the convergence is locally uniform in $\lambda$ over the complement of the spectrum in any open, simply connected and connected  region where the consistent splitting holds. 

This fact is associated to the existence of the exponential dichotomies for \eqref{evfastxi}. Indeed, \eqref{evfastxi} is an asymptotically autonomous equation, and $\lambda$ is not an eigenvalue  if and only if \eqref{evfastxi} has an exponential dichotomy on all of $\mathbb R$. By general results on dichotomies, cf.\ \cite{Coppel}, the dichotomy projections converge as $\xi \to \pm \infty$ to the respective spectral projections of  the matrices $A_f(\pm \infty, \lambda, \epsilon)$. In the notations just introduced that means that the projections onto the subspace $\Phi_+(\lambda,\tau)$, respectively, $\Phi_-(\lambda,\tau)$ converge to the projections onto the subspace  $\C^4 \ominus U_+(\lambda)$, respectively $U_-(\lambda)$ as $\tau \to 1$, respectively $\tau \to - 1$.

Let $K$ be a simple closed curve in the region of the complex plane  of the consistent splitting that  does not contain any of the eigenvalues of \eqref{evdelta}.  Let $K_0$ be the closed region bounded by $K$. For each $\lambda\in K_0$ we consider $2$-dimensional   subspaces $U_{\pm}(\lambda)$. These form  $2$-dimensional  bundles over $K_0$. 
In the framework of   $G^2(\mathbb C^4)$, the correspondence between $\lambda$ and the point  $\tilde U_{\pm}(\lambda)$  is continuous. 
Now consider a cylinder,
\begin{equation}\label{defcK}
\cK:=\big(K\times(-1,1)\big)\cup K^0\times\{-1\}\cup K^0\times\{+1\},
\end{equation}
given by $K\times(-1,1)$  capped with $K^0\times\{\pm1\}$. The arguments in \cite[p.180]{AGJ} show that the following maps  are well defined and continuous when $(\lambda,\tau)\in\cK$,
\begin{equation}\label{defG}
{{Q}}(\lambda,\tau)=\begin{cases}\tilde U_-(\lambda), \,\,\,\text{if }\,\,\, \lambda\in K^0, \,\,\, \tau=-1,\\
\tilde \Phi_-(\lambda,\tau), \,\,\,\text{if } \,\,\, \lambda\in K, \,\,\, \tau\in (-1,1),\\
\tilde U_+(\lambda), \,\,\,\text{if } \,\,\, \lambda\in K^0, \,\,\, \tau=1.
\end{cases}
\end{equation}
While this definition is given in terms of objects from $G^2(\mathbb C^4)$, one can translate it to the linear subspaces  of $\mathbb C^4$ by using the pullback of the canonical bundle $\Gamma_2(\mathbb C^4)$ over $G^2(\mathbb C^4)$ by the map ${{Q}}$.  Recall that the canonical bundle $\Gamma_2(\mathbb C^4)$ is given by the identity map over $G^2(\mathbb C^4)$.  This means that we distinguish between
 the linear subspaces in $\mathbb C^4$, and the points in the Grassmanian $G^2(\mathbb C^4)$, and use the identical map $i$ to attach to each point of  $G^2(\mathbb C^4)$ the respective linear subspace in  $\mathbb C^4$. 
The map ${{Q}}$ defined on the cylinder as in \eqref{defG} takes values in the Grassmanian $G^2(\mathbb C^4)$. Thus the composition 
$i\circ G$ is the pullback by ${{Q}}$ of the identical map $i$, that is, of the canonical bundle $\Gamma_2(\mathbb C^4)$ with the base $G^2(\mathbb C^4)$. The map
$i\circ {{Q}}$
is a map defined on the cylinder $\cK$,  and is taking values in $\Gamma_2(\mathbb C^4)$.
Following \cite[p.180]{AGJ}, we introduce the following notation.

\begin{Definition}\label{defEK} The augmented unstable bundle    $\mathcal E(K)$  for the eigenvalue problem \eqref{evdelta}   is the pullback of the canonical bundle $\Gamma_2(\mathbb C^4)$ by the map ${{Q}}$. 
\end{Definition}
Thus, the fibers of  $\mathcal E(K)$ are the two dimensional subspaces in $\C^4$ corresponding to the points $G(\lambda,\tau)\in G^2(\C^4)$ in \eqref{defG} while $\cK$ is the base of $\mathcal E(K)$. As we have mentioned earlier, all the objects in \eqref{defG} do depend on $\epsilon$ and $\delta$. To stress this fact, in what follows we will use notation $\mathcal E(K)[\eps,\delta]$, $U_\pm(\lambda)[\eps,\delta]$, $\Phi_-(\lambda,\tau)[\eps,\delta]$ when needed.


The stability index is then defined as the first Chern number $c_1$ of the bundle. The first Chern number is a topological invariant of the bundle that measures in some sense the non-triviality of the bundle. Instead of giving a formal definition we  refer to the  property in  \cite[Theorem, p. 173]{AGJ}  that for the augmented unstable bundle  $\mathcal E(K)$  the first Chern number  is equal to the number of zeroes of the Evans function over the contour $K$ and thus gives the count of eigenvalues inside of the contour $K$. 
In some situations, the bundle can be represented as  the Whitney sum of some subbundles.  The Whitney sum of bundles defined on a base space is the vector bundle over the same base space   with fibers given by the direct sum of the fibers of the bundles participating in the Whitney sum.  The important feature of the first Chern number is that it is additive in  the Whitney sum. 

When a system with a slow-fast structure is considered, then this structure is imprinted into the eigenvalue problem as well.  The phase space for the eigenvalue problem written as a dynamical system consists of the  manifold with the slow dynamics and the complementary regions  exhibiting fast dynamics. The augmented unstable bundle can naturally  be  decomposed \cite{AGJ, GJ90, kuehn} into the Whitney sum of the associated slow subbundle $\mathcal E_s(K)$  and fast subbundle $\mathcal E_f(K)$: 
\begin{equation}\mathcal E(K)=\mathcal E_s(K)\oplus \mathcal E_f(K).\label{154}\end{equation} 
Using the additivity of the first Chern number, then 
\begin{equation}c_1(\mathcal E(K))=c_1(\mathcal E_s(K))+  c_1(\mathcal E_f(K)).\label{155}\end{equation}

\subsection{Slow-fast structure of the Augmented Unstable Bundle}

In what follows, $0<\epsilon \ll \delta \ll 1$   is  assumed  to be such that a traveling front $(u_f,w_f)$ exists as stated by Theorem \ref{T:1}. We remind the reader that the limits of $(u_f,w_f)$   as  $\epsilon \to 0$ and $\delta \to 0$ do exist (see Section~\ref{RM}). 

{To make this information transparent and following the notations introduced in Section~\ref{SFEP},} we denote 
$(u_f,w_f)=(u_f{[\epsilon,\delta]},w_f{[\epsilon,\delta]})$ for $\eps,\delta>0$, and then set
\eq{ (u_f{[0,\delta]}, w_f{[0,\delta]})&=\lim_{\epsilon \to 0}(u_f{[\epsilon,\delta]}, w_f{[\epsilon,\delta]}),\\ (u_f{[0,0]}, w_f{[0,0]})&=\lim_{\delta \to 0}(u_f{[0,\delta]}, w_f{[0,\delta]}).
}{limdef}
We will use these notations regardless of what the  independent variable is, assuming that the independent variable may be easily identified  from the context.

 We consider the compactified version \eqref{slowdtau} of the eigenvalue problem written as a first-order system, and we multiply the 2nd equation  by $\epsilon$ to obtain
 \eq{
  \frac{dU_1}{d\zeta}&=U_2,\\
 \epsilon \frac{dU_2}{d\zeta}& = \left((\lambda - \frac{1}{\delta } f_u(u_f[\epsilon,\delta],w_f[\epsilon,\delta]) )U_1 -cU_{2}-  \frac{1}{\delta }f_w(u_f[\epsilon,\delta],w_f[\epsilon,\delta])W_1\right),\\
 \frac{dW_1}{d\zeta}&=W_2,\\
 \frac{dW_2}{d\zeta}&= -g_u(u_f[\epsilon,\delta],w_f[\epsilon,\delta])U_1+(\lambda - g_w(u_f[\epsilon,\delta],w_f[\epsilon,\delta]))W_1 -cW_{2} ,\\
  \frac{d \tau}{d\zeta} &=k(1-\tau^2).
}{slowdtau2}
{ We  here understand that implicitly this is an autonomous system, as these equations are a  subset of a system where the eigenvalue problem is coupled to the traveling wave equation described in Section~\ref{SFEP}. We will not provide details about the reduction of the traveling wave equation which is decoupled from \eqref{slowdtau2}  since it is done in detail in \cite{CGM}.} 
We  consider it along  with the equation \eqref{fastdtau},  which is \eqref{slowdtau2}  written in the fast variable $\xi=\zeta/\epsilon$ 
 \eq{
  \frac{dU_1}{d\xi}&=\epsilon U_2,\\
 \frac{dU_2}{d\xi}& = (\lambda - \frac{1}{\delta } f_u(u_f[\epsilon,\delta],w_f[\epsilon,\delta]) )U_1-cU_{2} -  \frac{1}{\delta }f_w(u_f[\epsilon,\delta],w_f[\epsilon,\delta])W_1,\\
 \frac{dW_1}{d\xi}&=\epsilon W_2,\\
 \frac{dW_2}{d\xi}&=  \epsilon \left(-g_u(u_f[\epsilon,\delta],w_f[\epsilon,\delta])U_1+(\lambda - g_w(u_f[\epsilon,\delta],w_f[\epsilon,\delta]))W_1-cW_{2} \right),\\
 \frac{d \tau}{d\xi}  &= \epsilon k(1-\tau^2).
}{fastdtau2} 
 In the limit $\epsilon\to 0$,  the dynamics of system \eqref{slowdtau2}  can be described as the ``slow''  flow  on the set 
  \eq{ M_0= \Big\{(U_1,U_2,W_1,W_2, \tau):\, U_{2}= \frac{1}{c}\left( (\lambda - \frac{1}{\delta } f_u(u_f[0,\delta],w_f[0,\delta]) )U_1 -  \frac{1}{\delta }f_w(u_f[0,\delta],w_f[0,\delta])W_1\right) \Big\}
 }{m}  
  given by 
\eq{
  \frac{dU_1}{d\zeta}&=\frac{1}{c}\left((\lambda - \frac{1}{\delta } f_u(u_f[0,\delta]w_f[0,\delta]) )U_1 -  \frac{1}{\delta }f_w(u_f[0,\delta],w_f[0,\delta])W_1\right) ,\\
 \frac{dW_1}{d\zeta}&=W_2,\\
 \frac{dW_2}{d\zeta}&= -g_u(u_f[0,\delta],w_f[0,\delta])U_1+(\lambda - g_w(u_f[0,\delta],w_f[0,\delta]))W_1 -cW_{2} ,\\
  \frac{d \tau}{d\zeta} &=k(1-\tau^2).
}{fastdtau0}
To understand the dynamics in the direction, transversal to the set \eqref{m},   we look at the limit of  \eqref{fastdtau2} as $\epsilon \to 0$ and obtain the  ``fast flow'' 
\eq{
  \frac{dU_1}{d\xi}&=0,\\
 \frac{dU_2}{d\xi}& = (\lambda - \frac{1}{\delta } f_u(u_f[0,\delta],w_f[0,\delta]) )U_1 -cU_{2}-  \frac{1}{\delta }f_w(u_f[0,\delta],w_f[0,\delta])W_1,\\
 \frac{dW_1}{d\xi}&=0,\\
 \frac{dW_2}{d\xi}&=  0,\\
 \frac{d \tau}{d\xi}  &= 0.
}{fastdtau20} 

The set $M_0$ given in  \eqref{m}  is  a set of equilibria for the system \eqref{fastdtau20}. This set is  normally hyperbolic and attracting.
Indeed,  all points of $M_0$ are fixed points for the flow of \eqref{fastdtau20} while each point 
outside of $M_0$ generates a trajectory of the flow \eqref{fastdtau20} that is being exponentially attracted to $M_0$ as the second equation in \eqref{fastdtau20} has the term $-cU_2$.  {{ Because of the normal hyperbolicity of $M_0$, by Fenichel theory, in the linear system \eqref{fastdtau2} which is  a perturbation of the limiting system  \eqref{fastdtau20}}}, the flow of \eqref{fastdtau2} for $\epsilon>0$ has an invariant manifold $M_\epsilon=M_0+O(\epsilon)$ as $\epsilon\to0$. Thus, the flow of the equivalent system  \eqref{slowdtau2} for $\epsilon>0$ also
has an invariant manifold $M_\epsilon=M_0+O(\epsilon)$. The set $M_\epsilon$ is also attracting. 
 The flow of \eqref{slowdtau2} on $M_\epsilon$ 
is an $\epsilon$-perturbation of the flow of \eqref{fastdtau0} on $M_0$, more precisely,
\eq{
  \frac{dU_1}{d\zeta}&=\frac{1}{c}\left((\lambda - \frac{1}{\delta } f_u(u_f[0,\delta],w_f[0,\delta]) )U_1 -  \frac{1}{\delta }f_w(u_f[0,\delta],w_f[0,\delta])W_1\right) +O(\epsilon) ,\\
 \frac{dW_1}{d\zeta}&=W_2,\\
 \frac{dW_2}{d\zeta}&= -g_u(u_f[0,\delta],w_f[0,\delta])U_1+(\lambda - g_w(u_f[0,\delta],w_f[0,\delta]))W_1 -cW_{2} +O(\epsilon),\\
  \frac{d \tau}{d\zeta} &=k(1-\tau^2).
}{fastdtau0pert}
{Here  and in $M_{\epsilon}$, the  $O(\epsilon)$ terms include  the small deformation of the front $$(u_f[\epsilon,\delta],w_f[\epsilon,\delta])= (u_f[0,\delta] +O(\epsilon), w_f[0,\delta] +O(\epsilon)).$$ This explains the appearance of  the $O(\epsilon)$ term appear in the 3rd equation. 
}

 Our conclusion is based on  a set of following lemmas.

\begin{Lemma} \label{L19} For every fixed $\delta$, there exists $\epsilon_0$ such that for any $\epsilon <\epsilon_0$ 
\begin{equation}\mathcal E(K)[\epsilon,\delta]=\mathcal E_{slow[\epsilon]}(K)[\epsilon,\delta]\label{156}\end{equation} 
or, equivalently $\mathcal E_{fast[\epsilon]}(K)[\epsilon,\delta]=\emptyset$. 
Then \begin{equation} c_1 (\mathcal E(K)[\epsilon,\delta]) = c_1(\mathcal E(K)[0,\delta]).
\end{equation}
\end{Lemma}
\begin{Proof}
We focus on the fast subbundle in  the decomposition \eqref{155}. The fast subbundle is build as a perturbation  of an unstable bundle in the limiting system \eqref{fastdtau20}. It is easy to see that   the equilibria of \eqref{fastdtau20} do not have unstable manifolds. More precisely, the subspace $U_1=W_1=W_2=0$ is an invariant, "fast" subspace where  the dynamics given by $U_2= -c U_2$, for each $\tau\in [-1,1]$.  Since $c>0$, there are no unstable eigenvectors, therefore the unstable bundle for this system is an empty set.  Another, more direct way to see it is by considering \eqref{fastdtau2}. Indeed, let's  denote  the system  \eqref{fastdtau2} as
\eq{\frac{d \vec{U}}{d\xi}  &=  A_f(\tau, \lambda, \epsilon) \vec{U}, \\ \frac{d \tau}{d\xi}  &= \epsilon k(1-\tau^2).}{u}
The matrix $A_f(\tau, \lambda, \epsilon)$ has  three slow and  one fast variable, moreover
$$A_f^+(\lambda, \epsilon)=\lim_{\tau\to\infty} A_f(\tau, \lambda, \epsilon)$$
 has  three eigenvalues of order $\epsilon$ and one eigenvalue of order $O(1)$. The eigenvalue of order one is strictly negative. 
  Therefore in the decomposition \eqref{155}, the fast unstable subbundle is empty and, thus, \eqref{156} holds. 
  
 Since the slow dynamics is reduced to the set $M_{\epsilon}$ and is given by   \eqref{fastdtau0pert}, the slow subbundle $\mathcal E_{slow[\epsilon]}(K)[\epsilon,\delta]$  is constructed as an $\epsilon$-order perturbation of the unstable bundle in the limiting system  \eqref{fastdtau0} and has the same first Chern number. What follows is an explanation of this fact adopted from \cite{GJPhase}.
 
  Let's look at  \eqref{fastdtau20}  again for a fixed value of $\lambda$.  System  \eqref{fastdtau20} has three zero eigenvalues.  The eigenvectors corresponding to these zero eigenvalues are 
  \eq{
  (1, \frac{1}{c} (\lambda - \frac{1}{\delta } f_u(u_f[0,\delta],w_f[0,\delta]) ), 0, 0),\,\,\, (0,   \frac{1}{c\delta }f_w(u_f[0,\delta],w_f[0,\delta]), 1 , 0, 0),\,\,\, 
(0,0,0,1).}{}
 For every fixed $\tau$  the subspaces spanned by these vectors are points  in $G^3(\mathbb C^4)$, so they form a curve of critical points parametrized by $\tau$. Since there is only one nonzero eigenvalue in the system \eqref{fastdtau20}, all of these points are normally hyperbolic and the curve perturbs to a slow manifold for sufficiently small $\epsilon$. By allowing  $\lambda$  to vary inside of the region $K^0$, we obtain a normally hyperbolic subset  of  $K^0\times (-1,1)\times G^3(\mathbb C^4)$.  For sufficiently small $\epsilon$, it perturbs to a 3-dimensional invariant set with respect to the flow of \eqref{fastdtau2}. This defines a 3 dimensional  subbundle of $K^0\times (-1,1)\times  \mathbb C^4$. Let's denote this subbundle $\Gamma_s(\epsilon)$.
  We now consider a projection of $\mathbb C^4$ onto the slow variables $U_1$, $W_1$, and $W_2$ (onto $\mathbb C^3$). When $\epsilon=0$, the slow flow on    $K^0\times (-1,1)\times \mathbb C^3$ is given by \eqref{fastdtau0}  and for sufficiently small  $\epsilon>0$, the flow is given by \eqref{fastdtau0pert}. 
  
The unstable subspaces in   $\Gamma_s(\epsilon) \times \{\tau=\pm 1\}$ are 2 -dimensional.   We can define the augmented unstable slow bundle $\mathcal E_{slow[\epsilon]}(K)[\epsilon,\delta]$  using the same map as in \eqref{defG}. Indeed,  the projectivised flow on $\mathbb C^3$ is embedded in   $\mathbb C^4$ through the slow variables. Notice that all of the components of  the  map  in the definition \eqref{defG}  depend on $\epsilon$. Moreover, this dependance is continuous.  When $\epsilon=0$, the bundle $\mathcal E_{slow[\epsilon]}(K)[0,\delta]$ is the augmented unstable bundle for the reduced problem \eqref{fastdtau0}. Since the map \eqref{defG} is continuous in $\epsilon$,  the first Chern numbers  of $\mathcal E_{slow[\epsilon]}(K)[0,\delta]$ and $\mathcal E_{slow[\epsilon]}(K)[\epsilon,\delta]$ are the same if $\epsilon$ is  sufficiently small.
 \end{Proof}

Next, we focus on the limiting  flow  generated by \eqref{fastdtau0} which we rewrite as
\eq{
\delta \frac{dU_1}{d\zeta}&=\frac{1}{c}
\left((\delta \lambda -  f_u(u_f[0,\delta],w_f[0,\delta]) )U_1 -  f_w(u_f[0,\delta],w_f[0,\delta])W_1\right) ,\\ \frac{dW_1}{d\zeta}&=W_2,\\ 
\frac{dW_2}{d\zeta}&=  -g_u(u_f[0,\delta],w_f[0,\delta])U_1+(\lambda - g_w(u_f[0,\delta],w_f[0,\delta]))W_1-cW_{2} ,\\ 
\frac{d \tau}{d\zeta} &=k(1-\tau^2). 
}{fastddelta}
We consider it together with the equivalent system in a the scaling $s=\zeta/\delta$, 
\eq{
 \frac{dU_1}{ds}&=\frac{1}{c}\left((\delta \lambda -  f_u(u_f[0,\delta],w_f[0,\delta]) )U_1 -  f_w(u_f[0,\delta],w_f{[0,\delta]})W_1\right) ,\\
 \frac{dW_1}{ds}&=\delta W_2,\\
 \frac{dW_2}{ds}&=  \delta( -g_u(u_f[0,\delta],w_f[0,\delta])U_1+(\lambda - g_w(u_f[0,\delta],w_f[0,\delta]))W_1-cW_{2}),\\
  \frac{d \tau}{ds} &=\delta k(1-\tau^2), 
}{slowdelta}
and exploit the smallness of $\delta$ in  both \eqref{fastddelta}  and \eqref{slowdelta}. More precisely, we prove the following Lemma. 
\begin{Lemma} \label{L20} There exists  $\delta_0$ such that  for any $\delta <\delta_0$ 
\begin{equation}\mathcal E(K)[0, \delta]=\mathcal E_{slow(\delta)} (K)[0,\delta] \oplus \mathcal E_{fast(\delta)}(K)[0,\delta], \label{189}\end{equation} 
where {$\mathcal E_{fast(\delta)}(K)[0,\delta]$} is trivial.
Therefore, 
\begin{equation} c_1 (\mathcal E(K)[0, \delta]) = c_1(\mathcal E_{slow(\delta)} (K)[0,0]).
\end{equation}
Moreover, 
$\mathcal E_{slow(\delta)} (K)[0,0] $ 
 is generated by the eigenvalue problem for the generalized Fisher-KPP equation   \eqref{reduction}  responsible for the existence of the limiting heteroclinic orbit.
\end{Lemma}

\begin{Proof}
When $\delta =0$, the system \eqref{slowdelta} becomes
\eq{
 \frac{dU_1}{ds}&=\frac{1}{c}\left((-  f_u(u_f[0,0],w_f[0,0]) )U_1 -  f_w(u_f[0,\delta],w_f{[0,0]})W_1\right) ,\\
 \frac{dW_1}{ds}&=0,\\
 \frac{dW_2}{ds}&=  0,\\
  \frac{d \tau}{ds} &=0, 
}{slowdelta0}
and it has a set of equilibria
\begin{equation}
\left\{(U_1, W_1,W_2, \tau):\,  U_1= - \frac{ f_w(u_f[0,0],w_f[0,0])}{ f_u(u_f[0,0],w_f[0,0]) }W_1\right\}.\label{M0d}
\end{equation}
The linearization of \eqref{slowdelta0} about each equilibrium point  has exactly one nonzero  eigenvalue 
\begin{equation}-\frac{ f_u(u_f[0,0],w_f[0,0])} {c} >0.\end{equation} 

 The corresponding 1-dimensional  eigenspace  is the span of the eigenvector $(1, 0, 0)$ or $U_1$-axis in $\mathbb C^3$. Let us vary $\lambda$ in $K^0$. A subbundle $K^0\times (-1,1)\times \mathrm{Span}(1,0,0)$ is a trivial bundle invariant for \eqref{slowdelta0}. For a fixed $\lambda$ and sufficiently small $\delta$ ,  the unstable manifold  $\mathrm{Span}(1,0,0)$ perturbs to strong unstable manifolds at  $\tau=\pm 1$   respectively, each of which can be considered  a point in $G^1(\mathbb C^3)$.

For each compact subset of \eqref{M0d} and  sufficiently small $\delta$ there exists an invariant set  for  \eqref{slowdelta} which is a $\delta$-order perturbation of that compact subset of  \eqref{M0d}. However, the construction of the unstable bundle is done over compact spaces such as  the space of the Grassmanian
manifolds 
or the projective spaces.   We see from  \eqref{slowdelta}  which  we rewrite as
\eq{
\frac{d \vec{\mathcal U}}{ds}  &= B(\tau, \lambda, \delta) \vec{\mathcal U},\\
  \frac{d \tau}{d\zeta} &=\delta k(1-\tau^2),
}{w}
where  $\vec{\mathcal U}^T =(U_1,W_1,W_2)$ and 
\begin{equation}
B= \begin{pmatrix}
\frac{1}{c}\left( \delta \lambda - f_u(u_f[0,\delta],w_f[0,\delta]) \right)& - \frac{1}{c}\ f_w(u_f[0,\delta],w_f[0,\delta])&0 \\
0&0&\delta\\
 -\delta g_u(u_f[0,\delta],w_f[0,\delta]) & \delta (\lambda - g_w(u_f[0,\delta],w_f[0,\delta]))& -c\delta
\end{pmatrix}
\end{equation}
that $U_1$ is a fast variable  and $W_1$ and  $W_2$ 
are the slow variables.  We will consider the projectivized flow generated by the system  \eqref{slowdelta}  on $G^1(\mathbb C^3)= \mathbb CP^2$, 
\eq{\frac{d \hat{\mathcal U}}{d\zeta}  &= B^{(1)}(\tau, \lambda, \delta; \hat{\mathcal  U}), \quad    \hat{\mathcal U} \in G^1(\mathbb C^3), \\ \frac{d \tau}{ds} &=\delta k(1-\tau^2).
}{ww}
When $\delta=0$, $\tau$ is treated as a parameter.  For each fixed $\tau \in [-1,1]$, matrix $B(\tau, \lambda, 0)$ has one nonzero, positive eigenvalue. The one-dimensional unstable eigenspace is a critical point for projectivized flow, moreover it is attracting for any $\tau$. As an attractor, it perturbs to a nearby attractor for \eqref{ww} for sufficiently small $\delta$. No eigenfunctions may be formed by following that perturbed, fast unstable manifold;  it is used to construct the fast  unstable bundle $\mathcal E_{fast(\delta)}(K)[0,\delta]$   
for  all $\lambda$ in $K$, and, thus, the fast unstable bundle is trivial. 

The complementary slow subbundle $\mathcal E_{slow(\delta)} (K)[0,\delta]$  is also constructed as a perturbation of a bundle associated with the limiting flow which we denote  $\mathcal E_{slow(\delta)} (K)[0,0]$. 

The system \eqref{slowdelta0} has two zero eigenvalues for each fixed $\tau$.  The eigenspace associated with  these eigenvalues is spanned by $(0,0,1)$ and $\left( - \frac{ f_w(u_f[0,0],w_f[0,0])}{ f_u(u_f[0,0],w_f[0,0])}, 1, 0\right)$. In $G^2(\mathbb C ^3)$ these form a curve parametrized by $\tau$.  Each point on the curve is a normally hyperbolic critical point of the flow generated by \eqref{slowdelta0}  in $G^2(\mathbb C ^3)$. We can append the equation $\frac{d \lambda}{d \zeta}$  and consider the flow on $K^0\times(-1,1)\times G^2(\mathbb C ^3)$.  We call the normally hyperbolic set of these critical points $\Psi_{slow (\delta)}(0)$. For sufficiently small $\delta$, $\Psi_{slow (\delta)}(0)$ perturbs to $\Psi_{slow (\delta)}(\delta)$. The latter 
  is associated with  a 2-dimensional  invariant  for  \eqref{slowdelta} subbundle in $K^0\times(-1,1)\times \mathbb C ^3$. We call it 
  $\Gamma _{slow(\delta)} (\delta)$. 
  We then project $\mathbb C^3$  onto the slow variables $(W_1,W_2)$, i. e. onto $\mathbb C^2$. 
More precisely,  the flow generated by the system \eqref{fastddelta} with $\delta =0$ is
\eq{
0&=\frac{1}{c}
\left(( -  f_u(u_f[0,0],w_f[0,0]) )U_1 -  f_w(u_f[0,0],w_f[0,0])W_1\right) ,\\ \frac{dW_1}{d\zeta}&=W_2,\\ 
\frac{dW_2}{d\zeta}&=  -g_u(u_f[0,0,w_f[0,0)U_1+(\lambda - g_w(u_f[0,0],w_f[0,0]))W_1-cW_{2} ,\\ 
\frac{d \tau}{d\zeta} &=k(1-\tau^2). 
}{fastddelta0}
  In this context the set \eqref{M0d} is called the slow manifold. 
  The   flow on the slow manifold  is given by 
\eq{
  \frac{dW_1}{d\zeta}&=W_2,\\
 \frac{dW_2}{d\zeta}&=  \left( (g_u(u_f[0,0],w_f[0,0]) \frac{ f_w(u_f[0,0],w_f[0,0])}{ f_u(u_f[0,0],w_f[0,0]) }+(\lambda - g_w(u_f[0,0],w_f[0,0]))\right)W_1-cW_{2} ,\\
  \frac{d \tau}{d\zeta} &=k(1-\tau^2).
}{limdelta0}
It can be considered a flow on $K^0\times (-1,1)\times \mathbb C^2$.  
For sufficiently small $\delta$, the flow  on the invariant set  \eqref{fastddelta}  is a $\delta$-order perturbation of \eqref{limdelta0}, 
\eq{
  \frac{dW_1}{d\zeta}&=W_2,\\
 \frac{dW_2}{d\zeta}&=  \left( (g_u(u_f[0,0],w_f[0,0]) \frac{ f_w(u_f[0,0],w_f[0,0])}{ f_u(u_f[0,0],w_f[0,0]) }+(\lambda - g_w(u_f[0,0],w_f[0,0]))\right)W_1-cW_{2} +O(\delta),\\
  \frac{d \tau}{d\zeta} &=k(1-\tau^2).
}{limdelta0d}
 For each fixed $\tau\in[-1,1]$,  the set \eqref{M0d}  is a 2-dimensional subset of the equilibria of \eqref{slowdelta0}. 

 As long as $c>2 \sqrt{\frac{1-\alpha}{1+\eta}}$, the  unstable manifold of the equilibrium $(W_1, W_2,)=(0,0)$   has one-dimensional unstable manifold at $\tau=\pm 1$. Therefore, the usual capping procedure may be used to define $\mathcal E_{slow(\delta)}(K)[0,\delta]$ within  $\Gamma _{slow(\delta)} (\delta)$.  When $\delta=0$, we call this bundle $\mathcal E_{slow(\delta)}(K)[0,0]$; it measures the number of zeros of the Evans function of the reduced problem. Since for sufficiently small $\delta$ the dependence on $\delta$   in \eqref{limdelta0d} is continuous, the first Chern numbers of  $\mathcal E_{slow(\delta)}(K)[0,0]$ and $\mathcal E_{slow(\delta)}(K)[0,\delta]$  are the same. 
 
We notice that the first two equations in \eqref{limdelta0} are equivalent to the equation of the second order
\begin{equation}\label{e:144}
\lambda W_1=\frac{d^2W_1}{d\zeta^2}+c\frac{dW_1}{d\zeta} + \left(  - g_u(u_f[0,0],w_f[0,0]) \frac{ f_w(u_f[0,0],w_f[0,0])}{ f_u(u_f[0,0],w_f[0,0]) }+ g_w(u_f[0,0],w_f[0,0])\right)W_1.
\end{equation}
Let us consider the reduced equation for the $w$ component of the front \eqref{reduction}. The equation \eqref{reduction} is  the equation
\begin{equation}\label{14}
\frac{d^2  w}{d\zeta^2} +c \frac{dw}{d\zeta} +g(u(w),w)=0,
\end{equation} 
where $u$ is the nonzero solution  of the equation $f(u,w)=0$.  The latter equation leads to $df= f_udu+f_wdw=0$ along the front $(u_f[0,0], w_f[0,0])$. Therefore $du/dw=-f_w/f_u$.
The linearization of   the term $g(u(w),w)$ from  \eqref{14} is then
\begin{equation}\left(- g_u(u_f[0,0],w_f[0,0]) \frac{ f_w(u_f[0,0],w_f[0,0])}{ f_u(u_f[0,0],w_f[0,0]) }+ g_w(u_f[0,0],w_f[0,0]) \right)W_1.\end{equation}
Therefore the equation  \eqref{e:144} is the linearization of the Fisher-KPP equation 
\begin{equation}\label{1444}
w_t= \frac{d^2  w}{d\zeta^2} +c \frac{dw}{d\zeta} +g(u(w),w)
\end{equation} 
about  the front solution characterized by the speed $c$. 
\end{Proof}

{{The general implications of Lemma~\ref{L20} are  as follows. Assume that the eigenvalue problem for the linearization of the reduced equation \eqref{1444} has an eigenvalue  $\lambda_0$ of multiplicity $m$.  There exists $\epsilon_0>0 $ such that for any $\epsilon<\epsilon_0$, there exists $\delta_0=\delta_0(\epsilon) > 0$, such that for any  $\delta <\delta_0$,  there are exactly $m$  (counting multiplicity) zeros of the Evans function associated  
with  \eqref{evdelta2}  in a small neighborhood of order $O(\epsilon, \delta)$  of $\lambda_0$.  

{\bf{Conclusion.}} Based on the results from the current and previous sections, we conclude the stability  of the fronts is governed by the stability of the fronts in \eqref{1444}. The stability of fronts in the scalar Fisher-KPP equations is well understood. 
According to \cite{Sattinger76}[Theorem 6.3], the fronts in the classical Fisher-KPP equation are spectrally stable  in the weighted $L^{\infty}$ norm with  the weight $1+ e^{\tilde \sigma \zeta}$, with some  specified constant $\tilde\sigma$. The stability results  described in \cite{Sattinger76}[Theorem 6.3]  for the classical Fisher-KPP equation without loss of generality hold for the generalized Fisher-KPP equations, such as \eqref{1444}. Indeed,  the spectral stability result   \cite{Sattinger76}[Theorem 6.3]  follows from  \cite{Sattinger76}[Theorem 4.1];  the asymptotic stability in the weighted $C^1$ norm  follows  from \cite{Sattinger76}[Theorem 4.3]. Both of these theorems are applicable to the generalized Fisher-KPP equations. This concludes the proof of Theorems~\ref{T:2} and  implies  the asymptotic stability of the  fronts in the weighted $C^1$-norm with the same weight in which the spectral stability is achieved. The weight removes the eigenvalue caused by the translational symmetry  at the origin, thus the stability is asymptotic and is not orbital.   
}}

\section{Numerics\label{S:Numerics}}

{{In this section, we perform numerical computations of the Evans function to explicitly identify regions in the parameter space where the fronts are spectrally stable in the appropriate weighted space for both case $\epsilon$ zero and nonzero.}} In effect, we are {{numerically providing evidence}} that the front solutions for those parameter values are spectrally stable 
 when considered in $L^2(\R)$ equipped with the weight given in \eqref{disW}, with $\sigma$ satisfying  \eqref{KPPE}.   As explained in Section \ref{EssSect}, the weight is introduced so that the essential spectrum is ``stabilized'' by being moved to the left side of the complex plane. As explained in Section \ref{EssSect}, the weight is introduced so that the essential spectrum is ``stabilized'' by being moved to the left side of the complex plane.   

%
%
 
 The general strategy in this section is to first use the bounds  \eqref{reim} ($\epsilon\neq 0$) and \eqref{delta0Ws} ($\epsilon= 0$) to determine a region of the complex plane containing any eigenvalue with positive real part.  Then, second, we numerically perform a winding number computation of the Evans function around that region to identify any zero that would correspond to an eigenvalue for the problem \eqref{evdelta2}. 
 
 The Evans function computations are performed using the MATLAB-based numerical library for Evans function computation called STABLAB \cite{StabLab}. 
 To perform those computations, we use the formulation \eqref{evslowzeta} of the  eigenvalue problem \eqref{evdelta2} written as a four-dimensional, first-order, linear system of the form 
 \eq{
\frac{d \vec{U}}{d\zeta}  &= A_s(\zeta, \lambda,\epsilon)\vec{U}, 
}{linear}
where $A_s$ is the $4\times 4$ matrix given in \eqref{Am}.  
Note that the linear system \eqref{linear} is referred as the ``slow'' system and it differs from from the ``fast'' system \eqref{fastzeta}  obtained  after a change to the fast variable $\xi =\zeta/\epsilon$ was applied.  The asymptotic behavior as $\zeta\rightarrow\pm\infty$ of the solutions to (\ref{linear}) 
is determined by the matrices 
\[
A_s^{\pm}(\lambda,\epsilon)=\lim_{\zeta\rightarrow \pm\infty}A_s(\zeta,\lambda,\epsilon),
\]
which are found by inserting the values $(u_f,w_f)=(1,0)$ or $(u_f,w_f)=(\alpha,1-\alpha^2)$ in $A_s$ and using \eqref{fderivlim}. Alternatively,  one can use \eqref{fsrel2} giving directly an expression for the asymptotic matrices of the slow system.  We find
\eq{
A_s^{+}(\lambda,\epsilon)= \begin{pmatrix}
 0&1  &0&0\\
 \displaystyle{
\frac{1}{\epsilon}\left(\lambda +\frac{1}{\delta }\right) }&\displaystyle{-\frac{c}{\epsilon}} & \displaystyle{\frac{1}{2\epsilon \delta }}&0\\
 0&0& 0&1\\
0 &0&  \lambda+ \displaystyle{\frac{\alpha-1}{1+\eta}}&- c
\end{pmatrix},}{Ap}
\eq{A_s^{-}(\lambda,\epsilon)=
\begin{pmatrix}
 0&1  &0&0\\
\displaystyle{\frac{1}{\epsilon}\left(\lambda + \frac{2\alpha^2}{\delta(1+\alpha) }\right)} &\displaystyle{-\frac{c}{\epsilon}} & \displaystyle{\frac{\alpha}{\epsilon \delta (1+\alpha)} }&0\\
 0&0& 0&1\\
\displaystyle{\frac{\alpha^2-1}{\eta+\alpha}} &0&  \lambda &- c
\end{pmatrix}.
}{Am2}
The eigenvalues of matrix $A_s^{+}$ are such that when $\lambda$ has positive real part and is on the right side of the parabola defined by the second equation of \eqref{curve++}, then two have positive real part and two have negative real parts.
The two eigenvalues with negative real parts are given by
\eq{\mu_{1}^+=\frac{1}{2\epsilon }\left(-c-\sqrt{c^2+4\epsilon\left(\lambda+\frac{1}{\delta}\right)}\right),\quad  \mu_{2}^+=\frac{1}{2}\left(-c-\sqrt{c^2+4\left(\lambda+\frac{\alpha-1}{1+\eta}\right)}\right).}{eigplus}
Note that those eigenvalues differ from the ones given in \eqref{APP} associated to the fast system by an $\epsilon$ factor. We abuse notation and use the same symbols for the ones above as the ones for the fast system given in \eqref{APP} that do not appear in this section.  

The eigenvalues of $A_s^{-}$ cannot be obtained explicitly. However, we compute the characteristic polynomial and obtain
\eq{
&{r}^{4}+{\frac {c \left( \epsilon+1 \right) {r}^{3}}{\epsilon}}+{
\frac { \left(  \left(  \left( -\epsilon-1 \right) \lambda+{c}^{2}
 \right)  \left( 1+\alpha \right) \delta-2\,{\alpha}^{2} \right) {r}^{
2}}{\delta\, \left( 1+\alpha \right) \epsilon}}-2\,{\frac {c \left( 
\lambda\,\delta\, \left( 1+\alpha \right) +{\alpha}^{2} \right) r}{
\delta\, \left( 1+\alpha \right) \epsilon}}\\
 &+{\frac { \left( 2\,\lambda
-1 \right) {\alpha}^{3}+\left( {\lambda}^{2}\delta+2\,\eta\,\lambda
 \right) {\alpha}^{2}+ \left( 1+\delta\, \left( \eta+1 \right) {
\lambda}^{2} \right) \alpha+\delta\,\eta\,{\lambda}^{2}}{ \left( \eta+
\alpha \right) \epsilon\,\delta\, \left( 1+\alpha \right) }}=0.}{charAm}
For example, in the case  $\alpha=0.75$, $\delta=0.1$, $\epsilon=0.01$, $\eta=3$, $c=1$ and $\lambda=3$, we have the following 
three eigenvalues
$$
r=1.319833682, 8.670594954, -2.314599976, -108.6758287.
$$
To study the solutions of \eqref{charAm} for large real values of $\lambda$, we scale $r$ as $r=\sqrt{\lambda}\rho$.
We substitute in \eqref{charAm}, divide by $\lambda^2$ and apply the limit $\lambda\to\infty$ to obtain
\eq{
&\epsilon{\rho}^{4}
-{ {  \left( \epsilon+1 \right)  }}{\rho}^{2}
+1
=0.
}{charAmlarge}
Since \eqref{charAmlarge} has two positive and two negative solutions ($\rho=\pm1/\epsilon$ and $\rho=\pm 1$), it implies that 
the same is true for \eqref{charAm} when $\lambda$ is real and large enough. 

Thus equation \eqref{charAm} has four real solutions, two positive and two negative, when $\lambda$ is real, positive, and large enough.
Since the signs of the real parts of the eigenvalues only can change for values of $\lambda$ on the left side of the complex plane, as shown by the 
continuous spectrum obtained in Lemmas \ref{3} and \ref{4}, we conclude that when $\lambda$ has positive real part, $A_s^{-}$ has two eigenvalues with positive real part and two eigenvalues with negative real part. 

We denote the two eigenvalues of 
$A_s^{+}$ with negative real parts as $\mu_{i}^+$ with eigenvectors $v_{i}^+$, $i=1$, $2$, and the two eigenvalues of 
$A_s^{-}$ with positive real parts as $\mu_{i}^-$ with eigenvectors $v_{i}^-$,  $i=1$, $2$. 
As a consequence, for $\lambda$ on the right side of the parabola defined by the second equation of \eqref{curve++} with $\realpart{\lambda}>0$, the system (\ref{linear}) has two linearly independent solutions $U_{1}^+$ and $U_{2}^+$, converging to zero as $\zeta\rightarrow\infty$ 
and two solutions {$U_{1}^-$ and $U_{2}^-$} converging to zero as $\zeta\rightarrow -\infty$, satisfying
$$
\lim_{\zeta\rightarrow \pm \infty} U_{i}^\pm e^{-\mu_{i}^\pm\zeta}=v_{i\pm},\;\;i=1,2.
$$
Clearly, $\lambda_0$ is an eigenvalue for the problem (\ref{evdelta}) if and only if 
 the space of solutions of (\ref{linear}) bounded as $\zeta \rightarrow + \infty$, spanned by $\{U_{1}^+,\;U_{2}^+\}$, and
the space of solutions bounded as $\zeta \rightarrow - \infty$, spanned
by $\{U_{1}^-,\;U_{2}^-\}$, have an intersection of strictly
positive dimension when $\lambda=\lambda_0$.
The most straightforward way to test whether these two spaces of
solutions intersect non-trivially is to calculate the determinant of
the two spanning sets, evaluated at some value of $\zeta$ (usually taken as
 $\zeta =0$).  For the purpose of the numerical computations, this is the function
we call the Evans function
\cite{Evans, Jones, Yanagida, AGJ, Pego, Gardner98, Kapitula98a, Sandstede, Kapitula00, Li00}. Although this definition differs from the definition used in Section \ref{EvansDef}, its zeroes correspond to the eigenvalues of \eqref{evdelta2}, 
it is analytic  to the right of the essential spectrum and it is real for $\lambda$
real. Numerically,   the
two solutions  $U_{1}^+$ and $U_{2}^+$ are obtained by integrating (\ref{linear}) backwards from a sufficiently large positive value of $\zeta$, with initial 
conditions in the $v_{1}^+$ and $v_{2}^+$ directions, respectively. Even though these two vectors are linearly independent, the numerical integration
 will lead to an alignment with the eigendirection corresponding to the
eigenvalue with the smallest real part. One way to circumvent this problem is to compute  the Evans function using the alternative 
definition involving exterior algebra
\cite{AfBr01, AlBr02, Br99, BrDeGo02, Br00, DeGo05, NgR79, skms, EvansWeberGroup}. That framework was used in Section \ref{EvansDef} to define the Evans function on the compactified version of the eigenvalue problem.

In our computations below, we need to extend the definition of the Evans function all the way to the imaginary axis. The argument below is similar to the one 
that lead to the extension of the Evans function given by Lemma \ref{lemma8.1}. Since we are using a different definition of the Evans function than in Section \ref{EvansDef}, we argue in detail how the Evans function given in this section can be extended analytically up to the imaginary axis.  
Remember that for the asymptotic matrix \eqref{Am2} at $\zeta\to\infty$,
in the case where $\lambda$ has positive real part and is on the right side of the parabola defined by the second equation of \eqref{curve++}, the matrix $A_s^{+}(\lambda,\epsilon)$ \eqref{Am2} has two eigenvalues with positive real parts and two with negative real part given in \eqref{eigplus}. This allows us to define the Evans function as it is in the discussion above. This reflects the fact that the Evans function is analytic inside any region that does not intersect with the essential spectrum. When $\lambda$ moves to the left of the boundary defined by \eqref{curve++}, the real part of one of the eigenvalues of   $A_s^{+}(\lambda,\epsilon)$ goes from being positive to being negative. Thus, the dimension of the space of solutions   of \eqref{linear} that converges to zero as $\zeta\to \infty$ goes from 2 to 3.  However, it follows from the expressions of the eigenvalues of $A_s^{+}(\lambda,\epsilon)$ that the two eigenvalues given in \eqref{eigplus} remain the two eigenvalues with the smallest real parts as long as the following condition is satisfied 
\eq{
\realpart{\lambda}\geq \max\left(\frac{1-\alpha}{\eta+1}-\frac{c^2}{4},\;-\frac{1}{\delta}-\frac{c^2}{4\epsilon}\right).
}{conda}
Since the right side of \eqref{conda} is negative by the condition on $c$ given by Theorem \ref{T:1}, the region defined by that inequality includes the whole right side of the complex plane. It thus follows from the {{Gap Lemma}} \cite{Gardner98,Kapitula98a}, or from results from the theory of ODEs \cite{codd},  that the Evans functions as we define it above is analytic for the whole region of the complex plane defined by $\realpart{\lambda}\geq 0$,  including the part on the left side of the curve determined by the second equation of \eqref{curve++}. 

Another point of view is to consider the fact that, with the weight introduced in \eqref{disW}, one, in effect, is computing the Evans function for the linear system
  \eq{
\frac{d \vec{U}}{d\zeta }  &= \tilde{A}(\zeta , \lambda)\vec{U}, 
}{linearT}
where $\tilde{A}$ is such that
\eqnn{
\tilde{A}=\left\{
\begin{array}{cl}
A_s-\sigma I&{\text{ for }}\zeta>Z,\\
A_s&{\text{ for }}\zeta<-Z.
\end{array}
\right.
}
where $A_s$ is given in \eqref{Am}. By Inequality \eqref{KPPE}, $\sigma$ is less than the ``rogue'' eigenvalue, i.e. the eigenvalue that goes from 
having a positive real part to a negative real part when crossing into the part of the essential spectrum on the right side of the complex plane. 
As a consequence, the number of eigenvalues with negative or positive real part for the asymptotic matrix associated to system \eqref{linearT} does not change for values of $\lambda$ such that $\realpart{\lambda}>0$, 
and the corresponding Evans function is automatically analytic on the right side of the complex plane. 

Moreover, since the derivative of the front converges to zero at the rate specified 
by \eqref{myre}, we have that it is not an element of the weighted space. Thus, the Evans function ``loses'' the eigenvalue $\lambda=0$ caused by the translation symmetry. We will thus be able below to perform our winding number computations along curves that go through the origin since the Evans function does not have an a priori zero there. 
 
As indicated earlier, we perform the numerical Evans function computations using the MATLAB-based numerical library for Evans function computation called STABLAB \cite{StabLab}. 
 In STABLAB, the Evans function is computed using the {\sl{polar-coordinate method}}, a method initially proposed by Humpherys and Zumbrun in \cite{Humpherys06}, which represents the unstable and stable manifolds using the continuous orthogonalization method of Drury \cite{Drury80} together with a scalar ODE that restores analyticity. Since we are interested in the zeroes of the Evans function, the standard method is to compute the integral of the logarithmic derivative of the Evans function on a given closed curve and obtain its winding number along that curve.  
In order to numerically verify that there are no zeroes of the Evans function inside a given region of the complex plane,  we choose a closed curve {whose interior encloses} the region.

\begin{figure}[h]
\vspace*{0mm}
\hspace{-0cm}
\scalebox{.5}{{\includegraphics{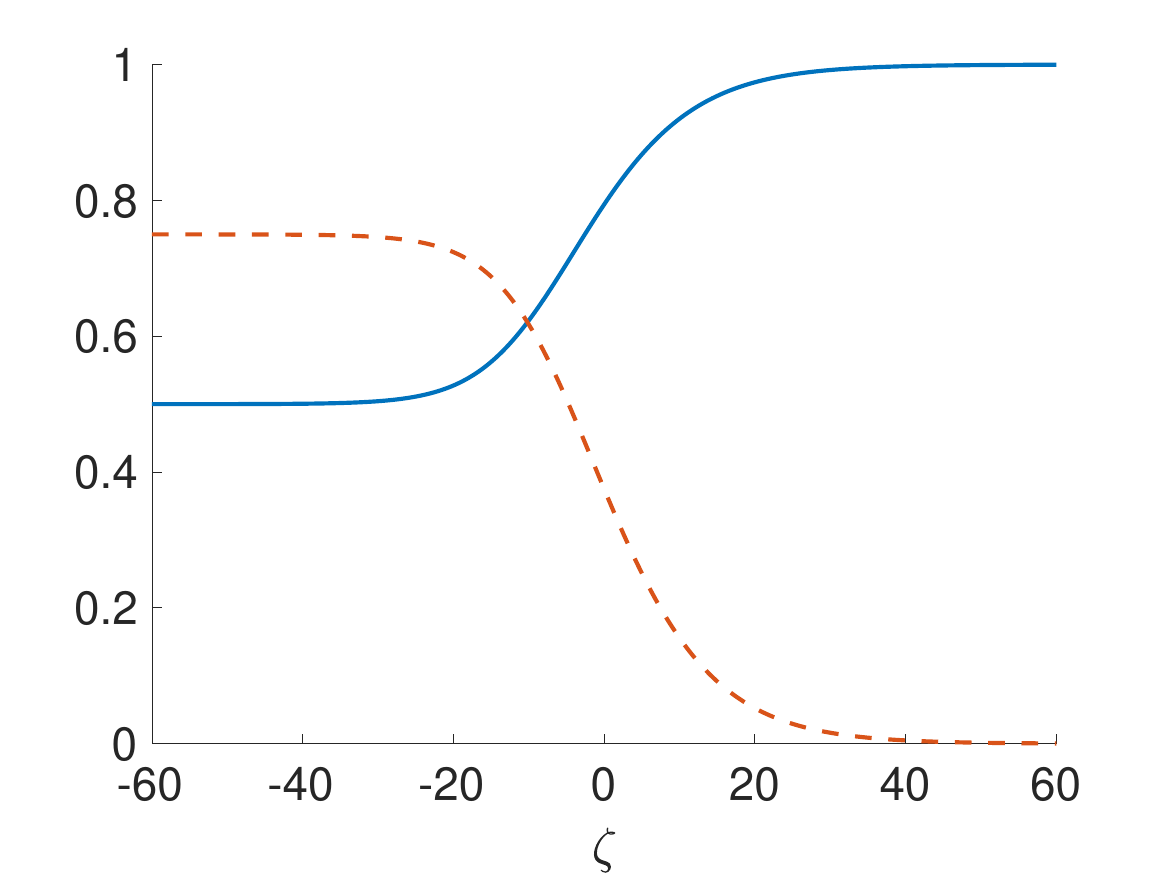}}}
\vspace*{00mm}
\caption{\label{Front1} Front solution described by Theorem \ref{T:1} also given as an heteroclinic of system \ref{e:tw} in the case  
$\alpha=0.5$, $\delta=0.1$, $\epsilon=0.05$, $\eta=2$, and $c=1.5$. The solid curve corresponds to $u$, while the dashed 
curve to $w$.
}
\end{figure}

STABLAB computes the front solution using the MATLAB fifth-order collocation package bvp5c, with relative and absolute errors equal to $10^{-8}$ and $10^{-9}$, respectively. 
In order to provide strong numerical evidence for the claim that there is a region in the parameter space for which the fronts are spectrally stable, we have performed STABLAB computations for 1000 cases for values of the parameters satisfying $2\leq \eta \leq 3, 0.1\leq \alpha \leq 0.8, 1\leq c\leq 2, \delta=0.1$, and $\epsilon=0.05$, all satisfying the conditions listed in Theorem \ref{T:1}. In each case, we numerically computed the bound prescribed by the RHS of \eqref{reim} in Lemma \ref{L:2}. The expression in \eqref{reim}  depends on four positive free parameters $\beta_i,\;i=1,2,3,4$. While, ideally, one would desire to obtain the values of those parameters that give the least bound, we found that task to be rather complicated. Instead, we took the pragmatic position to find values of the $\beta_i$'s that give a bound that is ``good enough''.  To be more precise, when computing the bound, we numerically found the values of the $\beta_i$'s between 0 and 1 that give the least bound.  We used a value slightly superior to that bound as the radius of the semicircle centered at the origin and contained in the right side of the complex plane to define a close curve. By Lemma \ref{L:2}, all the unstable eigenvalues are in the region contained inside that curve.   In each case, we used STABLAB to compute the winding number of the Evans function along that semicircle and found zero every time.  Figure \ref{Front1} shows the graph of the two-component front obtained through STABLAB in the specific case where $\alpha=0.5$, $\delta=0.1$, $\epsilon=0.05$, $\eta=2$, and $c=1.5$.

In the case $\epsilon=0$, we consider the eigenvalue problem \eqref{evdelta2}, which becomes a third-order system when $\epsilon$ is set to zero. We write it as a linear system of the form 
 \eq{
\frac{d \vec{U}}{d\zeta}  = A_0(\zeta, \lambda)\vec{U}, 
}{linear0} 
where
\eq{
A_0(\zeta, \lambda)\equiv \begin{pmatrix}
\displaystyle{\frac{\lambda}{c} - \frac{1}{\delta } f_u(u_f,w_f)} &\displaystyle {-  \frac{1}{\delta }f_w(u_f,w_f)}&0\\
 0&0& 1\\
-  g_u(u_f,w_f) &  \lambda - g_w(u_f,w_f)&- c
\end{pmatrix}.}{Amm3} 
The asymptotic  matrices are found like previously, by using \eqref{fderivlim}. The two matrices are given by 
\eq{
A_0^{+}(\lambda)= \begin{pmatrix}
 \displaystyle{
\frac{\lambda}{c} +\frac{1}{\delta }}& \displaystyle{\frac{1}{2 \delta }}&0\\
 0& 0&1\\
0 &  \lambda+ \displaystyle{\frac{\alpha-1}{1+\eta}}&- c
\end{pmatrix},}{Ap3}
\eq{A_0^{-}(\lambda)=
\begin{pmatrix}
\displaystyle{\frac{\lambda}{c} + \frac{2\alpha^2}{\delta(1+\alpha) }} & \displaystyle{\frac{\alpha}{\delta (1+\alpha)} }&0\\
 0& 0&1\\
\displaystyle{\frac{\alpha^2-1}{\eta+\alpha}} &  \lambda &- c
\end{pmatrix}.
}{Am23}
In the case where $\lambda$ has positive real part and is on the right side of the parabola defined by the second equation of \eqref{curve++}, the matrix $A_0^{+}(\lambda)$ has two eigenvalues with positive real parts and one with negative real part. The eigenvalue with negative real part 
is given by
\eqnn{\mu_{1}^+=\frac{1}{2}\left(-c-\sqrt{c^2+4\epsilon\left(\lambda+\frac{\alpha-1}{1+\eta}\right)}\right).}
Like in the case $\epsilon\neq 0$, the eigenvalues of $A_0^{-}(\lambda)$ cannot be obtained explicitly. However, their signs can be obtained by considering the case where $\lambda$ and $\eta$ are large. One then finds that two of the eigenvalues $\mu_{1}^-$ and $\mu_{2}^-$ have positive real part, while the third one has negative real parts. Thus, system \eqref{linear} with the matrix $A$ given by \eqref{Amm3}, has two (one) linearly independent solutions decaying to zero as $\zeta\rightarrow -\infty$ ($\zeta\rightarrow \infty$).

In order to numerically find a region in the parameter space for which the fronts are spectrally stable, we proceeded the same way as in the $\epsilon\neq 0$ case. That is, we have performed STABLAB computations for 1000 cases for values of the parameters satisfying $2\leq \eta \leq 3, 0.1\leq \alpha \leq 0.8, 1\leq c\leq 1, \delta=0.1$, and $\epsilon=0$, all satisfying the conditions listed in Lemma \ref{L:d} for the existence of the fronts. To define the radius of the semicircle along which we compute the winding number, we used a value slightly greater than the bound provided by Lemma \ref{L:3}. Once more, in each case, we found the winding number of the Evans function to be zero.

\section{Acknowledgements} 

The authors acknowledge hospitality and support of Banff International Research Station   where the plan for this project  was developed during a Research in Teams  workshop organized by the authors  in 2019. 

{The authors would like to thank the anonymous referee for the valuable and thoughtful comments and suggestions that led to the enhancement of the main result.}

Y.L. was supported by the NSF grant DMS-2106157, and would like to
thank the Courant Institute of Mathematical Sciences at NYU and
especially Prof.\ Lai-Sang Young for their hospitality. 

A.G. was in part supported by Miami University through Faculty Research Grants.

\end{document}